\documentclass[11pt]{scrartcl}
\addtokomafont{disposition}{\rmfamily}

\usepackage{amsmath, amssymb, amsthm}
\usepackage{quiver}
\usepackage{tikz}
\usetikzlibrary{shapes,arrows,intersections}
\usepackage{enumitem}
\usepackage{mathtools}
\usepackage{authblk} 
\usepackage{hyperref} 
\usepackage{cleveref} 
\usepackage[nottoc]{tocbibind} 
\usepackage{comment}

\usepackage{todonotes}

\usetikzlibrary{calc}
\usetikzlibrary{trees}

\tikzstyle{morphism}=[fill=white, draw=black, shape=rectangle]
\tikzstyle{medium box}=[fill=white, draw=black, shape=rectangle, minimum width=0.7cm, minimum height=0.7cm]
\tikzstyle{large morphism}=[fill=white, draw=black, shape=rectangle, minimum width=1.7cm, minimum height=1cm]
\tikzstyle{bn}=[fill=black, draw=black, shape=circle, inner sep=1.5pt]
\tikzstyle{state}=[fill=white, draw=black, regular polygon, regular polygon sides=3, minimum width=0.8cm, shape border rotate=180, inner sep=0pt]
\tikzstyle{medium state}=[fill=white, draw=black, regular polygon, regular polygon sides=3, minimum width=1.3cm, inner sep=0pt, shape border rotate=180]
\tikzstyle{large state}=[fill=white, draw=black, regular polygon, regular polygon sides=3, minimum width=2.2cm, shape border rotate=180, inner sep=0pt]
\tikzstyle{wide state}=[fill=white, draw=black, shape=isosceles triangle, minimum width=0.8cm, shape border rotate=270, inner sep=1.4pt, minimum height=0.5cm, isosceles triangle apex angle=80]
\tikzstyle{wn}=[fill=white, draw=black, shape=circle, inner sep=1.5pt]
\tikzstyle{blue morphism}=[fill=white, draw={rgb,255: red,15; green,0; blue,150}, shape=rectangle, text={rgb,255: red,15; green,0; blue,150}, tikzit category=blue]
\tikzstyle{red morphism}=[fill=white, draw={rgb,255: red,150; green,0; blue,2}, shape=rectangle, text={rgb,255: red,150; green,0; blue,2}, tikzit category=red]
\tikzstyle{blue state}=[fill=white, draw={rgb,255: red,15; green,0; blue,150}, shape=circle, regular polygon, regular polygon sides=3, minimum width=0.8cm, shape border rotate=180, inner sep=0pt, text={rgb,255: red,15; green,0; blue,150}, tikzit category=blue]
\tikzstyle{blue node}=[fill={rgb,255: red,15; green,0; blue,150}, draw={rgb,255: red,15; green,0; blue,150}, shape=circle, tikzit category=blue, inner sep=1.5pt]
\tikzstyle{blue}=[text={rgb,255: red,15; green,0; blue,150}, tikzit draw={rgb,255: red,191; green,191; blue,191}, tikzit category=blue, tikzit fill=white, inner sep=0mm]
\tikzstyle{blue wide state}=[fill=white, draw={rgb,255: red,15; green,0; blue,150}, text={rgb,255: red,15; green,0; blue,150}, shape=isosceles triangle, minimum width=0.8cm, shape border rotate=270, inner sep=1.4pt, minimum height=0.5cm, isosceles triangle apex angle=80]
\tikzstyle{red node}=[fill={rgb,255: red,150; green,0; blue,2}, draw={rgb,255: red,150; green,0; blue,2}, shape=circle, inner sep=1.5pt]
\tikzstyle{Purple node}=[fill={rgb,255: red,120; green,0; blue,120}, draw={rgb,255: red,120; green,0; blue,120}, text={rgb,255: red,120; green,0; blue,120}, shape=circle, inner sep=1.5pt]
\tikzstyle{red}=[text={rgb,255: red,150; green,0; blue,2}, inner sep=0mm, tikzit fill=white, tikzit draw={rgb,255: red,191; green,191; blue,191}]
\tikzstyle{purple}=[text={rgb,255: red,150; green,0; blue,150}, inner sep=0mm, tikzit fill=white, tikzit draw={rgb,255: red,191; green,191; blue,191}]
\tikzstyle{white morphism}=[fill=white, draw=white, shape=rectangle, tikzit draw={rgb,255: red,139; green,139; blue,139}]
\tikzstyle{leak morphism}=[fill=white, draw={rgb,255: red,120; green,0; blue,85}, shape=rectangle, text={rgb,255: red,120; green,0; blue,85}, tikzit category=leak]
\tikzstyle{leak}=[text={rgb,255: red,120; green,0; blue,85}, inner sep=0mm, tikzit fill=white, tikzit draw={rgb,255: red,191; green,191; blue,191}, tikzit category=leak]
\tikzstyle{leak node}=[fill={rgb,255: red,120; green,0; blue,85}, draw={rgb,255: red,120; green,0; blue,85}, shape=circle, inner sep=1.5pt, tikzit category=leak]
\tikzstyle{horiz state}=[fill=white, draw=black, regular polygon, regular polygon sides=3, minimum width=1cm, shape border rotate=90, inner sep=0pt]

\tikzstyle{arrow}=[->]
\tikzstyle{dashed box}=[-, dashed]
\tikzstyle{blue arrow}=[-, draw={rgb,255: red,15; green,0; blue,150}, tikzit category=blue]
\tikzstyle{red arrow}=[-, draw={rgb,255: red,150; green,0; blue,2}, tikzit category=red]
\tikzstyle{purple arrow}=[->, draw={rgb,255: red,120; green,0; blue,120}, >=stealth, shorten <=2pt, shorten >=2pt]
\tikzstyle{protected purple arrow}=[->, draw={rgb,255: red,120; green,0; blue,120}, >=stealth, shorten <=2pt, shorten >=2pt, preaction={line width=1.8pt, white, draw}]
\tikzstyle{mapsto}=[{|->}]
\tikzstyle{double wire}=[-, double]
\tikzstyle{curly brace}=[-, draw=none, tikzit draw={rgb,255: red,128; green,0; blue,128}]
\tikzstyle{protected}=[-, preaction={line width=1.8pt,white,draw}]
\tikzstyle{leak arrow}=[-, tikzit draw={rgb,255: red,150; green,0; blue,120}]
\tikzstyle{protected leak arrow}=[-, tikzit draw={rgb,255: red,150; green,0; blue,120}]
\tikzstyle{hollow arrow}=[-, very thin, white, preaction={line width=0.7pt,draw={rgb,255: red,120; green,0; blue,85}}, tikzit category=leak, tikzit draw={rgb,255: red,150; green,0; blue,120}]
\tikzstyle{protected hollow arrow}=[-, very thin, white, preaction={line width=0.7pt,draw={rgb,255: red,120; green,0; blue,85},preaction={line width=2.1pt,white,draw}}, tikzit category=leak, tikzit draw={rgb,255: red,150; green,0; blue,120}]
\tikzstyle{over arrow}=[-, black, preaction={draw=white, double}]
\usepackage{graphicx}
\usepackage{tikzit}

\usepackage{pict2e}

\usepackage{pifont} 
\newcommand{\cmark}{\ding{51}}%
\newcommand{\xmark}{\ding{55}}%

\makeatletter
\newcommand{\adjunction}[2]{%
	#1%
	\mathrel{\vcenter{%
			\offinterlineskip\m@th
			\ialign{%
				\hfil$##$\hfil\cr
				\longrightharpoonup\cr
				\noalign{\kern-.3ex}
				\smallbot\cr
				\longleftharpoondown\cr
			}%
	}}%
	#2
}
\newcommand{\longrightharpoonup}{\relbar\joinrel\rightharpoonup}
\newcommand{\longleftharpoondown}{\leftharpoondown\joinrel\relbar}
\newcommand{\smallbot}{%
	\begingroup\setlength\unitlength{.15em}%
	\begin{picture}(1,1)
		\roundcap
		\polyline(0,0)(1,0)
		\polyline(0.5,0)(0.5,1)
	\end{picture}%
	\endgroup
}
\makeatother

\newtheorem{theorem}{Theorem}[section]
\newtheorem{lemma}[theorem]{Lemma}
\newtheorem{corollary}[theorem]{Corollary}
\newtheorem{definition}[theorem]{Definition}

\newtheorem{proposition}[theorem]{Proposition}

\theoremstyle{definition}				
\newtheorem{example}[theorem]{Example}
\newtheorem{remark}[theorem]{Remark}

\newcommand{\cat}[1]{{\mathbf{#1}}} 
\newcommand{\samp}{\texttt{samp}} 
\renewcommand{\det}{\text{det}} 

\title{Categorical Algebra of Conditional Probability}

\author{Mika Bohinen and Paolo Perrone}
\affil{University of Oxford}
\date{2025}

\begin{document}
\maketitle

\begin{abstract}
	In the field of categorical probability, one uses concepts and techniques from category theory, such as monads and monoidal categories, to study the structures of probability and statistics.
	In this paper, we connect some ideas from categorical algebra, namely weakly cartesian functors and natural transformations, to the idea of conditioning in probability theory, using Markov categories and probability monads.

	First of all, we show that under some conditions, the monad associated to a Markov category with conditionals has a weakly cartesian functor and weakly cartesian multiplication.
	In particular, we show that this is the case for the Giry monad on standard Borel spaces.

	We then connect this theory to existing results on statistical experiments. We show that for deterministic statistical experiments, the so-called standard measure construction (which can be seen as a generalization of the ``hypernormalizations'' introduced by Jacobs) satisfies a universal property, allowing an equivalent definition which does not rely on the existence of conditionals.
\end{abstract}

\tableofcontents

\section{Introduction}

Categorical probability is an emerging field that applies category theory to
probability and statistics. This approach seeks to:
\begin{itemize}
	\item Reformulate probabilistic concepts using abstract categorical methods,
	      particularly through diagrammatic reasoning, rather than relying on
	      specific analytical properties like cardinality or separability;
	      
	\item Simplify complex probabilistic proofs by working at a higher level
	      of abstraction, making existing results more transparent and enabling
	      the discovery of new theorems that were previously intractable;
	      
	\item Bridge probability theory with other mathematical disciplines,
	      particularly theoretical computer science, through the shared language
	      of category theory.
\end{itemize}

Markov categories~\cite{Cho_2019,Fritz-2020-Asyntheticapproachto} provide the
primary mathematical framework for this approach (see
\cite[Remark~2.2]{fritz2022free} for historical context). This framework has
enabled categorical reformulations and generalizations of fundamental results in
probability theory, including de Finetti's theorem~\cite{fritz2021definetti},
the Kolmogorov and Hewitt-Savage zero-one laws~\cite{fritzrischel2019zeroone},
and the ergodic decomposition theorem~\cite{ensarguet2023ergodic}. In
statistics, Markov categories have formalized concepts like sufficient
statistics~\cite{Fritz-2020-Asyntheticapproachto} and, through their connection
with probability monads~\cite{jacobs2018probabilitymonads}, led to new insights
into Blackwell's theorem on statistical experiments~\cite{Fritzetal-2023-Representablemarkovc}.

A central concept in Markov categories is \emph{conditioning} (see
\Cref{conditional}). It seeks to capture the traditional operations of
conditioning---including disintegrations, conditional expectations, and
Bayesian inverses---using string diagrams, but without having to worry about
measure-theoretic subtleties.

This work develops a categorical-algebraic perspective on conditionals in
Markov categories. We establish that, given suitable compatibility conditions
between a Markov category and a probability monad (see \Cref{main_results}),
conditionals give rise to two key properties of the monad:
\begin{itemize}
	\item The underlying functor is weakly cartesian;
	\item The multiplication of the monad is weakly cartesian.
\end{itemize}
These conditions were previously established for the distribution monad on
sets~\cite{Constantinetal-2020-Partialevaluationsan}. Our results extend this
characterization to the Giry monad on standard Borel spaces
(\Cref{giryWC}).

To understand intuitively how conditionals relate to weak pullback
preservation, consider a weak pullback diagram in a category $\cat{C}$ (such
as $\cat{Set}$ or $\cat{Meas}$) and its image under a functor
$P:\cat{C}\to\cat{C}$:
\[
	\begin{tikzcd}
		A & B \\
		C & D
		\arrow["f", from=1-1, to=1-2]
		\arrow["g"', from=1-1, to=2-1]
		\arrow["m", from=1-2, to=2-2]
		\arrow["n"', from=2-1, to=2-2]
	\end{tikzcd}
	\qquad\qquad
	\begin{tikzcd}
		PA & PB \\
		PC & PD
		\arrow["Pf", from=1-1, to=1-2]
		\arrow["Pg"', from=1-1, to=2-1]
		\arrow["Pm", from=1-2, to=2-2]
		\arrow["Pn"', from=2-1, to=2-2]
	\end{tikzcd}
\]
For the right diagram to be a weak pullback, we must establish that for any
elements $p\in PB$ and $q\in PC$ with $Pm(p)=Pn(q)$, there exists some
(not necessarily unique) element $r\in PA$ satisfying $Pf(r)=p$ and
$Pg(r)=q$. When $P$ is a probability monad, $p$ and $q$ represent probability
distributions, and the existence of conditionals provides a canonical choice
for $r\in PA$: the conditional product of $p$ and $q$ given $D$. This
construction yields a measure under which the observations $f$ and $g$ are
conditionally independent given their common coarse-graining $m\circ f=n\circ
	g$. Under suitable conditions, this measure can be shown to have support in
$A$. For the complete technical development, see \Cref{main_results}.

These results lead to new insights in the theory of statistical experiments. In
\cite{Constantinetal-2020-Partialevaluationsan}, as well as in
\cite{FritzPerrone-2018-Monadspartialevaluat} and \cite{constantin2023weak}, it
was shown that monads with weakly cartesian functor and multiplication have a particularly
interesting \emph{bar construction}, a simplicial set (or more generally
simplicial object) formed by an algebra of the monad $(A,a)$ as follows:
\[
	\begin{tikzcd}[sep=large]
		\cdots \ar{r} & TTTA \ar[shift left=0.8cm]{r}{\mu} \ar[shift left=0.4cm]{r}{T\mu} \ar["TTa"{description}]{r} & TTA \ar[shift left=0.4cm]{r}{\mu} \ar["Ta"{description}]{r}  \ar[shift left=0.4cm]{l}{T\eta} \ar[shift left=0.8cm]{l}{TT\eta} & TA \ar[shift left=0.4cm]{l}{T\eta}
	\end{tikzcd}
\]
The 0-1 truncation of this simplicial set is called the \emph{partial
	evaluation relation}, and can be seen as the relation connecting a formal
expression to its result or partial results (see the references above and our
\Cref{pev}). When both conditions hold, this relation is transitive.

The partial evaluation relation for the case of probability monads was shown to
be equivalent to what in probability and statistics is known as
\emph{second-order stochastic dominance}
\cite{FritzPerrone-2018-Monadspartialevaluat,Perrone18}, which measures how
``spread'' probability measures over the real line (or a vector
space) are. In \cite{Fritzetal-2023-Representablemarkovc}, partial evaluations
were connected to Markov categories. There, a synthetic and more general
version of Blackwell's theorem on statistical experiments was proven, showing
equivalence between stochastic dominance of certain measures (the
\emph{standard measures}, see \Cref{stat_exp}) and an order on statistical
experiments measuring their ``informativeness''. 

In the final section, we establish a universal property for the standard
measures of deterministic statistical experiments (which generalize Jacobs'
``hypernormalizations'', see \Cref{hypernorm}). Namely, they can be seen as the
\emph{coarsest decomposition of a measure which is still finer than or equal to
	the ``partition'' induced by a statistical experiment}. (This idea of
``partition induced by a function'' may remind the reader of \emph{descent},
but we leave a more thorough investigation of this analogy to future work.)

The structures and techniques of categorical algebra reach far beyond weakly
cartesian functors and multiplications, and we hope that this work is just the first one of a
fruitful thread of research, paving the way for an even deeper structural
understanding of probability, possibly connecting it to fibrations and descent
theory.

\noindent\textbf{Outline.}
In \Cref{background} we establish the main concepts that are used in the rest
of the work. In particular, in \Cref{markov} we give an overview of Markov
categories and their relationship with probability monads. In \Cref{pev} we
explain the main ideas behind partial evaluations, especially in the context of
probability monads. We then turn to statistical experiments in \Cref{stat_exp},
where we also give our definition of ``hypernormalization'' as a special case
of a standard measure (\Cref{def:hypernormalization}).

Our main results are stated and proven in \Cref{main_results}. In
\Cref{wcartmu} we prove that for an a.s.-compatibly representable Markov
category with conditionals, the multiplication of the monad is weakly cartesian
(\Cref{thm:pullbacks}). This in particular applies to the Giry monad on
standard Borel spaces (\Cref{pullbacks_giry}). In \Cref{wcartp} we prove that
for a representable Markov category with conditionals and satisfying the
so-called ``equalizer principle'' (\Cref{def:equalizer_principle}), the functor
underlying the monad is weakly cartesian (\Cref{thm:main}). Once again, this in
particular applies to the Giry monad, which then has weakly cartesian functor
and multiplication (\Cref{giryWC}). In \Cref{hypernorm} we show that a hypernormalization satisfies a
universal property, namely that it is the coarsest measure (in the
stochastic dominance order) which is compatible with the experiment
(\Cref{thm:uni-std}).

Finally, in \Cref{appendix} we give some technical results used in our proofs.

\section{Markov categories and partial evaluations}\label{background}

\subsection{Markov categories and probability monads}\label{markov}

Let us briefly recall the main definitions of Markov categories and their relationship with probability monads.
For more details, see the original papers \cite{Cho_2019,Fritz-2020-Asyntheticapproachto,Fritzetal-2023-Representablemarkovc}.

\begin{definition}
	A \textbf{copy-discard} (\textbf{CD}) \textbf{category} is a symmetric monoidal category $(\cat{C},\otimes,I)$ where every object is equipped with a distinguished commutative comonoid structure, compatible with the tensor product. We denote the comonoid structure maps as follows:
	\begin{equation*}
		\texttt{copy} \;=
		\begin{tikzpicture}
	\begin{pgfonlayer}{nodelayer}
		\node [style=none] (0) at (0, -1) {};
		\node [style=none] (1) at (0, -1.5) {$X$};
		\node [style=bn] (2) at (0, 0) {};
		\node [style=none] (3) at (-1, 1) {};
		\node [style=none] (4) at (1, 1) {};
		\node [style=none] (5) at (-1, 1.5) {$X$};
		\node [style=none] (6) at (1, 1.5) {$X$};
	\end{pgfonlayer}
	\begin{pgfonlayer}{edgelayer}
		\draw (2) to (0.center);
		\draw [in=165, out=-90] (3.center) to (2);
		\draw [in=270, out=15] (2) to (4.center);
	\end{pgfonlayer}
\end{tikzpicture}
		\qquad\mbox{and}\qquad
		\texttt{del} \;=
		\begin{tikzpicture}
	\begin{pgfonlayer}{nodelayer}
		\node [style=none] (0) at (0, -0.5) {};
		\node [style=none] (1) at (0, -1) {$X$};
		\node [style=bn] (2) at (0, 1) {};
	\end{pgfonlayer}
	\begin{pgfonlayer}{edgelayer}
		\draw (2) to (0.center);
	\end{pgfonlayer}
\end{tikzpicture}

	\end{equation*}
	
	A \textbf{Markov category} is a CD category where for every $f:X \to Y$ in $\cat{C}$ the following equality holds,
	\begin{equation}
		\begin{tikzpicture}
	\begin{pgfonlayer}{nodelayer}
		\node [style=morphism] (0) at (-2.25, 0) {$f$};
		\node [style=bn] (1) at (-2.25, 1.25) {};
		\node [style=none] (2) at (-2.25, -1.25) {};
		\node [style=none] (3) at (0, 0) {$=$};
		\node [style=none] (4) at (2, -1.25) {};
		\node [style=bn] (5) at (2, 0.75) {};
	\end{pgfonlayer}
	\begin{pgfonlayer}{edgelayer}
		\draw (2.center) to (0);
		\draw (0) to (1);
		\draw (5) to (4.center);
	\end{pgfonlayer}
\end{tikzpicture}
		\label{eq:main-def}
	\end{equation}
	or equivalently, where the monoidal unit $I$ is terminal.
\end{definition}

Canonical examples of Markov categories are \emph{categories of Markov kernels}, hence the name.
These are also the most relevant examples for the purposes of the present paper.

\begin{example}
	The category \( \cat{Stoch} \) is specified via the following data.
	\begin{itemize}
		\item Objects are measurable spaces, i.e., pairs \( (X,
		      \mathcal{A}) \) where \( X \in \cat{Set} \) and \( \mathcal{A} \)
		      is a sigma-algebra on \( X \);
		\item Morphisms \( (X, \mathcal{A}) \to  (Y, \mathcal{B}) \) are Markov
		      kernels of entries \( k(B|x) \), for \( x \in  X \) and \( B \in
		      \mathcal{B} \). That is to say, $k(B|-)$ is a measurable function $X\to\mathbb{R}$ for all $B\in\mathcal{B}$, and \( k(-|x) \) is a probability	measure on \( (Y,\mathcal{B}) \) for all \( x \in  X \).
		\item The identity \( (X, \mathcal{A}) \to  (X, \mathcal{A}) \) is given by the Dirac delta, i.e.,
		      \begin{equation}
			      \delta(A|x) = \begin{cases}
				      1, & \text{ if }x \in  A \\
				      0, & \text{ otherwise.}
			      \end{cases}
			      \label{eq:diraciden}
		      \end{equation}
		\item The composition of two kernels \( k: (X, \mathcal{A}) \to (Y,
		      \mathcal{B}) \) and \( h:(Y, \mathcal{B}) \to  (Z, \mathcal{C})
		      \) is given by taking the Lebesgue integral, i.e.,
		      \begin{equation}
			      (h\circ k)(C|x) = \int_Y h(C|y)k(dy|x)
			      \label{eq:complebeg}
		      \end{equation}
		      for all \( x \in  X \) and \( C \in  \mathcal{C} \).
		\item The monoidal structure is given by the usual product of measurable spaces; Fubini's theorem is used to construct the product of kernels and measures.
		\item The ``copy'' map $\texttt{copy}:X\to X\otimes X$ is defined by
		      \[
			      \texttt{copy}(A\times B| x) \;=\; \begin{cases}
				      1 & \text{ if }x\in A\cap B ; \\
				      0 & \text{ otherwise.}
			      \end{cases}
		      \]
		\item The ``delete'' map $\texttt{del}:X\to I$ is the unique kernel to the one-point space.
	\end{itemize}
	
	The Markov category \( \cat{BorelStoch} \) is defined to be the full
	subcategory of \( \cat{Stoch} \) where the objects are standard
	Borel spaces.
\end{example}

\begin{definition}[\cite{Fritzetal-2023-Representablemarkovc}]
	Let \( \cat{C} \) be a Markov category. A morphism \( f: A \to X \) in \( \cat{C} \) is said to be \textbf{deterministic} if the following equality holds
	\begin{equation}
		\begin{tikzpicture}
	\begin{pgfonlayer}{nodelayer}
		\node [style=bn] (0) at (-3, 0) {};
		\node [style=none] (1) at (-2, 1) {};
		\node [style=none] (2) at (-2, 2.5) {$X$};
		\node [style=none] (3) at (-4, 1) {};
		\node [style=none] (4) at (-4, 2.5) {$X$};
		\node [style=none] (5) at (-3, -2.5) {$A$};
		\node [style=none] (6) at (-3, -2) {};
		\node [style=morphism] (7) at (-3, -1) {$f$};
		\node [style=bn] (8) at (3.5, -0.75) {};
		\node [style=none] (9) at (4.5, 2) {};
		\node [style=none] (10) at (4.5, 2.5) {$X$};
		\node [style=none] (11) at (2.5, 2) {};
		\node [style=none] (12) at (2.5, 2.5) {$X$};
		\node [style=none] (13) at (3.5, -2.5) {$A$};
		\node [style=none] (14) at (3.5, -2) {};
		\node [style=morphism] (15) at (2.5, 1) {$f$};
		\node [style=morphism] (16) at (4.5, 1) {$f$};
		\node [style=none] (17) at (0, 0) {$=$};
		\node [style=none] (18) at (-4, 2) {};
		\node [style=none] (19) at (-2, 2) {};
		\node [style=none] (20) at (2.5, 0.25) {};
		\node [style=none] (21) at (4.5, 0.25) {};
	\end{pgfonlayer}
	\begin{pgfonlayer}{edgelayer}
		\draw [in=270, out=165] (0) to (3.center);
		\draw [in=270, out=15] (0) to (1.center);
		\draw (0) to (7);
		\draw (7) to (6.center);
		\draw (15) to (11.center);
		\draw (16) to (9.center);
		\draw (8) to (14.center);
		\draw (18.center) to (3.center);
		\draw (19.center) to (1.center);
		\draw [in=165, out=-90] (20.center) to (8);
		\draw [in=270, out=15] (8) to (21.center);
		\draw (20.center) to (15);
		\draw (16) to (21.center);
	\end{pgfonlayer}
\end{tikzpicture}
		\label{eq:deterministic}
	\end{equation}
	The subcategory of \( \cat{C} \) that consists of only deterministic morphisms is denoted by \( \cat{C}_{\text{det}} \).
\end{definition}

\begin{example}
	The deterministic morphisms of $\cat{BorelStoch}$ are exactly those Markov kernels $k_f:X\to Y$ which are induced by a measurable function as follows:
	\[
		k_f(B|x) \;=\; \begin{cases}
			1 & \text{ if } f(x)\in B ; \\
			0 & \text{ otherwise. }
		\end{cases}
	\]
\end{example}

One of the most important pieces of structures for this paper is the idea of \emph{conditioning}.

\begin{definition}\label{conditional}
	Let \( f: A \to  X \otimes Y \) be a morphism in a Markov category \(
	\cat{C} \). A map \( f|_X: X \otimes A \to Y \) is called a
	\textbf{conditional} of \( f \) with respect to \( X \) if the following equation holds.
	\begin{equation}
		\begin{tikzpicture}
	\begin{pgfonlayer}{nodelayer}
		\node [style=medium box] (0) at (-2.75, 0.5) {$f$};
		\node [style=none] (1) at (-3.25, 3.75) {};
		\node [style=none] (2) at (-3.25, 4.25) {$X$};
		\node [style=none] (3) at (-2.25, 3.75) {};
		\node [style=none] (4) at (-2.25, 4.25) {$Y$};
		\node [style=none] (5) at (-2.75, -3.75) {$A$};
		\node [style=none] (6) at (-2.75, -3.25) {};
		\node [style=none] (7) at (-3.25, 1) {};
		\node [style=none] (8) at (-2.25, 1) {};
		\node [style=none] (9) at (0, 0) {$=$};
		\node [style=none] (10) at (4.5, -3.75) {$A$};
		\node [style=none] (11) at (4.5, -3.25) {};
		\node [style=bn] (13) at (4.5, -2.25) {};
		\node [style=medium box] (14) at (3.5, -0.75) {$f$};
		\node [style=medium box] (15) at (4.75, 2.5) {$\,f|_X\,$};
		\node [style=bn] (16) at (3, 0.75) {};
		\node [style=none] (17) at (4.75, 3.75) {};
		\node [style=none] (18) at (4.75, 4.25) {$Y$};
		\node [style=none] (19) at (4, 1.75) {};
		\node [style=none] (20) at (5.5, -1.25) {};
		\node [style=none] (21) at (2, 4.25) {$X$};
		\node [style=none] (22) at (2, 1.75) {};
		\node [style=bn] (23) at (4, 0.375) {};
		\node [style=none] (24) at (4, -0.75) {};
		\node [style=none] (25) at (3, -0.75) {};
		\node [style=none] (26) at (2, 3.75) {};
		\node [style=none] (27) at (5.5, 2.25) {};
		\node [style=none] (28) at (3.5, -1.25) {};
		\node [style=none] (29) at (4, 2.25) {};
	\end{pgfonlayer}
	\begin{pgfonlayer}{edgelayer}
		\draw (6.center) to (0);
		\draw (8.center) to (3.center);
		\draw (7.center) to (1.center);
		\draw (24.center) to (23);
		\draw (25.center) to (16);
		\draw (11.center) to (13);
		\draw [in=270, out=15] (13) to (20.center);
		\draw [in=270, out=15] (16) to (19.center);
		\draw (15) to (17.center);
		\draw [in=270, out=165] (16) to (22.center);
		\draw (26.center) to (22.center);
		\draw (27.center) to (20.center);
		\draw [in=165, out=-90] (28.center) to (13);
		\draw (14) to (28.center);
		\draw (29.center) to (19.center);
	\end{pgfonlayer}
\end{tikzpicture}
		\label{eq:conditional}
	\end{equation}
	
	A Markov category is said to \textbf{have all conditionals} (or more briefly, just \textbf{have conditionals}) if every morphism admits a conditional with respect to any of its outputs.
\end{definition}

\begin{example}
	The category $\cat{BorelStoch}$ has all conditionals, and they correspond to \textbf{regular conditional probability distributions} \cite[Example~11.3]{Fritz-2020-Asyntheticapproachto}.
	For example, given a probability measure $\psi: I \to X \otimes Y$ with
	$X$-marginal $\psi_X$, a conditional distribution is a kernel
	$\psi|_X: X \to Y$ satisfying
	\begin{equation*}
		\psi(S \times T) = \int_{x \in S}\psi|_X(T \mid x)\,\psi_X(dx) .
	\end{equation*}
	Such a kernel always exists if $X$ and $Y$ are standard Borel.
\end{example}

\begin{definition}
	\label{def:almost}
	Let $p : A \to X $ be a morphism and $f,g: X \to Y $ two parallel
	morphisms. We say that $f$ and $g$ are $p$-almost surely equal, denoted
	$f=_{p\text{-a.s.}} g$, if we have
	\begin{equation}
		\label{eq:aes}
		\begin{tikzpicture}
	\begin{pgfonlayer}{nodelayer}
		\node [style=bn] (0) at (-3.25, -0.25) {};
		\node [style=morphism] (1) at (-2.25, 1.25) {$f$};
		\node [style=none] (2) at (-2.25, 2.25) {};
		\node [style=none] (3) at (-2.25, 2.75) {$Y$};
		\node [style=none] (4) at (-4.25, 0.75) {};
		\node [style=none] (5) at (-4.25, 2.75) {$X$};
		\node [style=morphism] (6) at (-3.25, -1.375) {$p$};
		\node [style=none] (7) at (-3.25, -2.5) {};
		\node [style=none] (8) at (-3.25, -3) {$A$};
		\node [style=bn] (9) at (3, -0.25) {};
		\node [style=morphism] (10) at (4, 1.25) {$\vphantom{f}g$};
		\node [style=none] (11) at (4, 2.25) {};
		\node [style=none] (12) at (4, 2.75) {$Y$};
		\node [style=none] (13) at (2, 0.75) {};
		\node [style=none] (14) at (2, 2.75) {$X$};
		\node [style=morphism] (15) at (3, -1.375) {$p$};
		\node [style=none] (16) at (3, -2.5) {};
		\node [style=none] (17) at (3, -3) {$A$};
		\node [style=none] (18) at (0, 0) {$=$};
		\node [style=none] (19) at (-2.25, 0.75) {};
		\node [style=none] (20) at (-4.25, 2.25) {};
		\node [style=none] (21) at (4, 0.75) {};
		\node [style=none] (22) at (2, 2.25) {};
	\end{pgfonlayer}
	\begin{pgfonlayer}{edgelayer}
		\draw [in=630, out=165] (0) to (4.center);
		\draw (1) to (2.center);
		\draw (7.center) to (6);
		\draw (6) to (0);
		\draw [in=-90, out=180] (9) to (13.center);
		\draw (10) to (11.center);
		\draw (16.center) to (15);
		\draw (15) to (9);
		\draw (19.center) to (1);
		\draw [in=15, out=-90] (19.center) to (0);
		\draw (20.center) to (4.center);
		\draw (21.center) to (10);
		\draw [in=-90, out=15] (9) to (21.center);
		\draw (22.center) to (13.center);
	\end{pgfonlayer}
\end{tikzpicture}
	\end{equation}
\end{definition}

\begin{example}
	\label{ex:almostsure}
	In $\cat{BorelStoch}$, Definition~\ref{def:almost} gives the standard
	notion of almost surely equality. More specifically, given Markov kernels
	$f,g: (X, \Sigma_{A}) \to (Y, \Sigma_{B})$ and a probability measure $\nu:
		I \to X$, the relation $f =_{\nu\text{-a.s.}} g$ Equation~\ref{eq:aes}
	becomes the condition
	\begin{equation*}
		\int_S f(T|x)\nu(dx) = \int_S g(T|x)\nu(dx),
	\end{equation*}
	for all $S \in \Sigma_X$ and $T \in \Sigma_Y$. This is the same as saying that
	$f(T|-)$ and $g(T|-)$ are $\nu$-almost everywhere equal for all $T$.
\end{example}

An important example of Markov categories are those arising as Kleisli categories of a monad (sometimes called a \emph{probability monad}~\cite{jacobs2018probabilitymonads}).
The following phenomenon is well known.

\begin{proposition}
	\label{prop:importantproperty}
	Let \( (P, \mu, \delta) \) be a symmetric monoidal monad on some symmetric
	monoidal category $(D, \otimes, I)$. Then \( \cat{D}_P \) is a symmetric
	monoidal category, with:
	\begin{itemize}
		\item the same monoidal product as the one in \( \cat{D} \);
		      
		\item the tensor product of morphisms represented by \( f:A \to PX \)
		      and \( g:B \to PY \) being represented by the composite
		      \begin{equation}
			      A\otimes B \xrightarrow{f\otimes g} PX\otimes PY \xrightarrow{\nabla} P(X\otimes Y).
			      \label{eq:repr}
		      \end{equation}
	\end{itemize}
	Moreover, the inclusion \( \cat{D} \to  \cat{D}_P \) is strict symmetric monoidal.
\end{proposition}

The final statement about the inclusion implies that if $ X \in \cat{D} $ has a
distinguished comonoid structure, then so does $ X \in \cat{D}_P $.

\begin{definition}
	A monad $(P,\mu, \delta)$ on $\cat{D}$ is said to be
	\textbf{affine} if $PI \cong I$.
\end{definition}

Thus, if $P$ is affine and $I$ is terminal, then $PI \in \cat{D}_P$ is also
terminal and we get the following result.
\begin{corollary}
	Let \( (P,\mu ,\delta ) \) be a symmetric monoidal affine monad on a Markov
	category \( \cat{D} \). Then the Kleisli category \( \cat{D}_P \) is again a
	Markov category in a canonical way.
	\label{cor:affine}
\end{corollary}

\begin{example}
	$\cat{BorelStoch}$ can be seen as the Kleisli category of the Giry monad on standard Borel spaces.
\end{example}

Markov categories of this kind have a particular structure, very convenient for
the purposes of probability theory. This abstracts the idea that a Markov
kernel $X\to Y$ is equivalently specified by a measurable function $X\to PY$,
where $PY$ is the space of probability measures over $Y$, equipped with the
canonical sigma-algebra given by the Giry monad~\cite{giry}.

\begin{definition}
	Let $\cat{C} $ be a Markov category and $X\in \cat{C} $ an object. A
	\textbf{distribution object} for $X$ is an object $PX$ together with a
	morphism $\samp_X: PX \to X $ so that the induced map
	\begin{equation*}
		\samp_X \circ -: \cat{C}_{\text{det}}(A, PX) \to C(A, X)
	\end{equation*}
	is a bijection for all $A \in \cat{C} $. We denote the inverse of this map by
	\begin{equation*}
		(-)^\#: \cat{C}(A, X) \to \cat{C}_{\text{det}}(A, PX) ,
	\end{equation*}
	and we set $\delta_X=(1_X)^\#:X\to PX$.
\end{definition}

\begin{definition}
	\label{def:distfunctor}
	Let $\cat{C}$ be a Markov category. We say that $\cat{C}$ is
	\textbf{representable} if every object has a distribution object.
\end{definition}

Distribution objects, if they exist for all $X$, assemble together to give a
right adjoint $P:\cat{C} \to \cat{C}_{\text{det}} $ to the inclusion functor
$\cat{C}_{\text{det}} \hookrightarrow \cat{C} $. The unit of the adjunction is
given by the maps $\delta:X\to PX$, and the counit by $\samp:PX\to X$. The
resulting monad on $\cat{C}_{\text{det}}$ is given by $(P,P\samp,\delta)$.

\begin{example}
	The Giry monad on $\cat{BorelStoch}$ is in this form. In particular,
	\begin{itemize}
		\item The unit $\delta:X\to PX$ assigns to each point $x\in X$ the corresponding Dirac measure $\delta_x$;
		\item The counit of the adjunction $\samp:PX\to X$ is the kernel
		      \[
			      \samp(A|p) \;=\; p(A)
		      \]
		      for all $p\in PX$ and all measurable subsets $A\subseteq X$.
	\end{itemize}
\end{example}

\begin{definition}
	Let $\cat{C}$ be a Markov category. We say that $\cat{C}$ is
	\textbf{a.s.-compatibly representable} if it is representable and for any
	morphism $p: \Theta \to A$, the natural bijection
	\[
		\cat{C}_{\text{det}}(A,PX) \;\cong\; \cat{C}(A,X)
	\]
	respects almost sure equality. That is to say, for
	all $f,g: A \to X$, we have
	\begin{equation}
		\label{eq:almostsurely}
		f^\# =_{p\text{-a.s}} g^\# \iff f =_{p\text{-a.s.}} g.
	\end{equation}
\end{definition}

In later proofs we shall make use of another equivalent characterization of a.s.-compatibly representable Markov categories.
\begin{definition}
	Let $\cat{C}$ be a representable Markov category. It is said to
	satisfy the \textbf{sampling cancellation property} if, for any three
	morphisms $f,g:X\otimes A \to Y$ and $ p:A \to X$, the following
	implication holds
	\begin{equation}
		\label{eq:samp-canc}
		\begin{tikzpicture}
	\begin{pgfonlayer}{nodelayer}
		\node [style=none] (0) at (-10, -3.25) {$A$};
		\node [style=none] (1) at (-10, -2.75) {};
		\node [style=none] (2) at (-12, 4) {$X$};
		\node [style=none] (3) at (-12, 1.75) {};
		\node [style=none] (4) at (-9.5, 4) {$Y$};
		\node [style=none] (5) at (-9.5, 3.5) {};
		\node [style=bn] (6) at (-10, -1.75) {};
		\node [style=bn] (7) at (-11, 0.75) {};
		\node [style=morphism] (8) at (-11, -0.25) {$p$};
		\node [style=medium box] (9) at (-9.5, 2.25) {$f$};
		\node [style=none] (10) at (-10, 1.75) {};
		\node [style=none] (11) at (-9, -0.75) {};
		\node [style=none] (12) at (-7, 0) {$=$};
		\node [style=none] (13) at (-3.5, -3.25) {$A$};
		\node [style=none] (14) at (-3.5, -2.75) {};
		\node [style=none] (15) at (-5.5, 4) {$X$};
		\node [style=none] (16) at (-5.5, 1.75) {};
		\node [style=none] (17) at (-3, 4) {$Y$};
		\node [style=none] (18) at (-3, 3.5) {};
		\node [style=bn] (19) at (-3.5, -1.75) {};
		\node [style=bn] (20) at (-4.5, 0.75) {};
		\node [style=morphism] (21) at (-4.5, -0.25) {$p$};
		\node [style=medium box] (22) at (-3, 2.25) {$g$};
		\node [style=none] (23) at (-3.5, 1.75) {};
		\node [style=none] (24) at (-2.5, -0.75) {};
		\node [style=none] (25) at (0, 0) {$\implies$};
		\node [style=none] (26) at (4.5, -3.25) {$A$};
		\node [style=none] (27) at (4.5, -2.75) {};
		\node [style=none] (28) at (2.5, 4) {$X$};
		\node [style=none] (29) at (2.5, 1.75) {};
		\node [style=none] (30) at (5, 4) {$PY$};
		\node [style=none] (31) at (5, 3.5) {};
		\node [style=bn] (32) at (4.5, -1.75) {};
		\node [style=bn] (33) at (3.5, 0.75) {};
		\node [style=morphism] (34) at (3.5, -0.25) {$p$};
		\node [style=medium box] (35) at (5, 2.25) {$f^\#$};
		\node [style=none] (36) at (4.5, 1.75) {};
		\node [style=none] (37) at (5.5, -0.75) {};
		\node [style=none] (38) at (7.5, 0) {$=$};
		\node [style=none] (39) at (11, -3.25) {$A$};
		\node [style=none] (40) at (11, -2.75) {};
		\node [style=none] (41) at (9, 4) {$X$};
		\node [style=none] (42) at (9, 1.75) {};
		\node [style=none] (43) at (11.5, 4) {$PY$};
		\node [style=none] (44) at (11.5, 3.5) {};
		\node [style=bn] (45) at (11, -1.75) {};
		\node [style=bn] (46) at (10, 0.75) {};
		\node [style=morphism] (47) at (10, -0.25) {$p$};
		\node [style=medium box] (48) at (11.5, 2.25) {$g^\#$};
		\node [style=none] (49) at (11, 1.75) {};
		\node [style=none] (50) at (12, -0.75) {};
		\node [style=none] (51) at (-12, 3.5) {};
		\node [style=none] (52) at (-5.5, 3.5) {};
		\node [style=none] (53) at (2.5, 3.5) {};
		\node [style=none] (54) at (9, 3.5) {};
		\node [style=none] (55) at (-9, 1.75) {};
		\node [style=none] (56) at (-2.5, 1.75) {};
		\node [style=none] (57) at (5.5, 1.75) {};
		\node [style=none] (58) at (12, 1.75) {};
		\node [style=none] (59) at (-11, -0.75) {};
		\node [style=none] (60) at (-4.5, -0.75) {};
		\node [style=none] (61) at (3.5, -0.75) {};
		\node [style=none] (62) at (10, -0.75) {};
	\end{pgfonlayer}
	\begin{pgfonlayer}{edgelayer}
		\draw [in=270, out=15] (6) to (11.center);
		\draw [in=270, out=15] (7) to (10.center);
		\draw (7) to (8);
		\draw [in=270, out=165] (7) to (3.center);
		\draw (6) to (1.center);
		\draw (9) to (5.center);
		\draw [in=270, out=15] (19) to (24.center);
		\draw [in=270, out=15] (20) to (23.center);
		\draw (20) to (21);
		\draw [in=270, out=165] (20) to (16.center);
		\draw (19) to (14.center);
		\draw (22) to (18.center);
		\draw [in=270, out=15] (32) to (37.center);
		\draw [in=270, out=15] (33) to (36.center);
		\draw (33) to (34);
		\draw [in=270, out=165] (33) to (29.center);
		\draw (32) to (27.center);
		\draw (35) to (31.center);
		\draw [in=270, out=15] (45) to (50.center);
		\draw [in=270, out=15] (46) to (49.center);
		\draw (46) to (47);
		\draw [in=270, out=165] (46) to (42.center);
		\draw (45) to (40.center);
		\draw (48) to (44.center);
		\draw (51.center) to (3.center);
		\draw (52.center) to (16.center);
		\draw (53.center) to (29.center);
		\draw (54.center) to (42.center);
		\draw (55.center) to (11.center);
		\draw (56.center) to (24.center);
		\draw (57.center) to (37.center);
		\draw (58.center) to (50.center);
		\draw (8) to (59.center);
		\draw (21) to (60.center);
		\draw (34) to (61.center);
		\draw (47) to (62.center);
		\draw [in=165, out=-90] (59.center) to (6);
		\draw [in=165, out=-90] (60.center) to (19);
		\draw [in=165, out=-90] (61.center) to (32);
		\draw [in=165, out=-90] (62.center) to (45);
	\end{pgfonlayer}
\end{tikzpicture}
	\end{equation}
\end{definition}
\begin{remark}
	Remembering that $f = \samp\, \circ f^\# $ for all morphisms in
	$\cat{C}$, the condition above is named so because it amounts to canceling
	the $\samp$ from $f$ and $g$.
\end{remark}

\begin{proposition}[{\cite[Proposition~3.24]{Fritzetal-2023-Representablemarkovc}}]
	Let $\cat{C}$ be a representable Markov category. Then $\cat{C}$ is
	a.s.-compatibly representable if and only if it satisfies the sampling cancellation property.
\end{proposition}

\subsection{Partial evaluations}\label{pev}

We now turn our attention to partial evaluations. The main idea, quite simple, is that for example, a
formal expression like $1 + 2 + 3$ can be totally evaluated to $6$, but it can also
be partially evaluated to $3 + 3$ or $1 + 5$.

It is well known that partial evaluations for probability monads correspond exactly to the second-order stochastic dominance~\cite{FritzPerrone-2018-Monadspartialevaluat}, and their relationship with Markov categories was explored in \cite{Fritzetal-2023-Representablemarkovc}.
This paper extends that relationship, but before we do that, let us recall some of the main known ideas. For more details, see~\cite{FritzPerrone-2018-Monadspartialevaluat,Constantinetal-2020-Partialevaluationsan}.

\begin{definition}
	Let $\cat{C}$ be a category and $X \in \cat{C}$. An $S$-shaped
	\textbf{generalized element} of $X$ is a morphism $p:S \to X$ for some $S \in \cat{C}$.
	By abuse of notation we will also write $p \in X$ when $p$ is a generalized element.
\end{definition}

Putting $\cat{C} = \cat{Set}$ and $S = \{*\}$ recovers the usual notion of elements in a set.
Indeed, we are particularly interested in the case where $\cat{C}$ is a Markov category and $p$ is a state.

\begin{definition}
	Let $(T, \mu, \eta)$ be a monad on some category $\cat{C}$, and
	$A\in\cat{C}$. An $S$-shaped \textbf{generalized formal expression} on $A$
	is an $S$-shaped generalized element $p \in TA$ for some $S\in \cat{C}$.
\end{definition}

\begin{definition}
	Let $p,q \in TA$ be $S$-shaped generalized formal expressions on a $T$-algebra
	$(A,e)$. A \textbf{partial evaluation} from $p$ into $q$ is an $S$-shaped generalized
	element $k \in TTA$ such that the following diagram commutes
	\[
		\begin{tikzcd} & S \\ TA & TTA & TA. \arrow["p"', from=1-2, to=2-1]
			\arrow["k", from=1-2, to=2-2] \arrow["q", from=1-2, to=2-3]
			\arrow["\mu", from=2-2, to=2-1] \arrow["Te"', from=2-2, to=2-3]
		\end{tikzcd}
	\]
\end{definition}

\begin{definition}
	The \textbf{partial evaluation relation} on $TA$ is defined as follows: given $p,q\in TA$, we say that $p\to q$ or $p\le q$ if and only if there exists a partial evaluation from $p$ to $q$.
\end{definition}

This relation is always reflexive.
Here is a sufficient condition for it to be transitive.

Recall that a commutative diagram
\[\begin{tikzcd}
		A & B \\
		C & D
		\arrow["f", from=1-1, to=1-2]
		\arrow["g"', from=1-1, to=2-1]
		\arrow["m", from=1-2, to=2-2]
		\arrow["n"', from=2-1, to=2-2]
	\end{tikzcd}\]
is a \textbf{weak pullback} if, given $p: S \to C$ and $q:S \to B$ such that $n \circ p = m \circ q$ then there
exists $r: S \to A$ such that
\[\begin{tikzcd}
		S \\
		& A & B \\
		& C & D
		\arrow["r"{description}, dashed, from=1-1, to=2-2]
		\arrow["q", curve={height=-18pt}, from=1-1, to=2-3]
		\arrow["p"', curve={height=18pt}, from=1-1, to=3-2]
		\arrow["f", from=2-2, to=2-3]
		\arrow["g"', from=2-2, to=3-2]
		\arrow["m", from=2-3, to=3-3]
		\arrow["n"', from=3-2, to=3-3]
	\end{tikzcd}\]
is commutative.

(This is almost the same as a pullback except that we have dropped the uniqueness condition.)

\begin{definition}
	Let $(T, \mu, \eta)$ be a monad on some category $\cat{C}$.
	\begin{itemize}
		\item We say that the functor $T$ is \textbf{weakly Cartesian} if and only if it preserves weak pullbacks.
		\item We say that $\mu$
		      is \textbf{weakly Cartesian} if the diagram
		      \[\begin{tikzcd}
				      TTX & TTY \\
				      TX & TY
				      \arrow["TTf", from=1-1, to=1-2]
				      \arrow["{\mu_X}"', from=1-1, to=2-1]
				      \arrow["{\mu_Y}", from=1-2, to=2-2]
				      \arrow["Tf"', from=2-1, to=2-2]
			      \end{tikzcd}\]
		      is a weak pullback for all $X,Y \in \cat{C}$ and $f:X \to Y$.
	\end{itemize}
\end{definition}

\begin{remark}\label{rem:wc-terminology}
	When both conditions above hold, the monad is called \emph{weakly
	cartesian} by Weber~\cite[Definition~11.2]{weber2004generic}, who
	additionally requires $\eta$ to be weakly cartesian. We do not need
	this stronger condition for our results. The combined property of
	weakly cartesian functor and multiplication was studied under the
	name ``(BC)'' by Clementino--Hofmann--Janelidze~\cite{clementino2014monads}
	and Constantin--Fritz--Perrone--Shapiro~\cite{Constantinetal-2020-Partialevaluationsan};
	we avoid this terminology, as ``Beck--Chevalley condition'' has an
	established and different meaning in the context of fibrations and
	descent~\cite{benabou1970descent}.
\end{remark}

\begin{proposition}[{\cite[Proposition 4.1]{FritzPerrone-2018-Monadspartialevaluat}}]
	\label{prop:partialrel}
	Let $(T, \mu, \eta)$ be a monad on a category $\cat{C}$, $S \in
		\cat{C}$, and $A \in \cat{C}^T$ a $T$-algebra with $e: TA \to A$ the algebra
	map. The partial evaluation relation on $\cat{C}(S, TA)$ is always
	reflexive, and if the multiplication $\mu$ is weakly cartesian, the
	relation is also transitive.
\end{proposition}
\begin{proof}
	To see reflexivity, given $p:S\to TA$, take $T\eta\circ p:S\to TTA$. We
	have that 
	\[
		\mu\circ T\eta\circ p \;=\; p \;=\; Te\circ T\eta\circ p 
	\]
	using the monad and the algebra unit condition, and that shows that
	$T\eta\circ p$ is a partial evaluation from $p$ to itself.
	
	To see transitivity, suppose that $\mu$ is weakly cartesian, and consider
	``composable'' partial evaluations $r$ and $s$. That is, let
	$p,q,t:S\to TA$, and $r,s:S\to TTA$ be such that $\mu\circ r=p$,
	$Te\circ s=t$, and $Te\circ r= \mu\circ s=q$. We then have the solid
	arrows in the following commutative diagram:
	\[
		\begin{tikzcd}[sep=large]
			&& TTA \ar{r}{\mu} \ar{dr}[very near end,swap,inner sep=0.5mm]{Te} & TA \\
			S \ar[dashed]{r}{n} \ar[bend left=15]{urr}{r} \ar[bend right=15]{drr}[swap]{s} & TTTA \ar{ur}{\mu} \ar{r}{T\mu} \ar{dr}[swap,inner sep=0.5mm]{TTe} & TTA \ar{ur}[swap,near end]{\mu} \ar{dr}[very near end,inner sep=0.5mm]{Te} & TA \\
			&& TTA \ar{r}[swap]{Te} \ar{ur}[near end,inner sep=0.5mm]{\mu} & TA
		\end{tikzcd}
	\]
	Now notice that, by hypothesis, the diamond diagram involving $S$
	commutes ($Te\circ r= \mu\circ s$). Since the diamond involving $TTTA$
	is weakly cartesian (it is a naturality square for $\mu$), there exists
	a (possibly non-unique) arrow $n:S\to TTTA$ making the entire diagram
	commute. Forming now the partial evaluation $T\mu\circ n:S\to TTA$, we
	have that $\mu\circ T\mu\circ n = \mu\circ\mu\circ n=\mu\circ r=p$, and
	$Te\circ T\mu\circ n=Te\circ TTe\circ n = Te\circ s = t$. 
\end{proof}

\begin{remark}
	While in the case above the partial evaluation relation is a preorder,
	which can be seen as a truncation of a category, partial evaluations
	and their composition hardly ever form a category. The situation is
	more interesting when $T$ also preserves weak pullbacks, in which case, while we still do not have a
	category (or a quasi-category) in general, we have an interesting
	higher compositional structure. See
	\cite{Constantinetal-2020-Partialevaluationsan,constantin2023weak} for
	the details.
\end{remark}

From the point of view of the Kleisli category, these weak cartesian conditions
look as follows.

\begin{proposition}\label{prop:BC_Kleisli}
	Let $(T, \mu, \eta)$ be a monad on a category $\cat{C}$. A commutative
	square in the Kleisli category $\cat{C}_T$ is a weak pullback if and
	only if its image under the right-adjoint $R:\cat{C}_T\to\cat{C}$ is.
\end{proposition}
\begin{proof}
	Consider the bijection given by the Kleisli adjunction,
	\[
		\begin{tikzcd}[row sep=0]
			\cat{C}_T(X,A) = \cat{C}_T(LX,A) \ar{r}{\cong} & \cat{C}(X,RA) = \cat{C}(X,TA) \\
			f \ar[mapsto]{r} & f^\#
		\end{tikzcd}
	\]
	where $L$ and $R$ denote the left- and right-adjoints. (Note that on
	objects, $LX=X$ and $RA=TA$.) By naturality of the bijection above in
	$A$ (against morphisms of $\cat{C}_T$), we have that for all $p:X\to
		A$, $f:A\to B$ and $q:X\to B$ of $\cat{C}_T$, the triangle of
	$\cat{C}_T$ on the left commutes if and only if the triangle of
	$\cat{C}$ on the right does.
	\begin{equation}\label{eq:triangles}
		\begin{tikzcd}[row sep=small]
			& A \ar{dd}{f} \\
			X \ar{ur}{p} \ar{dr}[swap]{q} \\
			& B
		\end{tikzcd}
		\qquad\longleftrightarrow\qquad
		\begin{tikzcd}[row sep=small]
			& TA \ar{dd}{Rf} \\
			X \ar{ur}{p^\#} \ar{dr}[swap]{q^\#} \\
			& TB
		\end{tikzcd}
	\end{equation}
	Consider now a commutative square of $\cat{C}_T$ as the one on the left
	below.
	\[
		\begin{tikzcd}
			S \\
			& A & B \\
			& C & D
			\arrow["r"{description}, dashed, from=1-1, to=2-2]
			\arrow["q", curve={height=-18pt}, from=1-1, to=2-3]
			\arrow["p"', curve={height=18pt}, from=1-1, to=3-2]
			\arrow["f", from=2-2, to=2-3]
			\arrow["g"', from=2-2, to=3-2]
			\arrow["m", from=2-3, to=3-3]
			\arrow["n"', from=3-2, to=3-3]
		\end{tikzcd}
		\qquad\longleftrightarrow\qquad
		\begin{tikzcd}
			S \\
			& TA & TB \\
			& TC & TD
			\arrow["r^\#"{description}, dashed, from=1-1, to=2-2]
			\arrow["q^\#", curve={height=-18pt}, from=1-1, to=2-3]
			\arrow["p^\#"', curve={height=18pt}, from=1-1, to=3-2]
			\arrow["Rf", from=2-2, to=2-3]
			\arrow["Rg"', from=2-2, to=3-2]
			\arrow["Rm", from=2-3, to=3-3]
			\arrow["Rn"', from=3-2, to=3-3]
		\end{tikzcd}
	\]
	We have to prove that the square on the left is a weak pullback in
	$\cat{C}_T$ if and only if the one on the right is a weak pullback in
	$\cat{C}$.
	
	Let $S$ be an object of $\cat{C}$ (equivalently, of $\cat{C}_T$), and
	let $p:S\to C$ and $q:S\to B$ be morphisms of $\cat{C}_T$ such that
	$m\circ q=n\circ p$ (equivalently, by \eqref{eq:triangles}, let
	$p^\#:S\to TC$ and $q^\#:S\to TB$ be morphisms of $\cat{C}$ such that
	$Rm\circ q^\#=Rn\circ p^\#$). Again by \eqref{eq:triangles}, there
	exists $r:S\to A$ such that $f\circ r=q$ and $g\circ r= p$ if and only
	if there exists $r^\#:S\to TA$ such that $Rf\circ r^\#=q^\#$ and
	$Rg\circ r^\# = p^\#$.
\end{proof}

\begin{corollary}\label{cor:BC_Kleisli}
	Let $(T, \mu, \eta)$ be a monad on a category $\cat{C}$. Then
	\begin{enumerate}
		\item $T$ is weakly cartesian if and only if the left-adjoint $\cat{C}\to\cat{C}_T$ preserves weak pullbacks;
		\item $\mu$ is weakly cartesian if and only if the naturality diagram of the counit $\varepsilon$ in $\cat{C}_T$
		      \[\begin{tikzcd}
				      TX & TY \\
				      X & Y
				      \arrow["Tf", from=1-1, to=1-2]
				      \arrow["{\varepsilon_X}"', from=1-1, to=2-1]
				      \arrow["{\varepsilon_Y}", from=1-2, to=2-2]
				      \arrow["f"', from=2-1, to=2-2]
			      \end{tikzcd}\]
		      is a weak pullback for all $f:X \to Y$ of $\cat{C}$, recalling that $\mu=R\varepsilon$. (In Markov categories, $\varepsilon$ is the map $\samp$.)
	\end{enumerate}
\end{corollary}

Let us now look at partial evaluations from the point of view of representable Markov categories. (See \cite{Fritzetal-2023-Representablemarkovc} for all the details.)

Consider a $P$-algebra $e:PA \to A$ in $\cat{C}_{\text{det}}$. Using the Kleisli bijection
\begin{equation*}
	\cat{C}_{\text{det}}(\Theta, PA) \;\cong\; \cat{C}(\Theta, A)
\end{equation*}
we can compare the partial evaluation order on the left with a natural order on the right, which we define first.

\begin{definition}
	\label{def:second-order-dominance}
	Let $(A, e)$ be a $P$-algebra in $\cat{C}_{\text{det}}$, and let $\Theta$ be an object of $\cat{C}$. The \textbf{second-order dominance relation} on $\cat{C}(\Theta, A)$ is defined as follows: for $p, q : \Theta \to A$, we write $p \leq q$ if there exists $k : \Theta \to PA$ such that the diagram
	\begin{equation*}
		\begin{tikzcd}
			& \Theta \\
			A & PA & A
			\arrow["p"', from=1-2, to=2-1]
			\arrow["k", from=1-2, to=2-2]
			\arrow["q", from=1-2, to=2-3]
			\arrow["e"', from=2-2, to=2-3]
			\arrow["{\samp_A}", from=2-2, to=2-1]
		\end{tikzcd}
	\end{equation*}
	commutes in $\cat{C}$, i.e.\ $\samp_A \circ k = p$ and $e \circ k = q$.
\end{definition}

\begin{proposition}
	\label{prop:monotone}
	The Kleisli bijection
	\begin{equation*}
		\cat{C}_{\text{det}}(\Theta, PA) \;\cong\; \cat{C}(\Theta, A)
	\end{equation*}
	is monotone in both directions, where the order on $\cat{C}_{\text{det}}(\Theta, PA)$ is the partial evaluation order and the order on $\cat{C}(\Theta, A)$ is the second-order dominance relation of \Cref{def:second-order-dominance}.
\end{proposition}
\begin{proof}
	Consider $p,q:\Theta\to A$ in $\cat{C}$, equivalently $p^\#, q^\# : \Theta \to PA$ in $\cat{C}_{\text{det}}$. As in the proof of \Cref{prop:BC_Kleisli}, there exists $k^\# : \Theta \to PPA$ such that $\mu \circ k^\# = p^\#$ and $Pe \circ k^\# = q^\#$ if and only if there exists $k : \Theta \to PA$ such that $\samp_A \circ k = p$ and $e \circ k = q$.
\end{proof}

Setting $\Theta = I$ in $\cat{BorelStoch}$ one recovers usual notions of second-order stochastic dominance for random variables.
Indeed, when the $P$-algebra is the line segment $[0, 1]$ (the free $P$-algebra on the two-element set), we can interpret a partial evaluation as ``subdividing a probability measure and replacing the parts by their centers of mass''.
For free algebras of the form $A = PX$, whose elements are themselves probability measures on $X$, a partial evaluation is instead a coarse-graining of a mixture of probability measures into a coarser mixture.
Recall for example that a mixture of Gaussians is a probability measure that looks (for instance) as follows, and a random variable following such a distribution is called a mixture model:
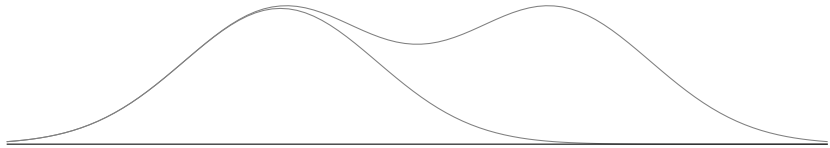
\begin{center}
	\begin{tikzpicture}[x={(0:1.8cm)},y={(90:1.8cm)},z={(180:0.5cm)},%
			bullet/.style={circle, fill, minimum size=2pt, inner sep=1.5pt, outer sep=0pt},%
		]
		\draw [color=gray,thick,domain=-3.5:3.5,samples=120,name path=normal] plot (\x,{exp(-(\x+1.5)*(\x+1.5)/1.28)+exp(-(\x-1.5)*(\x-1.5)/1.28)}) ;
		\draw [color=gray!60,dashed,domain=-3.5:3.5,samples=120] plot (\x,{exp(-(\x+1.5)*(\x+1.5)/1.28)}) ;
		\draw [color=gray!60,dashed,domain=-3.5:3.5,samples=120] plot (\x,{exp(-(\x-1.5)*(\x-1.5)/1.28)}) ;
		\draw (-3.5,0) -- (3.5,0) ;
	\end{tikzpicture}
\end{center}
Partial evaluations in this case can be considered ``coarse-grainings'' of these mixtures of measures:
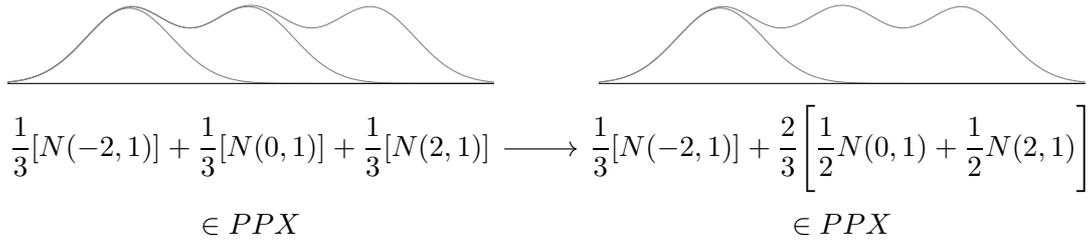
\begin{center}
	\resizebox{\linewidth}{!}{%
	\begin{tikzcd}[row sep=0,ampersand replacement=\&]
		\begin{tikzpicture}[x={(0:0.8cm)},y={(90:3cm)},z={(180:0.5cm)},%
				bullet/.style={circle, fill, minimum size=2pt, inner sep=1.5pt, outer sep=0pt},%
			]
			\draw [color=gray,domain=-4.5:4.5,samples=120] plot (\x,{(1/3)*exp(-(\x+2)*(\x+2)/2)}) ;
			\draw [color=gray,domain=-4.5:4.5,samples=120] plot (\x,{(1/3)*exp(-(\x)*(\x)/2)}) ;
			\draw [color=gray,domain=-4.5:4.5,samples=120] plot (\x,{(1/3)*exp(-(\x-2)*(\x-2)/2)}) ;
			\draw (-4.5,0) -- (4.5,0) ;
		\end{tikzpicture}
		\&
		\begin{tikzpicture}[x={(0:0.8cm)},y={(90:3cm)},z={(180:0.5cm)},%
				bullet/.style={circle, fill, minimum size=2pt, inner sep=1.5pt, outer sep=0pt},%
			]
			\draw [color=gray,domain=-4.5:4.5,samples=120] plot (\x,{(1/3)*exp(-(\x+2)*(\x+2)/2)}) ;
			\draw [color=gray,domain=-4.5:4.5,samples=120] plot (\x,{(1/3)*exp(-(\x)*(\x)/2)+(1/3)*exp(-(\x-2)*(\x-2)/2)}) ;
			\draw (-4.5,0) -- (4.5,0) ;
		\end{tikzpicture}
		\\
		\dfrac{1}{3} [N(-2,1)] + \dfrac{1}{3} [N(0,1)] + \dfrac{1}{3} [N(2,1)]
		\ar{r}
		\& 
		\dfrac{1}{3} [N(-2,1)] + \dfrac{2}{3} \Bigg[ \dfrac{1}{2} N(0,1) + \dfrac{1}{2} N(2,1) \Bigg]
		\\
		\in PPX \& \in PPX
	\end{tikzcd}%
	}
\end{center}
In particular, given a morphism $p:S\to PX$ (or a Kleisli morphism $p^\flat:S\to X$), we can view a morphism $d:S\to PPX$ such that $\mu\circ d=p$ (or equivalently a Kleisli morphism $d^\flat:S\to PX$ such that $\samp\circ d^\flat=p^\flat$) as a ``decomposition'' of $p$.

For more details we refer the reader to \cite{Fritzetal-2023-Representablemarkovc} and also \cite[Chapter 4]{Perrone18}.

\subsection{Statistical experiments and hypernormalizations}\label{stat_exp}

We can model a probability space in a Markov category as a pair $(\Theta,p)$ where $\Theta$ is an object on our category. and $p:I\to\Theta$.
A \textbf{statistical experiment} on $(\Theta,p)$ is a morphism $f:\Theta\to X$ up to $p$-almost sure equality. The idea is that
\begin{itemize}
	\item The space $\Theta$ contains states of the world that we cannot access or observe directly;
	\item The probabilistic state $p:I\to\Theta$ represents our incomplete knowledge on where we are in $\Theta$. From the Bayesian point of view, it is our prior distribution;
	\item The space $X$ encodes something we can observe, and is hence more ``coarse-grained'' than $\Theta$;
	\item The map $f:\Theta\to X$ performs this coarse-graining, or this observation, possibly in a noisy way, so that we allow it to be non-deterministic. (In $\cat{BorelStoch}$, this is a Markov kernel or a stochastic map.)
\end{itemize}

\begin{example}\label{ex:gaussian}
	Consider the following situation in $\cat{BorelStoch}$, where $(\Theta,p)$ is the real line with a normal distribution, and $X=\{1,2,3\}$ is a finite set.
	We now take a deterministic experiment $f:\Theta\to X$ which partitions $\Theta$ into regions labeled by the elements of $X$.
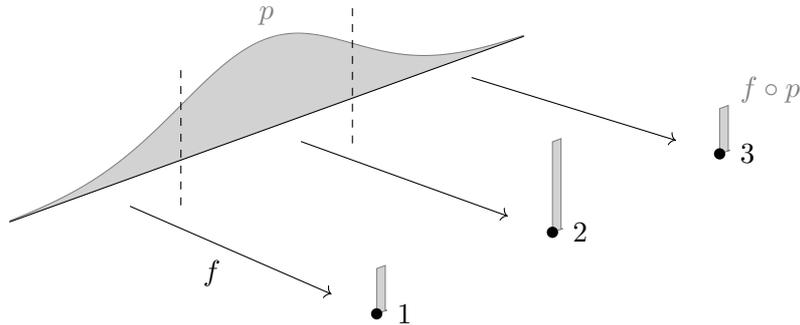
	\begin{center}
		\begin{tikzpicture}[x={(20:1.2cm)},y={(90:1.2cm)},z={(-20:4cm)},%
				bullet/.style={circle, fill, minimum size=2pt, inner sep=1.5pt, outer sep=0pt},%
			]
			\fill [gray!35, domain=-3:3, variable=\x]
			(-3, 0)
			-- plot (\x,{exp(-\x*\x/2)})
			-- (3, 0)
			-- cycle ;
			\draw [color=gray,domain=-3:3,samples=100] plot (\x,{exp(-\x*\x/2)}) node[midway,yshift=1.5cm] {$p$};
			\draw (-3,0) -- (3,0) ;
			
			\draw [black,dashed] (-1,-0.5) -- (-1,1) ;
			\draw [black,dashed] (1,-0.5) -- (1,1) ;
			
			\node [label=right:$1$] (1) at (-2.05,-0.2,1) {};
			\node [label=right:$2$] (2) at (0,0,1) {};
			\node [label=right:$3$] (3) at (1.95,0.2,1) {};
			
			\draw [->, shorten <=5mm, shorten >=5mm] (-2,0,0) -- (1) node[midway,below left] {$f$};
			\draw [->, shorten <=4.8mm, shorten >=4.8mm] (0,0,0) -- (2) ;
			\draw [->, shorten <=4.6mm, shorten >=4.6mm] (2,0,0) -- (3) ;
			
			\draw [draw=gray,fill=gray!35] (1) -- ++(0.1,0) -- ++(0,0.5) -- ++(-0.1,0) -- ++(0,-0.5) --cycle ;
			\draw [draw=gray,fill=gray!35] (2) -- ++(0.1,0) -- ++(0,1) -- ++(-0.1,0) -- ++(0,-1) --cycle ;
			\draw [draw=gray,fill=gray!35] (3) -- ++(0.1,0) -- ++(0,0.5) -- ++(-0.1,0) -- ++(0,-0.5) --cycle node[above right,yshift=0.5cm,color=gray] {$f\circ p$};
			
			\node [bullet] at (1) {};
			\node [bullet] at (2) {};
			\node [bullet] at (3) {};
		\end{tikzpicture}
	\end{center}
	We can interpret it in a Bayesian way as follows: we have a point on the real line, but we do not know where exactly, and the measure $p$ represent our belief about its position. $f$ is now an experiment that will tell us, deterministically, in which one of the three regions the point is, therefore improving our knowledge.
	
	In the picture above we have also drawn the pushforward measure $f\circ p$, which is the probability, according to our belief, that each point is chosen. (It is the same number as the measure, according to $p$, of the corresponding region of $\Theta$.)
	
	Note that, in $\cat{Stoch}$, $X$ is isomorphic to the real line equipped with the sigma-algebra induced by the partition (\cite[Proposition~2.7]{ensarguet2023ergodic}). The same can be said in the category of \emph{probability} spaces (see the same reference).
	We can therefore see a deterministic experiment equivalently as a coarse-graining of the sigma-algebra.
\end{example}

\begin{example}
	Similarly to the example above, we can view a non-deterministic
	experiment as a situation where, instead of partitioning the
	\emph{domain}, we more generally partition the \emph{mass} over the
	domain:
	
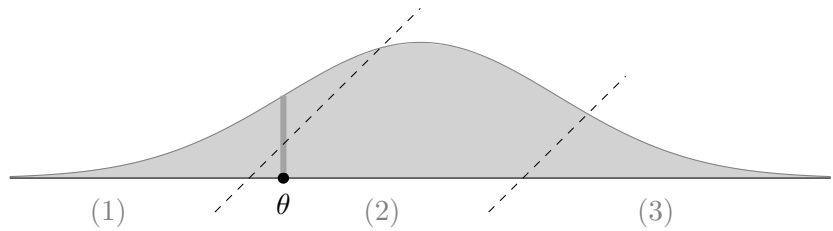
\begin{center}
		\begin{tikzpicture}[x={(0:1.8cm)},y={(90:1.8cm)},z={(180:0.5cm)},%
				bullet/.style={circle, fill, minimum size=2pt, inner sep=1.5pt, outer sep=0pt},%
			]
			\fill [gray!35, domain=-3:3, variable=\x]
			(-3, 0)
			-- plot (\x,{exp(-\x*\x/2)})
			-- (3, 0)
			-- cycle;
			\draw [color=gray,domain=-3:3,samples=100,name path=normal] plot (\x,{exp(-\x*\x/2)}) ;
			\draw (-3,0) -- (3,0) ;
			
			\node [bullet,label=below:$\theta$] (theta) at (-1,0) {};
			\draw [draw=none,name path=up] (theta) -- ++(0,1) ;
			\draw [line width=0.8mm,gray!72,name intersections={of=normal and up,by={i}}] (theta) -- (i);
			
			\draw [black,dashed] (-1.5,-0.25) -- (0.,1.25) ;
			\draw [black,dashed] (0.5,-0.25) -- (1.5,0.75) ;
			
			\node [label={[gray]below:$(1)$}] (1) at (-2,0,1) {};
			\node [label={[gray]below:$(2)$}] (2) at (0,0,1) {};
			\node [label={[gray]below:$(3)$}] (3) at (2,0,1) {};
		\end{tikzpicture}
	\end{center}
	The map $f$ this time takes each point $\theta$ and splits its output across the elements of $X$,
	with the probability of each output being proportional to the height of
	that element's region in the vertical column above $\theta$.
\end{example}

Now following again \Cref{ex:gaussian}, suppose that the result of the
experiment is the point $1$. Then we know that our point lies in the left
region of the curve. But we still do not know where exactly, within that region.
We can therefore update our belief to the following distribution,
$p(\theta|x=1)$:
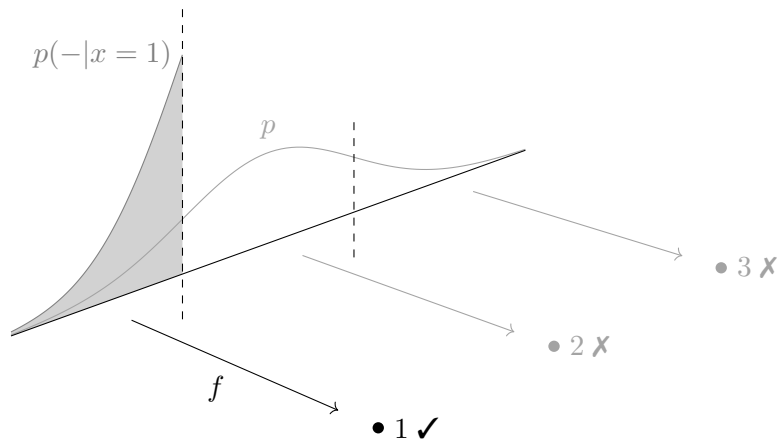
\begin{center}
	\begin{tikzpicture}[x={(20:1.2cm)},y={(90:1.2cm)},z={(-20:4cm)},%
			bullet/.style={circle, fill, minimum size=2pt, inner sep=1.5pt, outer sep=0pt},%
		]
		\fill [gray!35, domain=-3:-1, variable=\x]
		(-3, 0)
		-- plot (\x,{4*exp(-\x*\x/2)})
		-- (-1, 0)
		-- cycle ;
		\draw [color=gray!72,domain=-3:3,samples=100] plot (\x,{exp(-\x*\x/2)}) node[midway,yshift=1.5cm] {$p$};
		\draw [color=gray,domain=-3:-1,samples=100] plot (\x,{4*exp(-\x*\x/2)}) node[left] {$p(-|x=1)$};
		\draw (-3,0) -- (3,0) ;
		
		\draw [black,dashed] (-1,-0.5) -- (-1,3) ;
		\draw [black,dashed] (1,-0.5) -- (1,1) ;
		
		\node [bullet,label=right:$1\;$\cmark] (1) at (-2.05,-0.2,1) {};
		\node [gray!72,bullet,label={[gray!72]right:$2\;$\xmark}] (2) at (0,0,1) {};
		\node [gray!72,bullet,label={[gray!72]right:$3\;$\xmark}] (3) at (1.95,0.2,1) {};
		
		\draw [->, shorten <=5mm, shorten >=5mm] (-2,0,0) -- (1) node[midway,below left] {$f$};
		\draw [gray!72,->, shorten <=4.8mm, shorten >=4.8mm] (0,0,0) -- (2) ;
		\draw [gray!72,->, shorten <=4.6mm, shorten >=4.6mm] (2,0,0) -- (3) ;
	\end{tikzpicture}
\end{center}
\begin{itemize}
	\item The excluded regions $(2)$ and $(3)$ have now probability zero;
	\item The confirmed region $(1)$ has now probability one;
	\item Within region $(1)$, and only within there, we keep the probability proportional to the old measure $p$ -- but we need to renormalize it so that the total mass is one:
	      \[
		      p(A|x=1) \;=\; \dfrac{p(A\cap f^{-1}(1))}{p(f^{-1}(1))} .
	      \]
\end{itemize}
(A similar measure can also be constructed if the experiment is noisy.)

The measure $p(-|x=1)$ is called the \textbf{posterior distribution}, and the
process of moving from $p$ to $p(-|x=1)$ is called \textbf{Bayesian updating}.
Depending on the observed $x$, we get a different posterior distribution
$p(-|x)$. We can then form a probabilistic mapping $X\to\Theta$ assigning to
each $x$ its corresponding posterior distribution (or at least, a random
element distributed according to it). This morphism, when it exists, is called
a Bayesian inverse of $f$, and it is a special case of a conditional:

\begin{definition}
	Let $p:I \to \Theta $ and $ f: \Theta \to X $ be a statistical experiment in a Markov category. A
	\textbf{Bayesian inverse} of $ f $ with respect to $ p $ is a morphism $
		f^\dag_p: X \to \Theta $ such that the following equation holds,
	\begin{equation}
		\label{eq:bayesian0}
		\begin{tikzpicture}
	\begin{pgfonlayer}{nodelayer}
		\node [style=bn] (0) at (3, -1) {};
		\node [style=none] (3) at (4, 1.75) {};
		\node [style=none] (4) at (4, 2.25) {$\Theta$};
		\node [style=none] (5) at (2, 0) {};
		\node [style=none] (6) at (2, 2.25) {$X$};
		\node [style=bn] (8) at (-3, -1) {};
		\node [style=morphism] (9) at (-4, 0.5) {$f$};
		\node [style=none] (10) at (-4, 1.75) {};
		\node [style=none] (11) at (-4, 2.25) {$X$};
		\node [style=none] (12) at (-2, 0) {};
		\node [style=none] (13) at (-2, 2.25) {$\Theta$};
		\node [style=none] (15) at (0, -0.25) {$=$};
		\node [style=state] (16) at (3, -2) {$q$};
		\node [style=none] (21) at (4, 0) {};
		\node [style=medium box] (22) at (4, 0.5) {$f^\dagger_p$};
		\node [style=state] (25) at (-3, -2) {$p$};
		\node [style=none] (26) at (-2, 1.75) {};
		\node [style=none] (27) at (2, 1.75) {};
		\node [style=none] (29) at (-4, 0) {};
	\end{pgfonlayer}
	\begin{pgfonlayer}{edgelayer}
		\draw [in=-90, out=165] (0) to (5.center);
		\draw [in=-90, out=15] (8) to (12.center);
		\draw (9) to (10.center);
		\draw (16) to (0);
		\draw [in=15, out=-90] (21.center) to (0);
		\draw (3.center) to (22);
		\draw (8) to (25);
		\draw (12.center) to (26.center);
		\draw (27.center) to (5.center);
		\draw [in=165, out=-90] (29.center) to (8);
		\draw (29.center) to (9);
		\draw (21.center) to (22);
	\end{pgfonlayer}
\end{tikzpicture}
	\end{equation}
	where $ q = p \circ f $.
	
	More generally, when $p$ depends on a parameter $A$, a
	\textbf{Bayesian inverse} of $ f $ with respect to $ p:A\to\Theta $ is a morphism $
		f^\dag_p: A \otimes X \to \Theta $ such that the following equation holds.
	\begin{equation}
		\label{eq:bayesian}
		\begin{tikzpicture}
	\begin{pgfonlayer}{nodelayer}
		\node [style=bn] (0) at (3, 0) {};
		\node [style=none] (3) at (4.5, 2.75) {};
		\node [style=none] (4) at (4.5, 3.25) {$\Theta$};
		\node [style=none] (5) at (2, 1) {};
		\node [style=none] (6) at (2, 3.25) {$X$};
		\node [style=bn] (8) at (-3, 0) {};
		\node [style=morphism] (9) at (-4, 1.5) {$f$};
		\node [style=none] (10) at (-4, 2.75) {};
		\node [style=none] (11) at (-4, 3.25) {$X$};
		\node [style=none] (12) at (-2, 1) {};
		\node [style=none] (13) at (-2, 3.25) {$\Theta$};
		\node [style=none] (15) at (0, -0.25) {$=$};
		\node [style=morphism] (16) at (3, -1) {$q$};
		\node [style=none] (17) at (4, -3) {};
		\node [style=none] (18) at (4, -3.5) {$A$};
		\node [style=bn] (19) at (4, -2.25) {};
		\node [style=none] (20) at (5, -1.25) {};
		\node [style=none] (21) at (4, 1) {};
		\node [style=medium box] (22) at (4.5, 1.5) {$f^\dagger_p$};
		\node [style=none] (23) at (-3, -3) {};
		\node [style=none] (24) at (-3, -3.5) {$A$};
		\node [style=morphism] (25) at (-3, -1.5) {$p$};
		\node [style=none] (26) at (-2, 2.75) {};
		\node [style=none] (27) at (2, 2.75) {};
		\node [style=none] (28) at (5, 1) {};
		\node [style=none] (29) at (-4, 1) {};
		\node [style=none] (30) at (3, -1.25) {};
	\end{pgfonlayer}
	\begin{pgfonlayer}{edgelayer}
		\draw [in=-90, out=165] (0) to (5.center);
		\draw [in=-90, out=15] (8) to (12.center);
		\draw (9) to (10.center);
		\draw (16) to (0);
		\draw (19) to (17.center);
		\draw [in=15, out=-90] (20.center) to (19);
		\draw [in=15, out=-90] (21.center) to (0);
		\draw (3.center) to (22);
		\draw (25) to (23.center);
		\draw (8) to (25);
		\draw (12.center) to (26.center);
		\draw (27.center) to (5.center);
		\draw (28.center) to (20.center);
		\draw [in=165, out=-90] (29.center) to (8);
		\draw (29.center) to (9);
		\draw [in=165, out=-90] (30.center) to (19);
		\draw (16) to (30.center);
	\end{pgfonlayer}
\end{tikzpicture}
	\end{equation}\\
	where again $ q = p \circ f $.
	
	(When $p$ is clear from the context we will just write $f^\dag$.)
\end{definition}

Note that a Bayesian inverse of $f$ only depends on the a.s.\ equality class of
$f$, and if it exists, it is also unique almost surely.

In a representable Markov category we have access to the actual distribution
(via distribution objects), and so we can interpret the morphism
$(f^\dag_p)^\#:X\to P\Theta$ as assigning from a point $x$ its corresponding
posterior distribution (as a point of $P\Theta$). Mind the difference between:
\begin{itemize}
	\item The Bayesian inverse $f^\dag_p:X\to \Theta$: is a stochastic map,
	      which from $X$ gives us a probabilistic output distributed
	      according to the posterior distribution;
	\item The deterministic morphism $(f^\dag_p)^\#:X\to P\Theta$, which
	      from $X$ gives us the actual posterior distribution
	      (deterministically depending on $X$).
\end{itemize}
As usual, we can obtain the former by applying the map $\samp$ to the latter.
(But the former is also defined outside the representable case.)

As each posterior distribution can be seen as a (renormalized) ``piece'' of the
prior, conversely we can view the prior as a mixture of the posteriors. The
measure on $P\Theta$ which gives the mixing is called the \emph{standard
	measure}:

\begin{definition}
	Let $\cat{C}$ be a representable Markov category. Let $(\Theta,p)$ be a
	probability space, and let $f:\Theta\to X$ be a statistical experiment
	The \textbf{standard measure} of $f$ is the state $\hat{f}_p$ on
	$P\Theta$ given by
	\[
		\begin{tikzcd}
			I \ar{r}{p} & \Theta \ar{r}{f} & X \ar{r}{(f^\dag_p)^\#} & P\Theta .
		\end{tikzcd}
	\]
\end{definition}

\begin{remark}\label{rem:stdmeas-welldef}
	Although the Bayesian inverse $f^\dag_p$ is unique only up to
	$(f \circ p)$-almost sure
	equality~\cite[Proposition~13.6]{Fritz-2020-Asyntheticapproachto},
	the standard measure $\hat{f}_p$ is strictly well-defined,
	i.e.\ independent of the choice of Bayesian inverse.
	Indeed, let $f^\dag_1$ and $f^\dag_2$ be two Bayesian inverses of
	$f$ with respect to~$p$, so that
	$f^\dag_1 =_{(f \circ p)\text{-a.s.}} f^\dag_2$.
	By a.s.-compatible
	representability~\cite[Definition~3.21]{Fritzetal-2023-Representablemarkovc},
	the adjuncts satisfy
	$(f^\dag_1)^\# =_{(f \circ p)\text{-a.s.}} (f^\dag_2)^\#$
	as morphisms $X \to P\Theta$.
	Now, for any $g,h \colon X \to Y$ and any state $q \colon I \to X$,
	the relation $g =_{q\text{-a.s.}} h$ implies $g \circ q = h \circ q$
	strictly, since applying $\texttt{del}_X \otimes \mathrm{id}_Y$ to the
	defining equation of a.s.-equality
	(\Cref{def:almost}) and using the counit axiom
	of the copy map yields $g \circ q = h \circ q$.
	Setting $q = f \circ p$, $g = (f^\dag_1)^\#$, and
	$h = (f^\dag_2)^\#$ gives
	$(f^\dag_1)^\# \circ f \circ p = (f^\dag_2)^\# \circ f \circ p$,
	which is $\hat{f}_p^{(1)} = \hat{f}_p^{(2)}$.

	Note that the standard experiment
	$\hat{f} = (f^\dag_p)^\# \circ f \colon \Theta \to P\Theta$ is
	well-defined only up to $p$-a.s.\
	equality~\cite[{\S}5.3]{Fritzetal-2023-Representablemarkovc};
	the upgrade to strict equality for the standard measure uses
	the fact that $\hat{f}_p$ is a state (source~$I$).
	In particular, the statements of \Cref{blackwell2}, \Cref{stdmeas},
	and \Cref{thm:uni-std} involve only standard measures and are
	therefore unambiguous.
\end{remark}

For discrete $\Theta$ and deterministic $f$, we can equivalently view the
standard measure as follows:
\begin{itemize}
	\item $f$ is equivalently a (finite) partition of $\Theta$;
	\item The standard measure is a measure on $P\Theta$ such that each of
	      the elements (i.e.~measures) on its support are each supported
	      exactly on a cell of the partition. 
\end{itemize}
Such a construction, for the special case of discrete product projections was
called ``hypernormalization'' (up to different spelling conventions) in
\cite{hypernormalization}. Let us give here a partial generalization, with a
similar intuition.

\begin{definition}\label{def:hypernormalization}
	Let $(\Theta,p)$ be a probability space in a representable Markov
	category, and let $f:\Theta\to X$ be an a.s.~deterministic statistical
	experiment. The \textbf{hypernormalization} of $p$ with respect to $f$,
	if it exists, is the standard measure $\hat{f}_p$ on $P\Theta$. 
\end{definition}

The construction given in \cite{hypernormalization}, or at least its
instantiation for the case of the distribution monad, can be considered a
special case of the definition above. Let us see this:
\begin{example}
	In $\cat{Stoch}$, let us consider the special case of finite sets $X$
	and $Y$ with a probability distribution $p$ on $X\times Y$. Let now
	$\pi:X\times Y\to X$ be the product projection, which we can see as
	partitioning the product into $X$-many copies of $Y$, the fibers
	$\pi^{-1}(x)\subseteq X\times Y$. Instantiating our notion of
	hypernormalization, we are decomposing $p$ into a convex combination of
	measures $p_x\in P(X\times Y)$ as follows,
	\[
		p= \sum_{x\in X} q(x) \, p_x \qquad\Longleftrightarrow\qquad p(x',y)= \sum_{x\in X} q(x) \, p_x(x',y)
	\]
	where 
	\begin{itemize}
		\item $q$ is the marginal of $p$ on $X$;
		\item $p_x$ is really only appearing for those $x$ for which $q(x)>0$;
		\item Each measure $p_x$ (where $q(x)>0$) is supported on the fiber $\pi^{-1}(x)$;
		\item Each measure $p_x$ (where $q(x)>0$) is proportional to the restriction of $p$ to the fiber:
		      \[
			      p_x(x',y) \;=\;\begin{cases}
				      \frac{p(x',y)}{q(x)} & x'=x,      \\
				      0                    & x'\ne x .
			      \end{cases}
		      \]
	\end{itemize}
	This is equivalent to the construction given in
	\cite[Section~3]{hypernormalization} for the distribution monad.
\end{example}

Note that, according to our definitions, in order for the hypernormalization to
exist, we need the Bayesian inverse $f^\dag_p$ to exist. (In $\cat{Stoch}$,
this is guaranteed if $X$ is standard Borel.) We will see in \Cref{hypernorm}
that hypernormalizations satisfy a universal property. Therefore, one could
define them more generally outside the case where the necessary conditionals
exist (for example, in categories of topological spaces).

Let us now define an order on statistical experiments, deterministic or not.
Experiments are naturally ordered in terms of ``how informative they are'':

\begin{definition}
	Let $f:\Theta\to X$ and $g:\Theta\to Y$ be statistical experiments on
	$(\Theta,p)$. We say that $g\le f$ in the \textbf{Blackwell order} if
	and only if there exists a morphism $h:X\to Y$ such that $h\circ f=_p
		g$.
\end{definition}

In some sense, $f$ is ``more informative'' than $g$ if we can recover the
results of the experiment $g$ by processing the results of $f$ without having
any further access to $\Theta$.

\begin{theorem}[Blackwell-Sherman-Stein theorem for Markov categories, {\cite[Theorem~5.13]{Fritzetal-2023-Representablemarkovc}}]\label{thm:bss}
	Let $\cat{C}$ be an a.s.-compatibly representable Markov category. Let
	$f:\Theta\to X$ and $g:\Theta\to Y$ be statistical experiments on
	$(\Theta,p)$. Then $g\le f$ in the Blackwell order if and only if for
	their standard measures, $\hat{g}_p\ge\hat{f}_p$ in the stochastic
	dominance order of $\cat{C}(I,P\Theta)$.
\end{theorem}

\begin{example}
	As in \Cref{ex:gaussian}, let $(\Theta,p)$ be the real line with the
	normal distribution, and let $f:\Theta\to X$ with $X=\{1,2,3\}$.
	Consider now $Y=\{a,b\}$, and the (deterministic) function $h:X\to Y$
	given by $h(1)=a$, $h(2)=h(3)=b$. Then the function $g=h\circ
		f:\Theta\to Y$ can be seen as partitioning $\Theta$, but putting
	together the second and third region:
		\begin{center}
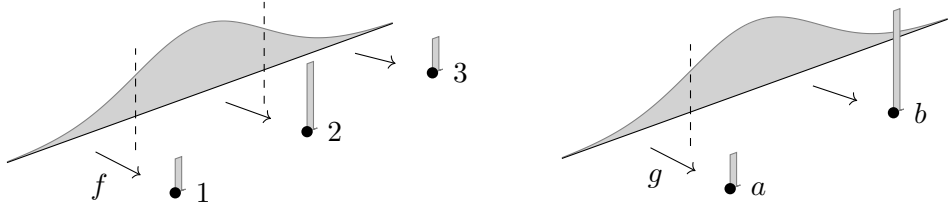
			\begin{minipage}{0.45\textwidth}
				\begin{center}
					\resizebox{\linewidth}{!}{%
					\begin{tikzpicture}[x={(20:0.9cm)},y={(90:0.9cm)},z={(-20:1.5cm)},%
						bullet/.style={circle, fill, minimum size=2pt, inner sep=1.5pt, outer sep=0pt},%
					]
					\fill [gray!35, domain=-3:3, variable=\x]
					(-3, 0)
					-- plot (\x,{exp(-\x*\x/2)})
					-- (3, 0)
					-- cycle ;
					\draw [color=gray,domain=-3:3,samples=100] plot (\x,{exp(-\x*\x/2)}) ;
					\draw (-3,0) -- (3,0) ;
					
					\draw [black,dashed] (-1,-0.5) -- (-1,1) ;
					\draw [black,dashed] (1,-0.5) -- (1,1) ;
					
					\node [label=right:$1$] (1) at (-2.05,-0.2,1) {};
					\node [label=right:$2$] (2) at (0,0,1) {};
					\node [label=right:$3$] (3) at (1.95,0.2,1) {};
					
					\draw [->, shorten <=3.6mm, shorten >=3.6mm] (-2,0,0) -- (1) node[midway,below left] {$f$};
					\draw [->, shorten <=3.6mm, shorten >=3.6mm] (0,0,0) -- (2) ;
					\draw [->, shorten <=3.6mm, shorten >=3.6mm] (2,0,0) -- (3) ;
					
					\draw [draw=gray,fill=gray!35] (1) -- ++(0.1,0) -- ++(0,0.5) -- ++(-0.1,0) -- ++(0,-0.5) --cycle ;
					\draw [draw=gray,fill=gray!35] (2) -- ++(0.1,0) -- ++(0,1) -- ++(-0.1,0) -- ++(0,-1) --cycle ;
					\draw [draw=gray,fill=gray!35] (3) -- ++(0.1,0) -- ++(0,0.5) -- ++(-0.1,0) -- ++(0,-0.5) --cycle ;
					
					\node [bullet] at (1) {};
					\node [bullet] at (2) {};
					\node [bullet] at (3) {};
					\end{tikzpicture}%
					}
				\end{center}
			\end{minipage}
			\begin{minipage}{0.45\textwidth}
				\begin{center}
					\resizebox{\linewidth}{!}{%
					\begin{tikzpicture}[x={(20:0.9cm)},y={(90:0.9cm)},z={(-20:1.5cm)},%
						bullet/.style={circle, fill, minimum size=2pt, inner sep=1.5pt, outer sep=0pt},%
					]
					\fill [gray!35, domain=-3:3, variable=\x]
					(-3, 0)
					-- plot (\x,{exp(-\x*\x/2)})
					-- (3, 0)
					-- cycle ;
					\draw [color=gray,domain=-3:3,samples=100] plot (\x,{exp(-\x*\x/2)}) ;
					\draw (-3,0) -- (3,0) ;
					
					\draw [black,dashed] (-1,-0.5) -- (-1,1) ;
					
					\node [label=right:$a$] (1) at (-2.05,-0.2,1) {};
					\node [label=right:$b$] (2) at (0.49,0.05,1) {};
					
					\draw [->, shorten <=3.6mm, shorten >=3.6mm] (-2,0,0) -- (1) node[midway,below left] {$g$};
					\draw [->, shorten <=3.6mm, shorten >=3.6mm] (0.5,0,0) -- (2) ;
					
					\draw [draw=gray,fill=gray!35] (1) -- ++(0.1,0) -- ++(0,0.5) -- ++(-0.1,0) -- ++(0,-0.5) --cycle ;
					\draw [draw=gray,fill=gray!35] (2) -- ++(0.1,0) -- ++(0,1.5) -- ++(-0.1,0) -- ++(0,-1.5) --cycle ;
					
					\node [bullet] at (1) {};
					\node [bullet] at (2) {};
					\end{tikzpicture}%
					}
				\end{center}
			\end{minipage}
		\end{center}
	Therefore the decomposition of $p$ is coarser than the one given by $f$.
\end{example}

The idea illustrated in this example will be made precise in \Cref{hypernorm}.

We conclude with a sort of converse to \Cref{thm:bss}, which seems to be new.

\begin{proposition}\label{blackwell2}
	Let $\cat{C}$ be an a.s.-compatibly representable Markov category with conditionals.
	Let $\pi,\tau:I\to P\Theta$. Then $\pi\le\tau$ in the order of stochastic dominance on $\cat{C}(I,P\Theta)$ if and only if both of the following conditions are satisfied:
	\begin{enumerate}
		\item $\samp\circ\pi=\samp\circ\tau$. (That is, they are decompositions of the same state, denote it by $p:I\to\Theta$.)
		\item There exist statistical experiments $f:\Theta\to X$ and $g:\Theta\to Y$ on $(\Theta,p)$ with $\hat{f}_p=\pi$, $\hat{g}_p=\tau$, and $f\ge g$ in the Blackwell order.
	\end{enumerate}
\end{proposition}

We will use the following auxiliary statement.
\begin{lemma}\label{stdmeas}
	Let $\pi:I\to P\Theta$, denote $\samp\circ\pi:I\to\Theta$ by $p$, and let $f:=\samp^\dag_\pi:\Theta\to P\Theta$, viewed as a statistical experiment on $(\Theta,p)$.
	Then the standard measure of $f$ is exactly $\pi$.
\end{lemma}
\begin{proof}[Proof of \Cref{stdmeas}]
	Marginalizing the Bayesian inverse equation defining $f=\samp^\dag_\pi$ gives $f\circ p=\pi$.
	By commutativity of copying, the same equation shows that $\samp$ is a Bayesian inverse of $f$ with respect to $p$, so that $f^\dag_p =_{\pi\text{-a.s.}} \samp$.
	By a.s.-compatible representability, $(f^\dag_p)^\sharp =_{\pi\text{-a.s.}} \samp^\sharp$, and composing with the state $\pi$ gives $(f^\dag_p)^\sharp\circ\pi = \samp^\sharp\circ\pi = \pi$ strictly, since $\samp^\sharp = \mathrm{id}_{P\Theta}$.
	The standard measure of $f$ is therefore
	\begin{align*}
		(f^\dag_p)^\sharp \circ f \circ p
		&= (f^\dag_p)^\sharp \circ\pi
		= \pi .
	\end{align*}
\end{proof}

\begin{proof}[Proof of \Cref{blackwell2}]
	One side of the implication is given by \Cref{thm:bss}. More in detail, suppose that $\samp\circ\pi=\samp\circ\tau$, and denote either side of the equation by $p:I\to\Theta$. Suppose now that there exist experiments $f:\Theta\to X$ and $g:\Theta\to Y$ on $(\Theta,p)$ with $\hat{f}_p=\pi$, $\hat{g}_p=\tau$, and $f\ge g$ in the Blackwell order. Then by \Cref{thm:bss}, $\pi\le\tau$.
	
	Conversely, suppose that $\pi\le\tau$ in the order of stochastic dominance. Then by definition there exists $\kappa:I\to PP\Theta$ such that $\samp\circ\kappa=\pi$ and $\mu\circ\kappa=\tau$. 
	By naturality of $\samp$ (and recalling that $\mu=P\samp$), we now have
	\[
		\samp\circ\pi\;=\;\samp\circ\samp\circ\kappa \;=\; \samp\circ\mu\circ\kappa \;=\; \samp\circ\tau .
	\]
	Denote now by $p:I\to\Theta$ either side of the equation above, and
	consider the statistical experiments $f=\samp^\dag_\pi$ and
	$g=\samp^\dag_\tau:\Theta\to P\Theta$ on $(\Theta,p)$. By \Cref{stdmeas},
	their standard measures are respectively $\pi$ and $\tau$. As
	$\pi\le\tau$, by \Cref{thm:bss} we have that $f\ge g$ in the Blackwell
	order.
\end{proof}

\section{Main results}\label{main_results}

\subsection{The multiplication map is weakly cartesian}\label{wcartmu}

\begin{theorem}
	\label{thm:pullbacks}
	Let $\cat{C}$ be an a.s.-compatibly representable Markov category with
	monad $(P,\mu,\delta)$. If $\cat{C}$ has conditionals, then $\mu$ is
	weakly Cartesian.
\end{theorem}

\begin{corollary}\label{pullbacks_giry}
	The Giry monad on standard Borel spaces has weakly cartesian multiplication.
\end{corollary}

In order to prove the theorem we use the following technical statement, which
can be seen as an equivalent characterization of a.s.-compatible
representability for the case when conditionals exist.

\begin{lemma}\label{lem:alt-asc}
	Let $\cat{C}$ be a representable Markov category.
	Then
	\begin{enumerate}
		\item If $p:I\to X$ and $f,g:X\to Y$ such that
		      $Pf=_\pi Pg$ for every $\pi:I\to PX$ with
		      $\samp\circ\pi=p$, then $f=_p g$:
		      \begin{equation}\label{eq:alt-asc0}
			      \begin{tikzpicture}
	\begin{pgfonlayer}{nodelayer}
		\node [style=bn] (1) at (-5, -0.75) {};
		\node [style=state] (6) at (-5, -1.75) {$\pi$};
		\node [style=none] (7) at (-6, 0.25) {};
		\node [style=none] (8) at (-4, 0.25) {};
		\node [style=none] (10) at (-4, 1.75) {};
		\node [style=none] (14) at (-6, 2.25) {$PX$};
		\node [style=none] (15) at (-6, 1.75) {};
		\node [style=none] (16) at (-4, 2.25) {$PY$};
		\node [style=morphism] (17) at (-4, 0.75) {$Pg$};
		\node [style=none] (18) at (0, 0) {$\Longrightarrow$};
		\node [style=bn] (20) at (-11.5, -0.75) {};
		\node [style=state] (25) at (-11.5, -1.75) {$\pi$};
		\node [style=none] (26) at (-12.5, 0.25) {};
		\node [style=none] (27) at (-10.5, 0.25) {};
		\node [style=none] (29) at (-10.5, 1.75) {};
		\node [style=none] (33) at (-12.5, 2.25) {$PX$};
		\node [style=none] (34) at (-12.5, 1.75) {};
		\node [style=none] (35) at (-10.5, 2.25) {$PY$};
		\node [style=morphism] (36) at (-10.5, 0.75) {$Pf$};
		\node [style=none] (37) at (5.75, 0) {$=$};
		\node [style=bn] (39) at (8.25, -0.75) {};
		\node [style=state] (44) at (8.25, -1.75) {$p$};
		\node [style=none] (45) at (7.25, 0.25) {};
		\node [style=none] (46) at (9.25, 0.25) {};
		\node [style=morphism] (50) at (9.25, 0.75) {$g$};
		\node [style=none] (51) at (9.25, 1.75) {};
		\node [style=none] (52) at (7.25, 2.25) {$X$};
		\node [style=none] (53) at (7.25, 1.75) {};
		\node [style=none] (54) at (9.25, 2.25) {$Y$};
		\node [style=bn] (57) at (2.75, -0.75) {};
		\node [style=state] (62) at (2.75, -1.75) {$p$};
		\node [style=none] (63) at (1.75, 0.25) {};
		\node [style=none] (64) at (3.75, 0.25) {};
		\node [style=morphism] (68) at (3.75, 0.75) {$f$};
		\node [style=none] (69) at (3.75, 1.75) {};
		\node [style=none] (70) at (1.75, 2.25) {$X$};
		\node [style=none] (71) at (1.75, 1.75) {};
		\node [style=none] (72) at (3.75, 2.25) {$Y$};
		\node [style=none] (74) at (-8, 0) {$=$};
		\node [style=none] (75) at (-2, 0) {$\forall\pi$};
	\end{pgfonlayer}
	\begin{pgfonlayer}{edgelayer}
		\draw [in=165, out=-90] (7.center) to (1);
		\draw [in=270, out=15] (1) to (8.center);
		\draw (10.center) to (8.center);
		\draw (15.center) to (7.center);
		\draw [in=165, out=-90] (26.center) to (20);
		\draw [in=270, out=15] (20) to (27.center);
		\draw (29.center) to (27.center);
		\draw (34.center) to (26.center);
		\draw [in=165, out=-90] (45.center) to (39);
		\draw [in=270, out=15] (39) to (46.center);
		\draw (51.center) to (50);
		\draw (53.center) to (45.center);
		\draw [in=165, out=-90] (63.center) to (57);
		\draw [in=270, out=15] (57) to (64.center);
		\draw (69.center) to (68);
		\draw (71.center) to (63.center);
		\draw (20) to (25);
		\draw (1) to (6);
		\draw (57) to (62);
		\draw (39) to (44);
		\draw (68) to (64.center);
		\draw (50) to (46.center);
	\end{pgfonlayer}
\end{tikzpicture}
		      \end{equation}
		\item If $p:A\to X$ and $f,g:X\otimes A\to Y$ such that
		      $Pf\circ\sigma_{X,A}=_\pi Pg\circ\sigma_{X,A}$ for
		      every $\pi:A\to PX$ with $\samp\circ\pi=p$, then $f=_p
			      g$: 
		      \begin{equation}\label{eq:alt-asc}
			      \resizebox{0.99\linewidth}{!}{\begin{tikzpicture}
	\begin{pgfonlayer}{nodelayer}
		\node [style=bn] (0) at (-5, -3) {};
		\node [style=bn] (1) at (-6, -0.5) {};
		\node [style=none] (2) at (-5, -4) {};
		\node [style=none] (3) at (-5, -4.5) {$A$};
		\node [style=none] (4) at (-6, -2) {};
		\node [style=none] (5) at (-4, -2) {};
		\node [style=morphism] (6) at (-6, -1.5) {$\pi$};
		\node [style=none] (7) at (-7, 0.5) {};
		\node [style=none] (8) at (-5, 0.5) {};
		\node [style=none] (9) at (-4, 0.5) {};
		\node [style=none] (10) at (-5, 1) {};
		\node [style=none] (11) at (-4, 1) {};
		\node [style=medium box] (12) at (-4.5, 1) {$\sigma_{X,A}$};
		\node [style=none] (13) at (-4.5, 3.75) {};
		\node [style=none] (14) at (-7, 4.25) {$PX$};
		\node [style=none] (15) at (-7, 3.75) {};
		\node [style=none] (16) at (-4.5, 4.25) {$PY$};
		\node [style=morphism] (17) at (-4.5, 2.75) {$Pg$};
		\node [style=none] (18) at (0.25, 0) {$\Longrightarrow$};
		\node [style=bn] (19) at (-11, -3) {};
		\node [style=bn] (20) at (-12, -0.5) {};
		\node [style=none] (21) at (-11, -4) {};
		\node [style=none] (22) at (-11, -4.5) {$A$};
		\node [style=none] (23) at (-12, -2) {};
		\node [style=none] (24) at (-10, -2) {};
		\node [style=morphism] (25) at (-12, -1.5) {$\pi$};
		\node [style=none] (26) at (-13, 0.5) {};
		\node [style=none] (27) at (-11, 0.5) {};
		\node [style=none] (28) at (-10, 0.5) {};
		\node [style=none] (29) at (-11, 1) {};
		\node [style=none] (30) at (-10, 1) {};
		\node [style=medium box] (31) at (-10.5, 1) {$\sigma_{X,A}$};
		\node [style=none] (32) at (-10.5, 3.75) {};
		\node [style=none] (33) at (-13, 4.25) {$PX$};
		\node [style=none] (34) at (-13, 3.75) {};
		\node [style=none] (35) at (-10.5, 4.25) {$PY$};
		\node [style=morphism] (36) at (-10.5, 2.75) {$Pf$};
		\node [style=none] (37) at (-8.25, 0) {$=$};
		\node [style=bn] (38) at (10, -3) {};
		\node [style=bn] (39) at (9, -0.5) {};
		\node [style=none] (40) at (10, -4) {};
		\node [style=none] (41) at (10, -4.5) {$A$};
		\node [style=none] (42) at (9, -2) {};
		\node [style=none] (43) at (11, -2) {};
		\node [style=morphism] (44) at (9, -1.5) {$p$};
		\node [style=none] (45) at (8, 0.5) {};
		\node [style=none] (46) at (10, 0.5) {};
		\node [style=none] (47) at (11, 0.5) {};
		\node [style=none] (48) at (10, 1.75) {};
		\node [style=none] (49) at (11, 1.75) {};
		\node [style=medium box] (50) at (10.5, 1.75) {$g$};
		\node [style=none] (51) at (10.5, 3.75) {};
		\node [style=none] (52) at (8, 4.25) {$X$};
		\node [style=none] (53) at (8, 3.75) {};
		\node [style=none] (54) at (10.5, 4.25) {$Y$};
		\node [style=bn] (56) at (4, -3) {};
		\node [style=bn] (57) at (3, -0.5) {};
		\node [style=none] (58) at (4, -4) {};
		\node [style=none] (59) at (4, -4.5) {$A$};
		\node [style=none] (60) at (3, -2) {};
		\node [style=none] (61) at (5, -2) {};
		\node [style=morphism] (62) at (3, -1.5) {$p$};
		\node [style=none] (63) at (2, 0.5) {};
		\node [style=none] (64) at (4, 0.5) {};
		\node [style=none] (65) at (5, 0.5) {};
		\node [style=none] (66) at (4, 1.75) {};
		\node [style=none] (67) at (5, 1.75) {};
		\node [style=medium box] (68) at (4.5, 1.75) {$f$};
		\node [style=none] (69) at (4.5, 3.75) {};
		\node [style=none] (70) at (2, 4.25) {$X$};
		\node [style=none] (71) at (2, 3.75) {};
		\node [style=none] (72) at (4.5, 4.25) {$Y$};
		\node [style=none] (74) at (6.75, 0) {$=$};
		\node [style=none] (75) at (-2, 0) {$\forall\pi$};
	\end{pgfonlayer}
	\begin{pgfonlayer}{edgelayer}
		\draw (0) to (2.center);
		\draw [in=165, out=-90] (4.center) to (0);
		\draw [in=-90, out=15] (0) to (5.center);
		\draw (1) to (4.center);
		\draw [in=165, out=-90] (7.center) to (1);
		\draw [in=270, out=15] (1) to (8.center);
		\draw (9.center) to (5.center);
		\draw (10.center) to (8.center);
		\draw (11.center) to (9.center);
		\draw (13.center) to (12);
		\draw (15.center) to (7.center);
		\draw (19) to (21.center);
		\draw [in=165, out=-90] (23.center) to (19);
		\draw [in=-90, out=15] (19) to (24.center);
		\draw (20) to (23.center);
		\draw [in=165, out=-90] (26.center) to (20);
		\draw [in=270, out=15] (20) to (27.center);
		\draw (28.center) to (24.center);
		\draw (29.center) to (27.center);
		\draw (30.center) to (28.center);
		\draw (32.center) to (31);
		\draw (34.center) to (26.center);
		\draw (38) to (40.center);
		\draw [in=165, out=-90] (42.center) to (38);
		\draw [in=-90, out=15] (38) to (43.center);
		\draw (39) to (42.center);
		\draw [in=165, out=-90] (45.center) to (39);
		\draw [in=270, out=15] (39) to (46.center);
		\draw (47.center) to (43.center);
		\draw (48.center) to (46.center);
		\draw (49.center) to (47.center);
		\draw (51.center) to (50);
		\draw (53.center) to (45.center);
		\draw (56) to (58.center);
		\draw [in=165, out=-90] (60.center) to (56);
		\draw [in=-90, out=15] (56) to (61.center);
		\draw (57) to (60.center);
		\draw [in=165, out=-90] (63.center) to (57);
		\draw [in=270, out=15] (57) to (64.center);
		\draw (65.center) to (61.center);
		\draw (66.center) to (64.center);
		\draw (67.center) to (65.center);
		\draw (69.center) to (68);
		\draw (71.center) to (63.center);
	\end{pgfonlayer}
\end{tikzpicture}}
		      \end{equation}
		      where $\sigma_{X,A}:PX\otimes A\to P(X\otimes A)$ is
		      the right strength of the monad $P$.
		      
		\item Suppose moreover that $\cat{C}$ has conditionals. Then
		      the converse to the implication above holds if and only
		      if $\cat{C}$ is a.s.-compatibly representable.
	\end{enumerate}
\end{lemma}

Note that the forward direction of parts 1 and~2 only requires the hypothesis for $\pi = \delta\circ p$; the universal quantification is used in the converse direction of part~3.

\begin{proof}[Proof of \Cref{lem:alt-asc}]
	\begin{enumerate}
		\item Note that this statement is implied by 2., setting $A=I$.
		      So let us prove 2.\ directly.

		\item Take $\pi=\delta\circ p$, and notice that
		      $\samp\circ\pi=\samp\circ\delta\circ p= p $. Applying
		      $\samp$ to both outputs of the first term in
		      \eqref{eq:alt-asc} we get:
		      
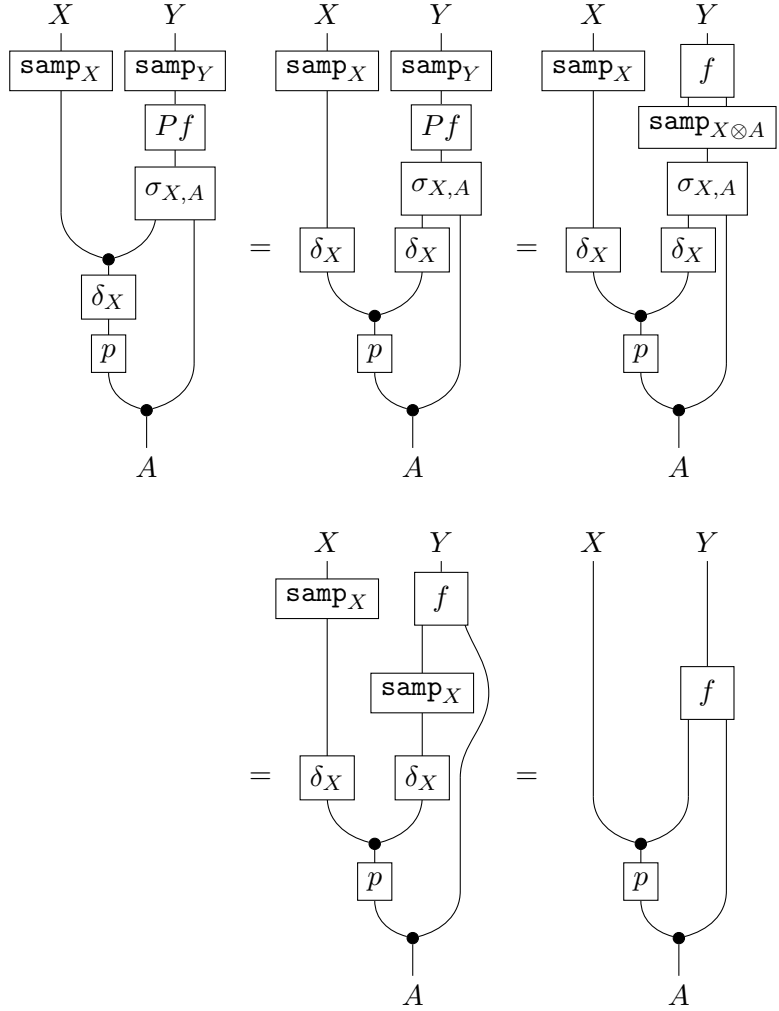
\[
			      \begin{tikzpicture}
	\begin{pgfonlayer}{nodelayer}
		\node [style=bn] (19) at (-3, 3.75) {};
		\node [style=bn] (20) at (-4, 7.75) {};
		\node [style=none] (21) at (-3, 2.75) {};
		\node [style=none] (22) at (-3, 2.25) {$A$};
		\node [style=none] (23) at (-4, 4.75) {};
		\node [style=none] (24) at (-1.75, 5) {};
		\node [style=morphism] (25) at (-4, 5.25) {$p$};
		\node [style=none] (26) at (-5.25, 9) {};
		\node [style=none] (27) at (-2.75, 9) {};
		\node [style=none] (28) at (-1.75, 9) {};
		\node [style=none] (29) at (-2.75, 9.5) {};
		\node [style=none] (30) at (-1.75, 9.5) {};
		\node [style=medium box] (31) at (-2.25, 9.5) {$\sigma_{X,A}$};
		\node [style=none] (32) at (-2.25, 13.75) {};
		\node [style=none] (33) at (-5.25, 14.25) {$X$};
		\node [style=none] (34) at (-5.25, 13.75) {};
		\node [style=none] (35) at (-2.25, 14.25) {$Y$};
		\node [style=morphism] (36) at (-2.25, 11.25) {$Pf$};
		\node [style=none] (37) at (0, 8) {$=$};
		\node [style=morphism] (38) at (-4, 6.75) {$\delta_X$};
		\node [style=morphism] (39) at (-2.25, 12.75) {$\samp_Y$};
		\node [style=morphism] (40) at (-5.25, 12.75) {$\samp_X$};
		\node [style=bn] (41) at (4, 3.75) {};
		\node [style=bn] (42) at (3, 6.25) {};
		\node [style=none] (43) at (4, 2.75) {};
		\node [style=none] (44) at (4, 2.25) {$A$};
		\node [style=none] (45) at (3, 4.75) {};
		\node [style=none] (46) at (5.25, 5) {};
		\node [style=morphism] (47) at (3, 5.25) {$p$};
		\node [style=none] (48) at (1.75, 7.5) {};
		\node [style=none] (49) at (4.25, 7.5) {};
		\node [style=none] (50) at (5.25, 9) {};
		\node [style=none] (51) at (4.25, 9.5) {};
		\node [style=none] (52) at (5.25, 9.5) {};
		\node [style=medium box] (53) at (4.75, 9.625) {$\sigma_{X,A}$};
		\node [style=none] (54) at (4.75, 13.75) {};
		\node [style=none] (55) at (1.75, 14.25) {$X$};
		\node [style=none] (56) at (1.75, 13.75) {};
		\node [style=none] (57) at (4.75, 14.25) {$Y$};
		\node [style=morphism] (58) at (4.75, 11.25) {$Pf$};
		\node [style=morphism] (59) at (1.75, 8) {$\delta_X$};
		\node [style=morphism] (60) at (4.75, 12.75) {$\samp_Y$};
		\node [style=morphism] (61) at (1.75, 12.75) {$\samp_X$};
		\node [style=morphism] (62) at (4.25, 8) {$\delta_X$};
		\node [style=none] (86) at (7, 8) {$=$};
		\node [style=bn] (87) at (11, 3.75) {};
		\node [style=bn] (88) at (10, 6.25) {};
		\node [style=none] (89) at (11, 2.75) {};
		\node [style=none] (90) at (11, 2.25) {$A$};
		\node [style=none] (91) at (10, 4.75) {};
		\node [style=none] (92) at (12.25, 5) {};
		\node [style=morphism] (93) at (10, 5.25) {$p$};
		\node [style=none] (94) at (8.75, 7.5) {};
		\node [style=none] (95) at (11.25, 7.5) {};
		\node [style=none] (96) at (12.25, 9) {};
		\node [style=none] (97) at (11.25, 9.5) {};
		\node [style=none] (98) at (12.25, 9.5) {};
		\node [style=medium box] (99) at (11.75, 9.625) {$\sigma_{X,A}$};
		\node [style=none] (100) at (11.75, 11.25) {};
		\node [style=none] (101) at (8.75, 14.25) {$X$};
		\node [style=none] (102) at (8.75, 13.75) {};
		\node [style=none] (103) at (11.75, 14.25) {$Y$};
		\node [style=medium box] (104) at (11.75, 12.75) {$f$};
		\node [style=morphism] (105) at (11.75, 11.25) {$\samp_{X\otimes A}$};
		\node [style=morphism] (106) at (11.25, 8) {$\delta_X$};
		\node [style=none] (107) at (0, -6) {$=$};
		\node [style=bn] (108) at (4, -10.25) {};
		\node [style=bn] (109) at (3, -7.75) {};
		\node [style=none] (110) at (4, -11.25) {};
		\node [style=none] (111) at (4, -11.75) {$A$};
		\node [style=none] (112) at (3, -9.25) {};
		\node [style=none] (113) at (5.25, -9) {};
		\node [style=morphism] (114) at (3, -8.75) {$p$};
		\node [style=none] (115) at (1.75, -6.5) {};
		\node [style=none] (116) at (4.25, -6.5) {};
		\node [style=none] (117) at (5.25, -6) {};
		\node [style=none] (118) at (4.25, -1.25) {};
		\node [style=none] (121) at (4.75, -0.25) {};
		\node [style=none] (122) at (1.75, 0.25) {$X$};
		\node [style=none] (123) at (1.75, -0.25) {};
		\node [style=none] (124) at (4.75, 0.25) {$Y$};
		\node [style=medium box] (125) at (4.75, -1.25) {$f$};
		\node [style=morphism] (126) at (4.25, -3.75) {$\samp_X$};
		\node [style=morphism] (127) at (4.25, -6) {$\delta_X$};
		\node [style=morphism] (131) at (8.75, 8) {$\delta_X$};
		\node [style=morphism] (132) at (8.75, 12.75) {$\samp_X$};
		\node [style=morphism] (133) at (1.75, -6) {$\delta_X$};
		\node [style=morphism] (134) at (1.75, -1.25) {$\samp_X$};
		\node [style=none] (135) at (11.75, 13.75) {};
		\node [style=none] (136) at (11.25, 12.75) {};
		\node [style=none] (137) at (12.25, 12.75) {};
		\node [style=none] (138) at (11.25, 11.25) {};
		\node [style=none] (139) at (12.25, 11.25) {};
		\node [style=none] (140) at (7, -6) {$=$};
		\node [style=bn] (141) at (11, -10.25) {};
		\node [style=bn] (142) at (10, -7.75) {};
		\node [style=none] (143) at (11, -11.25) {};
		\node [style=none] (144) at (11, -11.75) {$A$};
		\node [style=none] (145) at (10, -9.25) {};
		\node [style=none] (146) at (12.25, -9) {};
		\node [style=morphism] (147) at (10, -8.75) {$p$};
		\node [style=none] (148) at (8.75, -6.5) {};
		\node [style=none] (149) at (11.25, -6.5) {};
		\node [style=none] (150) at (12.25, -3.75) {};
		\node [style=none] (151) at (11.25, -3.75) {};
		\node [style=none] (152) at (11.75, -0.25) {};
		\node [style=none] (153) at (8.75, 0.25) {$X$};
		\node [style=none] (154) at (8.75, -0.25) {};
		\node [style=none] (155) at (11.75, 0.25) {$Y$};
		\node [style=medium box] (156) at (11.75, -3.75) {$f$};
		\node [style=none] (157) at (6, -3.75) {};
		\node [style=none] (158) at (5.25, -1.25) {};
	\end{pgfonlayer}
	\begin{pgfonlayer}{edgelayer}
		\draw (19) to (21.center);
		\draw [in=165, out=-90] (23.center) to (19);
		\draw [in=-90, out=15] (19) to (24.center);
		\draw (20) to (23.center);
		\draw [in=165, out=-90] (26.center) to (20);
		\draw [in=270, out=15] (20) to (27.center);
		\draw (28.center) to (24.center);
		\draw (29.center) to (27.center);
		\draw (30.center) to (28.center);
		\draw (32.center) to (31);
		\draw (34.center) to (26.center);
		\draw (41) to (43.center);
		\draw [in=165, out=-90] (45.center) to (41);
		\draw [in=-90, out=15] (41) to (46.center);
		\draw (42) to (45.center);
		\draw [in=165, out=-90] (48.center) to (42);
		\draw [in=270, out=15] (42) to (49.center);
		\draw (50.center) to (46.center);
		\draw (51.center) to (49.center);
		\draw (52.center) to (50.center);
		\draw (54.center) to (53);
		\draw (56.center) to (48.center);
		\draw (87) to (89.center);
		\draw [in=165, out=-90] (91.center) to (87);
		\draw [in=-90, out=15] (87) to (92.center);
		\draw (88) to (91.center);
		\draw [in=165, out=-90] (94.center) to (88);
		\draw [in=270, out=15] (88) to (95.center);
		\draw (96.center) to (92.center);
		\draw (97.center) to (95.center);
		\draw (98.center) to (96.center);
		\draw (100.center) to (99);
		\draw (102.center) to (94.center);
		\draw (108) to (110.center);
		\draw [in=165, out=-90] (112.center) to (108);
		\draw [in=-90, out=15] (108) to (113.center);
		\draw (109) to (112.center);
		\draw [in=165, out=-90] (115.center) to (109);
		\draw [in=270, out=15] (109) to (116.center);
		\draw (117.center) to (113.center);
		\draw (118.center) to (116.center);
		\draw (123.center) to (115.center);
		\draw (135.center) to (104);
		\draw (136.center) to (138.center);
		\draw (137.center) to (139.center);
		\draw (121.center) to (125);
		\draw (141) to (143.center);
		\draw [in=165, out=-90] (145.center) to (141);
		\draw [in=-90, out=15] (141) to (146.center);
		\draw (142) to (145.center);
		\draw [in=165, out=-90] (148.center) to (142);
		\draw [in=270, out=15] (142) to (149.center);
		\draw (150.center) to (146.center);
		\draw (151.center) to (149.center);
		\draw (154.center) to (148.center);
		\draw (152.center) to (156);
		\draw [in=90, out=-90] (157.center) to (117.center);
		\draw [in=90, out=-90] (158.center) to (157.center);
	\end{pgfonlayer}
\end{tikzpicture}
		      \]
		      using determinism of $\delta$, naturality of $\samp$, monoidality of $\samp$ (see \cite[Remark 3.16]{Fritzetal-2023-Representablemarkovc}), and $\samp\circ\delta=1$. The same can be done with $g$, and so, by the left equality in \eqref{eq:alt-asc} we get the right equality.

		\item Suppose now that $\cat{C}$ has conditionals.
		      Suppose first of all that for every $p:A\to X$, $f,g:X\otimes A\to Y$ and $\pi:A\to PX$ we have the reverse implication to \eqref{eq:alt-asc}. Take again $\pi=\delta\circ p$, and applying $\samp$ to the first output of the first term in \eqref{eq:alt-asc}, we get
		      
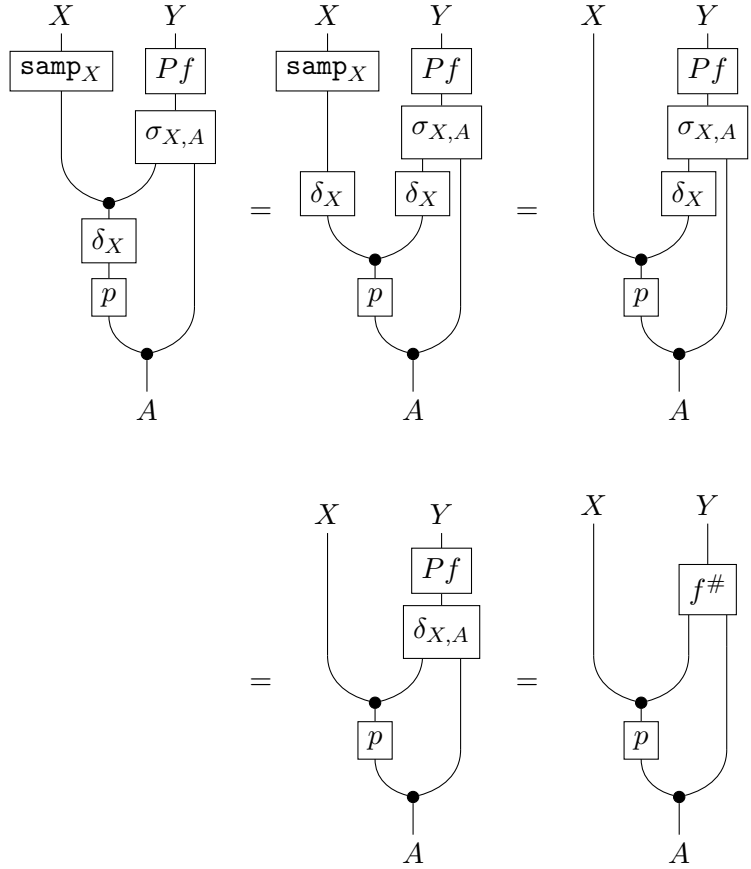
\[
			      \begin{tikzpicture}
	\begin{pgfonlayer}{nodelayer}
		\node [style=bn] (19) at (-7, -3.75) {};
		\node [style=bn] (20) at (-8, 0.25) {};
		\node [style=none] (21) at (-7, -4.75) {};
		\node [style=none] (22) at (-7, -5.25) {$A$};
		\node [style=none] (23) at (-8, -2.75) {};
		\node [style=none] (24) at (-5.75, -2.5) {};
		\node [style=morphism] (25) at (-8, -2.25) {$p$};
		\node [style=none] (26) at (-9.25, 1.5) {};
		\node [style=none] (27) at (-6.75, 1.5) {};
		\node [style=none] (28) at (-5.75, 1.5) {};
		\node [style=none] (29) at (-6.75, 2) {};
		\node [style=none] (30) at (-5.75, 2) {};
		\node [style=medium box] (31) at (-6.25, 2) {$\sigma_{X,A}$};
		\node [style=none] (32) at (-6.25, 4.75) {};
		\node [style=none] (33) at (-9.25, 5.25) {$X$};
		\node [style=none] (34) at (-9.25, 4.75) {};
		\node [style=none] (35) at (-6.25, 5.25) {$PY$};
		\node [style=morphism] (36) at (-6.25, 3.75) {$Pf$};
		\node [style=none] (37) at (-4, 0) {$=$};
		\node [style=morphism] (38) at (-8, -0.75) {$\delta_X$};
		\node [style=morphism] (40) at (-9.25, 3.75) {$\samp_X$};
		\node [style=bn] (41) at (0, -3.75) {};
		\node [style=bn] (42) at (-1, -1.25) {};
		\node [style=none] (43) at (0, -4.75) {};
		\node [style=none] (44) at (0, -5.25) {$A$};
		\node [style=none] (45) at (-1, -2.75) {};
		\node [style=none] (46) at (1.25, -2.5) {};
		\node [style=morphism] (47) at (-1, -2.25) {$p$};
		\node [style=none] (48) at (-2.25, 0) {};
		\node [style=none] (49) at (0.25, 0) {};
		\node [style=none] (50) at (1.25, 1.5) {};
		\node [style=none] (51) at (0.25, 2) {};
		\node [style=none] (52) at (1.25, 2) {};
		\node [style=medium box] (53) at (0.75, 2.125) {$\sigma_{X,A}$};
		\node [style=none] (54) at (0.75, 4.75) {};
		\node [style=none] (55) at (-2.25, 5.25) {$X$};
		\node [style=none] (56) at (-2.25, 4.75) {};
		\node [style=none] (57) at (0.75, 5.25) {$PY$};
		\node [style=morphism] (58) at (0.75, 3.75) {$Pf$};
		\node [style=morphism] (59) at (-2.25, 0.5) {$\delta_X$};
		\node [style=morphism] (61) at (-2.25, 3.75) {$\samp_X$};
		\node [style=morphism] (62) at (0.25, 0.5) {$\delta_X$};
		\node [style=none] (86) at (3, 0) {$=$};
		\node [style=bn] (159) at (7, -3.75) {};
		\node [style=bn] (160) at (6, -1.25) {};
		\node [style=none] (161) at (7, -4.75) {};
		\node [style=none] (162) at (7, -5.25) {$A$};
		\node [style=none] (163) at (6, -2.75) {};
		\node [style=none] (164) at (8.25, -2.5) {};
		\node [style=morphism] (165) at (6, -2.25) {$p$};
		\node [style=none] (166) at (4.75, 0) {};
		\node [style=none] (167) at (7.25, 0) {};
		\node [style=none] (168) at (8.25, 1.5) {};
		\node [style=none] (169) at (7.25, 2) {};
		\node [style=none] (170) at (8.25, 2) {};
		\node [style=medium box] (171) at (7.75, 2.125) {$\sigma_{X,A}$};
		\node [style=none] (172) at (7.75, 4.75) {};
		\node [style=none] (173) at (4.75, 5.25) {$X$};
		\node [style=none] (174) at (4.75, 4.75) {};
		\node [style=none] (175) at (7.75, 5.25) {$PY$};
		\node [style=morphism] (176) at (7.75, 3.75) {$Pf$};
		\node [style=morphism] (179) at (7.25, 0.5) {$\delta_X$};
		\node [style=none] (180) at (-4, -12.5) {$=$};
		\node [style=bn] (181) at (0, -15.5) {};
		\node [style=bn] (182) at (-1, -13) {};
		\node [style=none] (183) at (0, -16.5) {};
		\node [style=none] (184) at (0, -17) {$A$};
		\node [style=none] (185) at (-1, -14.5) {};
		\node [style=none] (186) at (1.25, -14.25) {};
		\node [style=morphism] (187) at (-1, -14) {$p$};
		\node [style=none] (188) at (-2.25, -11.75) {};
		\node [style=none] (189) at (0.25, -11.75) {};
		\node [style=none] (190) at (1.25, -11.75) {};
		\node [style=none] (191) at (0.25, -11.25) {};
		\node [style=none] (192) at (1.25, -11.25) {};
		\node [style=medium box] (193) at (0.75, -11.125) {$\delta_{X,A}$};
		\node [style=none] (194) at (0.75, -8.5) {};
		\node [style=none] (195) at (-2.25, -8) {$X$};
		\node [style=none] (196) at (-2.25, -8.5) {};
		\node [style=none] (197) at (0.75, -8) {$PY$};
		\node [style=morphism] (198) at (0.75, -9.5) {$Pf$};
		\node [style=none] (199) at (3, -12.5) {$=$};
		\node [style=bn] (200) at (7, -15.5) {};
		\node [style=bn] (201) at (6, -13) {};
		\node [style=none] (202) at (7, -16.5) {};
		\node [style=none] (203) at (7, -17) {$A$};
		\node [style=none] (204) at (6, -14.5) {};
		\node [style=none] (205) at (8.25, -14.25) {};
		\node [style=morphism] (206) at (6, -14) {$p$};
		\node [style=none] (207) at (4.75, -11.75) {};
		\node [style=none] (208) at (7.25, -11.75) {};
		\node [style=none] (209) at (8.25, -11.5) {};
		\node [style=none] (210) at (7.25, -10) {};
		\node [style=none] (211) at (8.25, -10) {};
		\node [style=none] (213) at (7.75, -8.25) {};
		\node [style=none] (214) at (4.75, -7.75) {$X$};
		\node [style=none] (215) at (4.75, -8.25) {};
		\node [style=none] (216) at (7.75, -7.75) {$PY$};
		\node [style=morphism] (217) at (7.75, -10) {$f^\#$};
	\end{pgfonlayer}
	\begin{pgfonlayer}{edgelayer}
		\draw (19) to (21.center);
		\draw [in=165, out=-90] (23.center) to (19);
		\draw [in=-90, out=15] (19) to (24.center);
		\draw (20) to (23.center);
		\draw [in=165, out=-90] (26.center) to (20);
		\draw [in=270, out=15] (20) to (27.center);
		\draw (28.center) to (24.center);
		\draw (29.center) to (27.center);
		\draw (30.center) to (28.center);
		\draw (32.center) to (31);
		\draw (34.center) to (26.center);
		\draw (41) to (43.center);
		\draw [in=165, out=-90] (45.center) to (41);
		\draw [in=-90, out=15] (41) to (46.center);
		\draw (42) to (45.center);
		\draw [in=165, out=-90] (48.center) to (42);
		\draw [in=270, out=15] (42) to (49.center);
		\draw (50.center) to (46.center);
		\draw (51.center) to (49.center);
		\draw (52.center) to (50.center);
		\draw (54.center) to (53);
		\draw (56.center) to (48.center);
		\draw (159) to (161.center);
		\draw [in=165, out=-90] (163.center) to (159);
		\draw [in=-90, out=15] (159) to (164.center);
		\draw (160) to (163.center);
		\draw [in=165, out=-90] (166.center) to (160);
		\draw [in=270, out=15] (160) to (167.center);
		\draw (168.center) to (164.center);
		\draw (169.center) to (167.center);
		\draw (170.center) to (168.center);
		\draw (172.center) to (171);
		\draw (174.center) to (166.center);
		\draw (181) to (183.center);
		\draw [in=165, out=-90] (185.center) to (181);
		\draw [in=-90, out=15] (181) to (186.center);
		\draw (182) to (185.center);
		\draw [in=165, out=-90] (188.center) to (182);
		\draw [in=270, out=15] (182) to (189.center);
		\draw (190.center) to (186.center);
		\draw (191.center) to (189.center);
		\draw (192.center) to (190.center);
		\draw (194.center) to (193);
		\draw (196.center) to (188.center);
		\draw (200) to (202.center);
		\draw [in=165, out=-90] (204.center) to (200);
		\draw [in=-90, out=15] (200) to (205.center);
		\draw (201) to (204.center);
		\draw [in=165, out=-90] (207.center) to (201);
		\draw [in=270, out=15] (201) to (208.center);
		\draw (209.center) to (205.center);
		\draw (210.center) to (208.center);
		\draw (211.center) to (209.center);
		\draw (215.center) to (207.center);
		\draw (213.center) to (217);
	\end{pgfonlayer}
\end{tikzpicture}
		      \]
		      using determinism of $\delta$, $\samp\circ\delta=1$, monoidality of $\delta$, and the definition of $f^\sharp$. The same can be done with $g$, and so, by the right equality in \eqref{eq:alt-asc}, we get exactly the sampling cancellation property \eqref{eq:samp-canc}.

		      Conversely, suppose that $\cat{C}$ is a.s.-compatibly representable (and has conditionals).
		      Let $p:A\to X$, $f,g:X\otimes A\to Y$, $\pi:A\to PX$ such that $\samp\circ\pi=p$, and suppose that the right equality in \eqref{eq:alt-asc} holds.
		      To prove the left equality, by the sampling cancellation property it suffices to show that
		      \begin{equation}\label{eq:alt-asc-proof3}
			      \begin{tikzpicture}
	\begin{pgfonlayer}{nodelayer}
		\node [style=bn] (0) at (4.25, -3.5) {};
		\node [style=bn] (1) at (3.25, -1) {};
		\node [style=none] (2) at (4.25, -4.5) {};
		\node [style=none] (3) at (4.25, -5) {$A$};
		\node [style=none] (4) at (3.25, -2.5) {};
		\node [style=none] (5) at (5.25, -2.5) {};
		\node [style=morphism] (6) at (3.25, -2) {$\pi$};
		\node [style=none] (7) at (2.25, 0) {};
		\node [style=none] (8) at (4.25, 0) {};
		\node [style=none] (9) at (5.25, 0) {};
		\node [style=none] (10) at (4.25, 0.5) {};
		\node [style=none] (11) at (5.25, 0.5) {};
		\node [style=medium box] (12) at (4.75, 0.5) {$\sigma_{X,A}$};
		\node [style=none] (13) at (4.75, 4.75) {};
		\node [style=none] (14) at (2.25, 5.25) {$PX$};
		\node [style=none] (15) at (2.25, 4.75) {};
		\node [style=none] (16) at (4.75, 5.25) {$Y$};
		\node [style=morphism] (17) at (4.75, 2.25) {$Pg$};
		\node [style=bn] (19) at (-3.75, -3.5) {};
		\node [style=bn] (20) at (-4.75, -1) {};
		\node [style=none] (21) at (-3.75, -4.5) {};
		\node [style=none] (22) at (-3.75, -5) {$A$};
		\node [style=none] (23) at (-4.75, -2.5) {};
		\node [style=none] (24) at (-2.75, -2.5) {};
		\node [style=morphism] (25) at (-4.75, -2) {$\pi$};
		\node [style=none] (26) at (-5.75, 0) {};
		\node [style=none] (27) at (-3.75, 0) {};
		\node [style=none] (28) at (-2.75, 0) {};
		\node [style=none] (29) at (-3.75, 0.5) {};
		\node [style=none] (30) at (-2.75, 0.5) {};
		\node [style=medium box] (31) at (-3.25, 0.5) {$\sigma_{X,A}$};
		\node [style=none] (32) at (-3.25, 4.75) {};
		\node [style=none] (33) at (-5.75, 5.25) {$PX$};
		\node [style=none] (34) at (-5.75, 4.75) {};
		\node [style=none] (35) at (-3.25, 5.25) {$Y$};
		\node [style=morphism] (36) at (-3.25, 2.25) {$Pf$};
		\node [style=none] (37) at (0, 0) {$=$};
		\node [style=morphism] (38) at (-3.25, 3.75) {$\samp_Y$};
		\node [style=morphism] (39) at (4.75, 3.75) {$\samp_Y$};
	\end{pgfonlayer}
	\begin{pgfonlayer}{edgelayer}
		\draw (0) to (2.center);
		\draw [in=165, out=-90] (4.center) to (0);
		\draw [in=-90, out=15] (0) to (5.center);
		\draw (1) to (4.center);
		\draw [in=165, out=-90] (7.center) to (1);
		\draw [in=270, out=15] (1) to (8.center);
		\draw (9.center) to (5.center);
		\draw (10.center) to (8.center);
		\draw (11.center) to (9.center);
		\draw (13.center) to (12);
		\draw (15.center) to (7.center);
		\draw (19) to (21.center);
		\draw [in=165, out=-90] (23.center) to (19);
		\draw [in=-90, out=15] (19) to (24.center);
		\draw (20) to (23.center);
		\draw [in=165, out=-90] (26.center) to (20);
		\draw [in=270, out=15] (20) to (27.center);
		\draw (28.center) to (24.center);
		\draw (29.center) to (27.center);
		\draw (30.center) to (28.center);
		\draw (32.center) to (31);
		\draw (34.center) to (26.center);
	\end{pgfonlayer}
\end{tikzpicture}
		      \end{equation}
		      Consider now the Bayesian inverse $(\samp_X)^\dag_\pi$ appearing in the following equation.
		      \begin{equation}\label{eq:bayesian-samp}
			      \begin{tikzpicture}
	\begin{pgfonlayer}{nodelayer}
		\node [style=bn] (0) at (3, 0) {};
		\node [style=none] (3) at (4.5, 2.75) {};
		\node [style=none] (4) at (4.5, 3.25) {$PX$};
		\node [style=none] (5) at (2, 1) {};
		\node [style=none] (6) at (2, 3.25) {$X$};
		\node [style=bn] (8) at (-3, 0) {};
		\node [style=morphism] (9) at (-4, 1.5) {$\samp_X$};
		\node [style=none] (10) at (-4, 2.75) {};
		\node [style=none] (11) at (-4, 3.25) {$X$};
		\node [style=none] (12) at (-2, 1) {};
		\node [style=none] (13) at (-2, 3.25) {$PX$};
		\node [style=none] (15) at (0, -0.25) {$=$};
		\node [style=morphism] (16) at (3, -1) {$p$};
		\node [style=none] (17) at (4, -3) {};
		\node [style=none] (18) at (4, -3.5) {$A$};
		\node [style=bn] (19) at (4, -2.25) {};
		\node [style=none] (20) at (5, -1.25) {};
		\node [style=none] (21) at (4, 1) {};
		\node [style=medium box] (22) at (4.5, 1.5) {$(\samp_X)^\dagger_\pi$};
		\node [style=none] (23) at (-3, -3) {};
		\node [style=none] (24) at (-3, -3.5) {$A$};
		\node [style=morphism] (25) at (-3, -1.5) {$\pi$};
		\node [style=none] (26) at (-2, 2.75) {};
		\node [style=none] (27) at (2, 2.75) {};
		\node [style=none] (28) at (5, 1) {};
		\node [style=none] (29) at (-4, 1) {};
		\node [style=none] (30) at (3, -1.25) {};
	\end{pgfonlayer}
	\begin{pgfonlayer}{edgelayer}
		\draw [in=-90, out=165] (0) to (5.center);
		\draw [in=-90, out=15] (8) to (12.center);
		\draw (9) to (10.center);
		\draw (16) to (0);
		\draw (19) to (17.center);
		\draw [in=15, out=-90] (20.center) to (19);
		\draw [in=15, out=-90] (21.center) to (0);
		\draw (3.center) to (22);
		\draw (25) to (23.center);
		\draw (8) to (25);
		\draw (12.center) to (26.center);
		\draw (27.center) to (5.center);
		\draw (28.center) to (20.center);
		\draw [in=165, out=-90] (29.center) to (8);
		\draw (29.center) to (9);
		\draw [in=165, out=-90] (30.center) to (19);
		\draw (16) to (30.center);
	\end{pgfonlayer}
\end{tikzpicture}
		      \end{equation}
		      The left-hand side of \eqref{eq:alt-asc-proof3} can be transformed as follows,
		      \[
			      \input{alt-asc-proof4.tikz}
		      \]
		      using naturality of $\samp$, monoidality of $\samp$, \eqref{eq:bayesian-samp}, the fact that $\samp\circ\pi=p$, and the comonoid axioms.
		      The same can be done for $g$. Therefore in the last term above we can substitute $f$ by $g$ by the right equality in \eqref{eq:alt-asc}, and we obtain \eqref{eq:alt-asc-proof3}, which implies the left equality in \eqref{eq:alt-asc} by the sampling cancellation property. \qedhere
	\end{enumerate}
\end{proof}

Let us now prove the main statement.

\begin{proof}[Proof of \Cref{thm:pullbacks}]
	Let \(f:X\to Y\) be deterministic.
	Using \Cref{cor:BC_Kleisli}, let \( p:A \to X \), and \( q: A \to PY \) be
	morphisms in $\cat{C}$ such that the following diagram is
	commutative.
	\begin{equation}\label{eq:outer}
		\begin{tikzcd}
			A \\
			& PX & PY \\
			& X & Y
			\arrow["q", curve={height=-18pt}, from=1-1, to=2-3]
			\arrow["p"', curve={height=18pt}, from=1-1, to=3-2]
			\arrow["Pf", from=2-2, to=2-3]
			\arrow["\samp"', from=2-2, to=3-2]
			\arrow["\samp", from=2-3, to=3-3]
			\arrow["f"', from=3-2, to=3-3]
		\end{tikzcd}
	\end{equation}
	We need to find a map \( r: A \to PX \) in $\cat{C}$ such that the following diagram commutes.
	\[
		\begin{tikzcd}
			A \\
			& PX & PY \\
			& X & Y
			\arrow["r"{description}, dashed, from=1-1, to=2-2]
			\arrow["q", curve={height=-18pt}, from=1-1, to=2-3]
			\arrow["p"', curve={height=18pt}, from=1-1, to=3-2]
			\arrow["Pf", from=2-2, to=2-3]
			\arrow["\samp"', from=2-2, to=3-2]
			\arrow["\samp", from=2-3, to=3-3]
			\arrow["f"', from=3-2, to=3-3]
		\end{tikzcd}
	\]
	Consider therefore the following morphism in $\cat{C}$,
	\begin{equation}
		\begin{tikzpicture}
	\begin{pgfonlayer}{nodelayer}
		\node [style=none] (0) at (-3, -4.25) {};
		\node [style=none] (1) at (-3, -4.75) {$A$};
		\node [style=none] (2) at (-3, 5.25) {$PX$};
		\node [style=none] (3) at (-3, 4.75) {};
		\node [style=medium box] (4) at (-3, 0) {$r$};
		\node [style=none] (5) at (0, 0) {$:=$};
		\node [style=none] (6) at (4.5, -4.25) {};
		\node [style=none] (7) at (4.5, -4.75) {$A$};
		\node [style=none] (8) at (4.5, 4.75) {};
		\node [style=bn] (10) at (4.5, -3) {};
		\node [style=morphism] (12) at (4.5, 3.5) {$P(f^\dagger_p)$};
		\node [style=none] (15) at (5.5, 0.75) {};
		\node [style=none] (16) at (3.5, -2) {};
		\node [style=none] (19) at (4.5, 5.25) {$PX$};
		\node [style=medium box] (20) at (4.5, 0.75) {$\;\sigma_{Y,A}\;$};
		\node [style=none] (21) at (5.5, -2) {};
		\node [style=none] (22) at (3.5, 0.75) {};
		\node [style=none] (23) at (6.5, 2.125) {$P(Y\otimes A)$};
		\node [style=none] (24) at (2.5, -0.5) {$PY$};
		\node [style=morphism] (25) at (3.5, -1.5) {$q$};
	\end{pgfonlayer}
	\begin{pgfonlayer}{edgelayer}
		\draw (4) to (0.center);
		\draw (10) to (6.center);
		\draw (3.center) to (4);
		\draw (8.center) to (12);
		\draw [in=165, out=-90] (16.center) to (10);
		\draw (12) to (20);
		\draw [in=15, out=-90] (21.center) to (10);
		\draw (22.center) to (16.center);
		\draw (15.center) to (21.center);
	\end{pgfonlayer}
\end{tikzpicture}
		\label{eq:r-def}
	\end{equation}
	where $f^\dag:Y\otimes A\to X$ is the Bayesian inverse of $f$ with respect to $p$,
	\begin{equation}
		\label{eq:bayesian2}
		\begin{tikzpicture}
	\begin{pgfonlayer}{nodelayer}
		\node [style=bn] (0) at (3, 0.25) {};
		\node [style=none] (3) at (4.5, 3) {};
		\node [style=none] (4) at (4.5, 3.5) {$X$};
		\node [style=none] (5) at (2, 1.25) {};
		\node [style=none] (6) at (2, 3.5) {$Y$};
		\node [style=bn] (8) at (-3, -0.5) {};
		\node [style=morphism] (9) at (-4, 1.5) {$f$};
		\node [style=none] (10) at (-4, 3) {};
		\node [style=none] (11) at (-4, 3.5) {$Y$};
		\node [style=none] (12) at (-2, 0.5) {};
		\node [style=none] (13) at (-2, 3.5) {$X$};
		\node [style=none] (15) at (0, -0.25) {$=$};
		\node [style=morphism] (16) at (3, -2.25) {$p$};
		\node [style=none] (17) at (4, -4.5) {};
		\node [style=none] (18) at (4, -5) {$A$};
		\node [style=bn] (19) at (4, -3.75) {};
		\node [style=none] (20) at (5, -2.75) {};
		\node [style=none] (21) at (4, 1.25) {};
		\node [style=medium box] (22) at (4.5, 1.75) {$f^\dagger_p$};
		\node [style=none] (23) at (-3, -4.5) {};
		\node [style=none] (24) at (-3, -5) {$A$};
		\node [style=morphism] (25) at (-3, -2.5) {$p$};
		\node [style=none] (26) at (-2, 3) {};
		\node [style=none] (27) at (2, 3) {};
		\node [style=none] (28) at (5, 1.25) {};
		\node [style=none] (29) at (-4, 0.5) {};
		\node [style=none] (30) at (3, -2.75) {};
		\node [style=morphism] (32) at (3, -0.75) {$f$};
	\end{pgfonlayer}
	\begin{pgfonlayer}{edgelayer}
		\draw [in=-90, out=165] (0) to (5.center);
		\draw [in=-90, out=15] (8) to (12.center);
		\draw (9) to (10.center);
		\draw (16) to (0);
		\draw (19) to (17.center);
		\draw [in=15, out=-90] (20.center) to (19);
		\draw [in=15, out=-90] (21.center) to (0);
		\draw (3.center) to (22);
		\draw (25) to (23.center);
		\draw (8) to (25);
		\draw (12.center) to (26.center);
		\draw (27.center) to (5.center);
		\draw (28.center) to (20.center);
		\draw [in=165, out=-90] (29.center) to (8);
		\draw (29.center) to (9);
		\draw [in=165, out=-90] (30.center) to (19);
		\draw (16) to (30.center);
	\end{pgfonlayer}
\end{tikzpicture}
	\end{equation}
	and
	$\sigma_{Y,A}:PY \otimes A \to P(A \otimes Y)$ is the right strength of the monad $P$.
	
	To show that the left triangle commutes, i.e.\ that $\samp_X\circ\mu\circ r^\sharp = p $ we have that
	
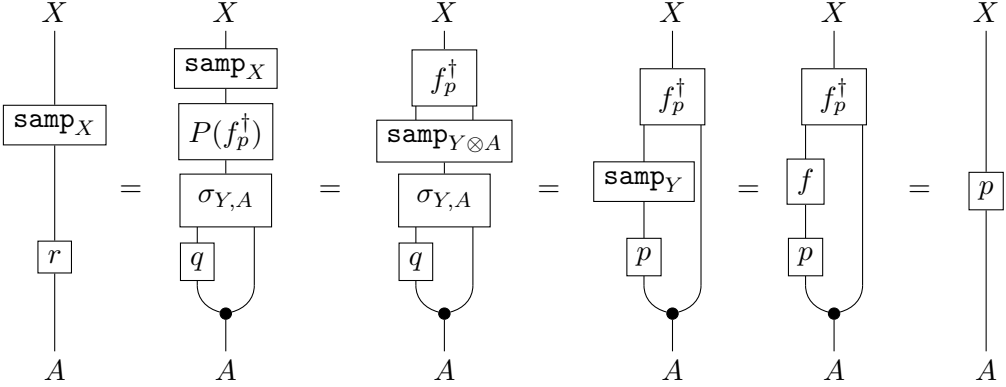
\begin{equation*}
		\begin{tikzpicture}
	\begin{pgfonlayer}{nodelayer}
		\node [style=none] (84) at (-14, -4.25) {};
		\node [style=none] (85) at (-14, -4.75) {$A$};
		\node [style=none] (86) at (-14, 4.75) {$X$};
		\node [style=none] (87) at (-14, 4.25) {};
		\node [style=morphism] (88) at (-14, -1.75) {$r$};
		\node [style=morphism] (89) at (-14, 1.75) {$\texttt{samp}_X$};
		\node [style=none] (90) at (-9.5, -4.25) {};
		\node [style=none] (91) at (-9.5, -4.75) {$A$};
		\node [style=none] (92) at (-9.5, 4.25) {};
		\node [style=bn] (93) at (-9.5, -3.25) {};
		\node [style=morphism] (94) at (-9.5, 1.575) {$P(f^\dagger_p)$};
		\node [style=morphism] (95) at (-10.25, -1.875) {$q$};
		\node [style=none] (96) at (-8.75, -0.25) {};
		\node [style=none] (97) at (-10.25, -2.5) {};
		\node [style=medium box] (98) at (-9.5, -0.25) {$\;\sigma_{Y,A}\;$};
		\node [style=none] (99) at (-8.75, -2.5) {};
		\node [style=none] (100) at (-10.25, -0.25) {};
		\node [style=none] (103) at (-12, 0) {$=$};
		\node [style=morphism] (104) at (-9.5, 3.25) {$\texttt{samp}_X$};
		\node [style=none] (105) at (-9.5, 4.75) {$X$};
		\node [style=none] (106) at (-3.75, -4.25) {};
		\node [style=none] (107) at (-3.75, -4.75) {$A$};
		\node [style=none] (108) at (-3.75, 4.25) {};
		\node [style=bn] (109) at (-3.75, -3.25) {};
		\node [style=morphism] (110) at (-3.75, 3) {$\;f^\dagger_p\;$};
		\node [style=morphism] (111) at (-4.5, -1.875) {$q$};
		\node [style=none] (112) at (-3, -0.25) {};
		\node [style=none] (113) at (-4.5, -2.5) {};
		\node [style=medium box] (114) at (-3.75, -0.25) {$\;\sigma_{Y,A}\;$};
		\node [style=none] (115) at (-3, -2.5) {};
		\node [style=none] (116) at (-4.5, -0.25) {};
		\node [style=none] (119) at (-6.75, 0) {$=$};
		\node [style=morphism] (120) at (-3.75, 1.325) {$\texttt{samp}_{Y\otimes A}$};
		\node [style=none] (121) at (-3.75, 4.75) {$X$};
		\node [style=none] (122) at (-4.5, 3) {};
		\node [style=none] (123) at (-4.5, 1.325) {};
		\node [style=none] (124) at (-3, 3) {};
		\node [style=none] (125) at (-3, 1.325) {};
		\node [style=none] (150) at (6.5, -4.25) {};
		\node [style=none] (151) at (6.5, -4.75) {$A$};
		\node [style=none] (152) at (6.5, 4.25) {};
		\node [style=bn] (153) at (6.5, -3.25) {};
		\node [style=morphism] (154) at (6.5, 2.5) {$\;f^\dagger_p\;$};
		\node [style=morphism] (155) at (5.75, -1.75) {$p$};
		\node [style=none] (156) at (7.25, 2.5) {};
		\node [style=none] (157) at (5.75, -2.5) {};
		\node [style=none] (158) at (7.25, -2.5) {};
		\node [style=none] (159) at (5.75, 2.5) {};
		\node [style=none] (161) at (4.25, 0) {$=$};
		\node [style=morphism] (162) at (5.75, 0.25) {$f$};
		\node [style=none] (163) at (6.5, 4.75) {$X$};
		\node [style=none] (166) at (10.5, -4.25) {};
		\node [style=none] (167) at (10.5, -4.75) {$A$};
		\node [style=none] (168) at (10.5, 4.75) {$X$};
		\node [style=none] (169) at (10.5, 4.25) {};
		\node [style=morphism] (170) at (10.5, 0) {$p$};
		\node [style=none] (172) at (8.75, 0) {$=$};
		\node [style=none] (173) at (2.25, -4.25) {};
		\node [style=none] (174) at (2.25, -4.75) {$A$};
		\node [style=none] (175) at (2.25, 4.25) {};
		\node [style=bn] (176) at (2.25, -3.25) {};
		\node [style=morphism] (177) at (2.25, 2.5) {$\;f^\dagger_p\;$};
		\node [style=morphism] (178) at (1.5, -1.75) {$q$};
		\node [style=none] (179) at (3, 2.5) {};
		\node [style=none] (180) at (1.5, -2.5) {};
		\node [style=none] (181) at (3, -2.5) {};
		\node [style=none] (182) at (1.5, 2.5) {};
		\node [style=none] (183) at (-1, 0) {$=$};
		\node [style=morphism] (184) at (1.5, 0.25) {$\samp_Y$};
		\node [style=none] (185) at (2.25, 4.75) {$X$};
	\end{pgfonlayer}
	\begin{pgfonlayer}{edgelayer}
		\draw (88) to (84.center);
		\draw (89) to (87.center);
		\draw (89) to (88);
		\draw (93) to (90.center);
		\draw (92.center) to (94);
		\draw [in=165, out=-90] (97.center) to (93);
		\draw (94) to (98);
		\draw [in=15, out=-90] (99.center) to (93);
		\draw (100.center) to (97.center);
		\draw (109) to (106.center);
		\draw (108.center) to (110);
		\draw [in=165, out=-90] (113.center) to (109);
		\draw [in=15, out=-90] (115.center) to (109);
		\draw (116.center) to (113.center);
		\draw (120) to (114);
		\draw (122.center) to (123.center);
		\draw (124.center) to (125.center);
		\draw (153) to (150.center);
		\draw (152.center) to (154);
		\draw [in=165, out=-90] (157.center) to (153);
		\draw [in=15, out=-90] (158.center) to (153);
		\draw (159.center) to (157.center);
		\draw (170) to (166.center);
		\draw (169.center) to (170);
		\draw (96.center) to (99.center);
		\draw (112.center) to (115.center);
		\draw (156.center) to (158.center);
		\draw (176) to (173.center);
		\draw (175.center) to (177);
		\draw [in=165, out=-90] (180.center) to (176);
		\draw [in=15, out=-90] (181.center) to (176);
		\draw (182.center) to (180.center);
		\draw (179.center) to (181.center);
	\end{pgfonlayer}
\end{tikzpicture}
	\end{equation*}
	using the definition of $r$, naturality of $\samp$, monoidality of $\samp$, commutativity of the outer diagram~\eqref{eq:outer}, and the second marginal of \eqref{eq:bayesian2}.
	
	Let us now show that the top triangle commutes, i.e.\ $Pf\circ r = q$.
	By \Cref{lem:f-as-det}, $f\circ f^\dag_p =_{f\circ p} 1_Y\otimes\texttt{del}_A$ in the sense of \Cref{lem:alt-asc}(2), with $X=Y$ and $p$ replaced by $f\circ p$.
	Since $\samp\circ q = f\circ p$ (by commutativity of the outer diagram~\eqref{eq:outer}), the converse direction of \Cref{lem:alt-asc}(2), instantiated with $\pi=q$, gives
	\[
		\begin{tikzpicture}
	\begin{pgfonlayer}{nodelayer}
		\node [style=none] (0) at (0, 0) {$=$};
		\node [style=none] (1) at (-3, -4) {};
		\node [style=none] (2) at (-3, -4.5) {$A$};
		\node [style=bn] (3) at (-3, -3) {};
		\node [style=morphism] (4) at (-4, -1.625) {$q$};
		\node [style=none] (5) at (-2, 1) {};
		\node [style=none] (6) at (-4, -2) {};
		\node [style=medium box] (7) at (-2.5, 1) {$\;\sigma_{Y,A}\;$};
		\node [style=none] (8) at (-2, -2) {};
		\node [style=bn] (10) at (-4, -0.5) {};
		\node [style=none] (12) at (-3, 0.5) {};
		\node [style=morphism] (13) at (-2.5, 2.875) {$P(f_p^\dagger)$};
		\node [style=morphism] (14) at (-2.5, 4.5) {$Pf$};
		\node [style=none] (15) at (-2.5, 5.5) {};
		\node [style=none] (18) at (-2.5, 6) {$PY$};
		\node [style=morphism] (32) at (4.25, 3.5) {$P(1_Y\!\otimes\!\texttt{del}_A)$};
		\node [style=none] (37) at (-3, 1) {};
		\node [style=none] (39) at (-5.25, 0.75) {};
		\node [style=none] (40) at (-5.25, 5.5) {};
		\node [style=none] (41) at (-5.25, 6) {$PY$};
		\node [style=none] (42) at (3.75, -4) {};
		\node [style=none] (43) at (3.75, -4.5) {$A$};
		\node [style=bn] (44) at (3.75, -3) {};
		\node [style=morphism] (45) at (2.75, -1.625) {$q$};
		\node [style=none] (46) at (4.75, 1) {};
		\node [style=none] (47) at (2.75, -2) {};
		\node [style=medium box] (48) at (4.25, 1) {$\;\sigma_{Y,A}\;$};
		\node [style=none] (49) at (4.75, -2) {};
		\node [style=bn] (50) at (2.75, -0.5) {};
		\node [style=none] (51) at (3.75, 0.5) {};
		\node [style=none] (54) at (4.25, 5.5) {};
		\node [style=none] (55) at (4.25, 6) {$PY$};
		\node [style=none] (56) at (3.75, 1) {};
		\node [style=none] (57) at (1.5, 0.75) {};
		\node [style=none] (58) at (1.5, 5.5) {};
		\node [style=none] (59) at (1.5, 6) {$PY$};
		\node [style=none] (60) at (7, 0) {$=$};
		\node [style=none] (62) at (10.75, -4) {};
		\node [style=none] (63) at (10.75, -4.5) {$A$};
		\node [style=bn] (64) at (10.75, -3) {};
		\node [style=morphism] (65) at (9.75, -1.625) {$q$};
		\node [style=none] (67) at (9.75, -2) {};
		\node [style=none] (69) at (11.75, -2) {};
		\node [style=bn] (70) at (9.75, -0.5) {};
		\node [style=none] (71) at (10.75, 0.5) {};
		\node [style=none] (73) at (10.75, 6) {$PY$};
		\node [style=none] (74) at (10.75, 5.5) {};
		\node [style=none] (75) at (8.75, 0.5) {};
		\node [style=none] (76) at (8.75, 5.5) {};
		\node [style=none] (77) at (8.75, 6) {$PY$};
		\node [style=bn] (78) at (11.75, 3) {};
	\end{pgfonlayer}
	\begin{pgfonlayer}{edgelayer}
		\draw (3) to (1.center);
		\draw [in=165, out=-90] (6.center) to (3);
		\draw [in=15, out=-90] (8.center) to (3);
		\draw [in=270, out=15] (10) to (12.center);
		\draw (5.center) to (8.center);
		\draw (10) to (6.center);
		\draw (37.center) to (12.center);
		\draw (7) to (15.center);
		\draw [in=165, out=-90] (39.center) to (10);
		\draw (40.center) to (39.center);
		\draw (44) to (42.center);
		\draw [in=165, out=-90] (47.center) to (44);
		\draw [in=15, out=-90] (49.center) to (44);
		\draw [in=270, out=15] (50) to (51.center);
		\draw (46.center) to (49.center);
		\draw (50) to (47.center);
		\draw (56.center) to (51.center);
		\draw (48) to (54.center);
		\draw [in=165, out=-90] (57.center) to (50);
		\draw (58.center) to (57.center);
		\draw (64) to (62.center);
		\draw [in=165, out=-90] (67.center) to (64);
		\draw [in=15, out=-90] (69.center) to (64);
		\draw [in=270, out=15] (70) to (71.center);
		\draw (70) to (67.center);
		\draw (74.center) to (71.center);
		\draw [in=165, out=-90] (75.center) to (70);
		\draw (76.center) to (75.center);
		\draw (78) to (69.center);
	\end{pgfonlayer}
\end{tikzpicture}
	\]
	using also naturality of the strength. Taking the second marginal, we get exactly $Pf\circ r = q$.
\end{proof}

In light of Proposition~\ref{prop:partialrel} we also have the following property on the
second-order stochastic dominance relation.
\begin{corollary}
	\label{cor:dominance}
	Let $\cat{C}$ be an a.s.-compatibly representable Markov category with
	conditionals. Then the second-order dominance
	relation on $\cat{C}(\Theta, A)$ is transitive for all $\Theta\in\cat{C}$ and all $P$-algebras $A$.
\end{corollary}

\subsection{The functor preserves weak pullbacks}\label{wcartp}

Let us now give a condition for when probability monads preserve weak pullbacks.
Before we begin, we need a condition relating equalizers in the category
$\cat{C}_{\det}$ (and hence finite limits, since $\cat{C}_{\det}$ has finite
products) and almost sure equality in the Markov category sense (which is
important for conditionals). 

Recall first of all that in traditional probability theory, two random
variables $f$ and $g$ on a probability space $(X,p)$ are almost surely equal if
and only if they agree on a set of probability one. This condition is
equivalent to saying that the equalizer of $f$ and $g$ has full measure, or again
equivalently, that the measure $p$, seen as a kernel $I\to X$, factors through
the equalizer of $f$ and $g$.

With this example in mind, let us define the following property of compatibility
between equalizers and almost sure equality in a Markov category.

\begin{definition}[{\cite[Definition~3.5.1]{Fritzetal-2023-Absolutecontinuitysu}}]\label{def:equalizer_principle}
	A Markov category $\cat{C}$ is said to satisfy the \textbf{equalizer principle} if:
	\begin{enumerate}[label=(\roman*)]
		\item Equalizers in $\cat{C}_{\text{det}}$ exist;
		\item For every equalizer diagram
		      \begin{equation*}
			      \begin{tikzcd}
				      E & X & Y
				      \arrow["{\cat{eq}}", from=1-1, to=1-2]
				      \arrow["f", shift left, from=1-2, to=1-3]
				      \arrow["g"', shift right, from=1-2, to=1-3]
			      \end{tikzcd}
		      \end{equation*}
		      in $\cat{C}_{\text{det}}$, every $p:A \to X$ in $\cat{C}$
		      satisfying $f =_p g$ factors uniquely across $\cat{eq}$.
	\end{enumerate}
\end{definition}

\begin{example}[{\cite[Proposition~3.5.4]{Fritzetal-2023-Absolutecontinuitysu}}]
	\label{lem:boreleq}
	$\cat{BorelStoch}$ satisfies the equalizer principle.
\end{example}

Here is now our main statement.

\begin{theorem}
	\label{thm:main}
	Let $\cat{C}$ be a representable Markov category with monad $(P,\mu,\delta)$.
	Suppose moreover that
	\begin{itemize}
		\item $\cat{C}$ has conditionals;
		\item $\cat{C}$ satisfies the equalizer principle.
	\end{itemize}
	Then the functor $P$ preserves weak pullbacks.
\end{theorem}

\begin{corollary}
	In the hypotheses of the theorem above, suppose $\cat{C}$ is not just representable, but also a.s.-compatibly so.
	Then, together with \Cref{thm:pullbacks}, the monad $(P,\mu,\delta)$ has weakly cartesian functor and multiplication.
\end{corollary}

\begin{corollary}\label{giryWC}
	The Giry monad on standard Borel spaces has weakly cartesian functor and multiplication.
\end{corollary}

We will use the following auxiliary statement, which is well known (see for example \cite{Constantinetal-2020-Partialevaluationsan}).

\begin{lemma}
	\label{lem:weakpullback}
	Let $\cat{C}$ be a category with pullbacks and $(P, \mu, \delta)$ a monad on $\cat{C}$.
	If $P$ turns pullbacks into weak pullbacks then $P$ preserves weak
	pullbacks.
\end{lemma}

\begin{remark}\label{rem:rel-pos}
	Recall that a Markov category $\cat{C}$ is called
	\emph{relatively positive}
	(called \emph{strictly positive}
	in~\cite[Definition~13.16]{Fritz-2020-Asyntheticapproachto};
	renamed in~\cite[Section~2.5]{Fritzetal-2023-Dilationsandinformat})
	if the positivity axiom holds with both antecedent and consequent
	relativized to $p$-almost sure equality.
	Existence of conditionals implies relative
	positivity~\cite[Lemma~13.17]{Fritz-2020-Asyntheticapproachto}.
	Equivalently, if $g \circ f$ is
	$p$-a.s.\ deterministic and $\cat{C}$ is relatively positive,
	then $g$ is $(f \circ p)$-a.s.\ deterministic.
\end{remark}

\begin{proof}[Proof of \Cref{thm:main}]
	First of all, since $\cat{C}_{\text{det}}$ has finite products (it is cartesian monoidal) and equalizers (by the equalizer principle), then it has all finite limits; in particular, it has all pullbacks.
	By \Cref{lem:weakpullback} it then suffices to show that $P$ turns
	pullbacks into weak pullbacks. Thus, suppose we have a pullback
	\[\begin{tikzcd}
			{X\times_Z Y} & Y \\
			X & Z.
			\arrow["{g^*f}", from=1-1, to=1-2]
			\arrow["{f^*g}"', from=1-1, to=2-1]
			\arrow["g", from=1-2, to=2-2]
			\arrow["f"', from=2-1, to=2-2]
		\end{tikzcd}\]
	in $\cat{C}_{\text{det}}$. Using \Cref{prop:BC_Kleisli}, we equivalently want to show that the diagram above is also a weak pullback in the Kleisli category. Suppose therefore that we have maps \( p:A \to X \)
	and \( q: A \to Y \), not necessarily deterministic, such that the following diagram commutes.
	\[\begin{tikzcd}
			A \\
			& {X\times_Z Y} & Y \\
			& X & Z.
			\arrow["{q}", curve={height=-18pt}, from=1-1, to=2-3]
			\arrow["{p}"', curve={height=18pt}, from=1-1, to=3-2]
			\arrow["{g^*f}", from=2-2, to=2-3]
			\arrow["{f^*g}"', from=2-2, to=3-2]
			\arrow["g", from=2-3, to=3-3]
			\arrow["f"', from=3-2, to=3-3]
		\end{tikzcd}\]
	Now form $\rho: A \to X \otimes Y$ as the conditional product
	\begin{equation*}
		\begin{tikzpicture}
	\begin{pgfonlayer}{nodelayer}
		\node [style=none] (0) at (-1, 0) {$:=$};
		\node [style=none] (1) at (-3, 0) {$\rho$};
		\node [style=none] (2) at (4, -3.25) {};
		\node [style=none] (3) at (4, -3.75) {$A$};
		\node [style=none] (4) at (2.5, 3) {};
		\node [style=none] (5) at (2.5, 3.5) {$X$};
		\node [style=none] (6) at (5.5, 3) {};
		\node [style=none] (7) at (5.5, 3.5) {$Y$};
		\node [style=bn] (8) at (3.25, -0.5) {};
		\node [style=bn] (9) at (4.75, -0.5) {};
		\node [style=bn] (10) at (4, -2.25) {};
		\node [style=medium box] (11) at (2.5, 1.5) {$\,\;f^\dagger_p\;\,$};
		\node [style=medium box] (12) at (5.5, 1.5) {$\,\;g^\dagger_q\;\,$};
		\node [style=none] (13) at (4.75, 0.75) {};
		\node [style=none] (14) at (6.25, 0.75) {};
		\node [style=none] (15) at (1.75, 0.75) {};
		\node [style=none] (16) at (3.25, 0.75) {};
		\node [style=none] (17) at (3.25, -1.5) {};
		\node [style=none] (18) at (4.75, -1.5) {};
		\node [style=morphism] (19) at (3.25, -1.25) {$s$};
		\node [style=none] (20) at (1.75, 1.25) {};
		\node [style=none] (21) at (3.25, 1.5) {};
		\node [style=none] (22) at (4.75, 1.5) {};
		\node [style=none] (23) at (6.25, 1.5) {};
	\end{pgfonlayer}
	\begin{pgfonlayer}{edgelayer}
		\draw [in=165, out=-90] (15.center) to (8);
		\draw [in=165, out=-90] (16.center) to (9);
		\draw [in=15, out=-90] (13.center) to (8);
		\draw (10) to (2.center);
		\draw (4.center) to (11);
		\draw (6.center) to (12);
		\draw [in=15, out=-90] (14.center) to (9);
		\draw (17.center) to (8);
		\draw (9) to (18.center);
		\draw [in=165, out=-90] (17.center) to (10);
		\draw [in=-90, out=15] (10) to (18.center);
		\draw (20.center) to (15.center);
		\draw (21.center) to (16.center);
		\draw (22.center) to (13.center);
		\draw (23.center) to (14.center);
	\end{pgfonlayer}
\end{tikzpicture}
	\end{equation*}
	where \( s \) is the common composition \( f \circ p = g \circ q \).
	As one can readily check, the marginals of $\rho$ are $p$ and $q$ respectively.
	In order to prove the theorem, it then suffices to show that it factors through $X\times_Z Y$ (which, in $\cat{C}_{\text{det}}$, is a subobject of $X\times Y=X\otimes Y$).
	
	Recall now that the pullback $X\times_Z Y$ in $\cat{C}_{\text{det}}$ can be expressed as the following equalizer, also in $\cat{C}_{\text{det}}$:
	\[
		\begin{tikzcd}
			X\times_Z Y \ar[hook]{r} & X\times Y \ar[shift left]{rr}{f\circ\pi_X} \ar[shift right]{rr}[swap]{g\circ\pi_Y} && Z
		\end{tikzcd}
	\]
	The equalizer principle tells us that a sufficient condition for $\rho$ to factor through $X\times_Z Y$ is that the maps $f\circ\pi_X$ and $g\circ\pi_Y$ are $\rho$-almost surely equal, i.e. that the following equation holds.
	\begin{equation}
		\label{eq:eq-prin-1}
		\begin{tikzpicture}
	\begin{pgfonlayer}{nodelayer}
		\node [style=none] (0) at (0, 0) {$=$};
		\node [style=none] (2) at (-5, -3.25) {};
		\node [style=none] (3) at (-5, -3.75) {$A$};
		\node [style=bn] (5) at (-5.75, 0) {};
		\node [style=bn] (6) at (-4.25, 0) {};
		\node [style=morphism] (7) at (-7.5, 1.75) {$f$};
		\node [style=bn] (8) at (-6, 2) {};
		\node [style=none] (9) at (-7.5, 3.25) {$Z$};
		\node [style=none] (10) at (-7.5, 2.75) {};
		\node [style=none] (11) at (-4, 3.25) {$X$};
		\node [style=none] (12) at (-4, 1.25) {};
		\node [style=none] (13) at (-2.5, 3.25) {$Y$};
		\node [style=none] (14) at (-2.5, 1.25) {};
		\node [style=medium box] (15) at (-5, -1.5) {$\quad\rho\quad$};
		\node [style=none] (16) at (-5.75, -1.5) {};
		\node [style=none] (17) at (-4.25, -1.5) {};
		\node [style=none] (35) at (-6, 1.25) {};
		\node [style=none] (36) at (-7.5, 1.25) {};
		\node [style=none] (37) at (-4, 2.75) {};
		\node [style=none] (38) at (-2.5, 2.75) {};
		\node [style=none] (39) at (5.5, -3.25) {};
		\node [style=none] (40) at (5.5, -3.75) {$A$};
		\node [style=bn] (41) at (4.75, 0) {};
		\node [style=bn] (42) at (6.25, 0) {};
		\node [style=morphism] (43) at (3, 1.75) {$g$};
		\node [style=bn] (44) at (4.5, 2) {};
		\node [style=none] (45) at (3, 3.25) {$Z$};
		\node [style=none] (46) at (3, 2.75) {};
		\node [style=none] (47) at (6.5, 3.25) {$X$};
		\node [style=none] (48) at (6.5, 1.25) {};
		\node [style=none] (49) at (8, 3.25) {$Y$};
		\node [style=none] (50) at (8, 1.25) {};
		\node [style=medium box] (51) at (5.5, -1.5) {$\quad\rho\quad$};
		\node [style=none] (52) at (4.75, -1.5) {};
		\node [style=none] (53) at (6.25, -1.5) {};
		\node [style=none] (54) at (4.5, 1.25) {};
		\node [style=none] (55) at (3, 1.25) {};
		\node [style=none] (56) at (6.5, 2.75) {};
		\node [style=none] (57) at (8, 2.75) {};
	\end{pgfonlayer}
	\begin{pgfonlayer}{edgelayer}
		\draw [in=270, out=15] (5) to (12.center);
		\draw [in=270, out=15] (6) to (14.center);
		\draw (17.center) to (6);
		\draw (16.center) to (5);
		\draw (15) to (2.center);
		\draw (8) to (35.center);
		\draw [in=270, out=165] (6) to (35.center);
		\draw [in=165, out=-90] (36.center) to (5);
		\draw (10.center) to (36.center);
		\draw (37.center) to (12.center);
		\draw (38.center) to (14.center);
		\draw [in=270, out=15] (41) to (48.center);
		\draw [in=270, out=15] (42) to (50.center);
		\draw (53.center) to (42);
		\draw (52.center) to (41);
		\draw (51) to (39.center);
		\draw (44) to (54.center);
		\draw [in=270, out=165] (42) to (54.center);
		\draw [in=165, out=-90] (55.center) to (41);
		\draw (46.center) to (55.center);
		\draw (56.center) to (48.center);
		\draw (57.center) to (50.center);
	\end{pgfonlayer}
\end{tikzpicture}
	\end{equation}
	The left-hand side of \eqref{eq:eq-prin-1} is now equal to the following,
	\begin{equation*}
		\resizebox{\linewidth}{!}{\begin{tikzpicture}
	\begin{pgfonlayer}{nodelayer}
		\node [style=none] (0) at (0, 0) {$=$};
		\node [style=bn] (2) at (-5.5, 1.75) {};
		\node [style=morphism] (4) at (-6.5, 3.25) {$f$};
		\node [style=none] (6) at (-6.5, 5) {$Z$};
		\node [style=none] (7) at (-6.5, 4.5) {};
		\node [style=none] (8) at (-4.5, 5) {$X$};
		\node [style=none] (9) at (-4.5, 2.75) {};
		\node [style=none] (10) at (-2.5, 5) {$Y$};
		\node [style=none] (11) at (-2.5, 4.5) {};
		\node [style=bn] (28) at (2.25, 0.5) {};
		\node [style=none] (32) at (1.5, 5) {$Z$};
		\node [style=none] (33) at (1.5, 1.25) {};
		\node [style=none] (34) at (3.75, 5) {$X$};
		\node [style=none] (35) at (3.75, 4.5) {};
		\node [style=none] (36) at (6, 5) {$Y$};
		\node [style=none] (37) at (6, 4.5) {};
		\node [style=none] (39) at (6.25, 0.25) {};
		\node [style=none] (40) at (4.5, -4.5) {};
		\node [style=none] (41) at (4.5, -5) {$A$};
		\node [style=bn] (44) at (3.75, -1.75) {};
		\node [style=bn] (45) at (5.25, -1.75) {};
		\node [style=bn] (46) at (4.5, -3.5) {};
		\node [style=medium box] (47) at (3.75, 2.5) {$\,\;f^\dagger_p\;\,$};
		\node [style=medium box] (48) at (6, 0.25) {$\,\;g^\dagger_q\;\,$};
		\node [style=none] (49) at (5.25, -0.25) {};
		\node [style=none] (50) at (6.75, -0.25) {};
		\node [style=none] (52) at (4.5, -1) {};
		\node [style=none] (54) at (3, 1.25) {};
		\node [style=none] (55) at (-7.75, 0) {$=$};
		\node [style=none] (75) at (10, 5) {$Z$};
		\node [style=none] (76) at (10, 1.75) {};
		\node [style=none] (77) at (12, 5) {$X$};
		\node [style=none] (78) at (12, 4.5) {};
		\node [style=none] (79) at (14.5, 5) {$Y$};
		\node [style=none] (80) at (14.5, 4.5) {};
		\node [style=none] (82) at (13.25, -4.5) {};
		\node [style=none] (83) at (13.25, -5) {$A$};
		\node [style=none] (84) at (14.5, 2.25) {};
		\node [style=bn] (85) at (12.5, 0) {};
		\node [style=bn] (86) at (14, 0) {};
		\node [style=bn] (87) at (13.25, -2.75) {};
		\node [style=medium box] (88) at (12, 2.25) {$\,\;f^\dagger_p\;\,$};
		\node [style=medium box] (89) at (14.5, 2.25) {$\,\;g^\dagger_q\;\,$};
		\node [style=none] (90) at (13.75, 1.5) {};
		\node [style=none] (91) at (15.25, 1.5) {};
		\node [style=none] (92) at (12.75, 1.5) {};
		\node [style=none] (94) at (11.25, 1.5) {};
		\node [style=none] (95) at (8.75, 0) {$=$};
		\node [style=none] (96) at (-6.5, 2.75) {};
		\node [style=none] (97) at (-4.5, 4.5) {};
		\node [style=none] (98) at (-4, -4.5) {};
		\node [style=none] (99) at (-4, -5) {$A$};
		\node [style=bn] (100) at (-4.75, -1.75) {};
		\node [style=bn] (101) at (-3.25, -1.75) {};
		\node [style=bn] (102) at (-4, -3.5) {};
		\node [style=medium box] (103) at (-5.5, 0.25) {$\,\;f^\dagger_p\;\,$};
		\node [style=medium box] (104) at (-2.5, 0.25) {$\,\;g^\dagger_q\;\,$};
		\node [style=none] (105) at (-3.25, -0.5) {};
		\node [style=none] (106) at (-1.75, -0.5) {};
		\node [style=none] (107) at (-6.25, -0.5) {};
		\node [style=none] (108) at (-4.75, -0.5) {};
		\node [style=none] (109) at (-4.75, -2.75) {};
		\node [style=none] (110) at (-3.25, -2.75) {};
		\node [style=morphism] (111) at (-4.75, -2.5) {$s$};
		\node [style=none] (112) at (-6.25, 0) {};
		\node [style=none] (113) at (-4.75, 0.25) {};
		\node [style=none] (114) at (-3.25, 0.25) {};
		\node [style=none] (115) at (-1.75, 0.25) {};
		\node [style=none] (116) at (-11, -4.5) {};
		\node [style=none] (117) at (-11, -5) {$A$};
		\node [style=bn] (118) at (-11.75, -0.5) {};
		\node [style=bn] (119) at (-10.25, -0.5) {};
		\node [style=morphism] (120) at (-13, 2.5) {$f$};
		\node [style=bn] (121) at (-11.5, 2.5) {};
		\node [style=none] (122) at (-13, 5) {$Z$};
		\node [style=none] (123) at (-13, 4.5) {};
		\node [style=none] (124) at (-10.5, 5) {$X$};
		\node [style=none] (125) at (-10.5, 0.75) {};
		\node [style=none] (126) at (-9.25, 5) {$Y$};
		\node [style=none] (127) at (-9.25, 0.75) {};
		\node [style=medium box] (128) at (-11, -2.5) {$\quad\rho\quad$};
		\node [style=none] (129) at (-11.75, -2.5) {};
		\node [style=none] (130) at (-10.25, -2.5) {};
		\node [style=none] (131) at (-11.5, 0.75) {};
		\node [style=none] (132) at (-13, 0.75) {};
		\node [style=none] (133) at (-10.5, 4.5) {};
		\node [style=none] (134) at (-9.25, 4.5) {};
		\node [style=none] (135) at (-5.5, 0.25) {};
		\node [style=none] (136) at (5.25, 0.25) {};
		\node [style=none] (137) at (6.75, 0.25) {};
		\node [style=morphism] (138) at (3.75, -2.5) {$s$};
		\node [style=none] (139) at (3.75, -2.75) {};
		\node [style=none] (140) at (5.25, -2.75) {};
		\node [style=none] (141) at (2.25, -0.25) {};
		\node [style=none] (142) at (4.5, 2.5) {};
		\node [style=none] (143) at (3, 2.5) {};
		\node [style=none] (144) at (1.5, 4.5) {};
		\node [style=none] (145) at (13.75, 2.25) {};
		\node [style=none] (146) at (15.25, 2.25) {};
		\node [style=morphism] (147) at (12.5, -1.25) {$s$};
		\node [style=none] (148) at (12.5, -2) {};
		\node [style=none] (149) at (14, -2) {};
		\node [style=none] (150) at (11.25, 2.25) {};
		\node [style=none] (151) at (12.75, 2.25) {};
		\node [style=none] (152) at (10, 4.5) {};
	\end{pgfonlayer}
	\begin{pgfonlayer}{edgelayer}
		\draw (7.center) to (4);
		\draw [in=270, out=15] (2) to (9.center);
		\draw [in=165, out=-90] (52.center) to (45);
		\draw [in=15, out=-90] (49.center) to (44);
		\draw (46) to (40.center);
		\draw [in=15, out=-90] (50.center) to (45);
		\draw (37.center) to (48);
		\draw [in=270, out=0] (28) to (54.center);
		\draw [in=180, out=-90] (33.center) to (28);
		\draw (35.center) to (47);
		\draw [in=165, out=-90] (92.center) to (86);
		\draw [in=15, out=-90] (90.center) to (85);
		\draw (87) to (82.center);
		\draw [in=15, out=-90] (91.center) to (86);
		\draw (80.center) to (89);
		\draw (78.center) to (88);
		\draw [in=165, out=-90] (94.center) to (85);
		\draw [in=270, out=180] (85) to (76.center);
		\draw [in=165, out=-90] (96.center) to (2);
		\draw (97.center) to (9.center);
		\draw [in=165, out=-90] (107.center) to (100);
		\draw [in=165, out=-90] (108.center) to (101);
		\draw [in=15, out=-90] (105.center) to (100);
		\draw (102) to (98.center);
		\draw [in=15, out=-90] (106.center) to (101);
		\draw (109.center) to (100);
		\draw (101) to (110.center);
		\draw [in=165, out=-90] (109.center) to (102);
		\draw [in=-90, out=15] (102) to (110.center);
		\draw (112.center) to (107.center);
		\draw (113.center) to (108.center);
		\draw (114.center) to (105.center);
		\draw (115.center) to (106.center);
		\draw [in=270, out=15] (118) to (125.center);
		\draw [in=270, out=15] (119) to (127.center);
		\draw (130.center) to (119);
		\draw (129.center) to (118);
		\draw (128) to (116.center);
		\draw (121) to (131.center);
		\draw [in=270, out=165] (119) to (131.center);
		\draw [in=165, out=-90] (132.center) to (118);
		\draw (123.center) to (132.center);
		\draw (133.center) to (125.center);
		\draw (134.center) to (127.center);
		\draw (2) to (135.center);
		\draw (11.center) to (104);
		\draw (137.center) to (50.center);
		\draw (136.center) to (49.center);
		\draw [in=165, out=-90] (139.center) to (46);
		\draw [in=15, out=-90] (140.center) to (46);
		\draw (45) to (140.center);
		\draw (44) to (139.center);
		\draw [in=165, out=-90] (141.center) to (44);
		\draw (28) to (141.center);
		\draw (142.center) to (52.center);
		\draw (143.center) to (54.center);
		\draw (144.center) to (33.center);
		\draw [in=165, out=-90] (148.center) to (87);
		\draw [in=15, out=-90] (149.center) to (87);
		\draw (149.center) to (86);
		\draw (85) to (148.center);
		\draw (150.center) to (94.center);
		\draw (151.center) to (92.center);
		\draw (145.center) to (90.center);
		\draw (146.center) to (91.center);
		\draw (152.center) to (76.center);
		\draw (96.center) to (4);
	\end{pgfonlayer}
\end{tikzpicture}}
	\end{equation*}
	using the definition of $\rho$, relative positivity of $\cat{C}$ (\Cref{rem:rel-pos}) together with the fact that $f
	\circ f^\dag$ is $(s,1_A)$-a.s.\ deterministic (\Cref{lem:f-as-det}), and the comonoid axioms.
	The same procedure can be done for the right-hand side of \eqref{eq:eq-prin-1}, and so the two sides are equal.
	Thus, by the equalizer principle, there must exist a map \( r: A \to X \times_Z Y \) in \(
	\mathcal{C} \) such that
	\begin{equation*}
		\rho = (f^*g, g^*f) \circ r ,
	\end{equation*}
	so that $(f^*g)\circ r=p$ and $(g^*f)\circ r= q$.
\end{proof}

\subsection{The universal property of hypernormalizations}\label{hypernorm}

The purpose of this section is to make the following intuition precise:
\begin{itemize}
	\item Consider a deterministic experiment $f$ on $(\Theta,p)$, which we can consider a ``partition'' or ``coarse-graining'';
	\item Form the standard measure on $P\Theta$, or ``hypernormalization'' (see \Cref{stat_exp});
	\item We can view this standard measure as the ``coarsest'' decomposition of $p$ which is still finer than the ``partition'' $f$.
\end{itemize}
From a somewhat different point of view, we can view the standard measure as \emph{the} decomposition of $p$ (in the sense of partial evaluations) induced by the partitioning $f$. 

Here is the precise statement. 

\begin{theorem}\label{thm:uni-std}
	Let $\cat{C}$ be an a.s.-compatibly representable Markov category.
	Let	$(\Theta,p)$ be a probability space, let $f:\Theta\to X$ $p$-a.s.\ deterministic, suppose its Bayesian inversion exists, and consider the standard measure $\hat{f}_p$.
	Denote also by $q$ the pushforward measure $f\circ p$ on $X$.
	
	Then for every $\pi:I \to P\Theta$ such that
	\begin{itemize}
		\item $\texttt{samp}_\Theta \circ \pi = p$ (i.e.~$\pi$ is a decomposition of $p$);
		\item $Pf\circ \pi = \delta_X \circ q$ (i.e.~each sample of $\pi$ is concentrated on a single fibre of $f$),
	\end{itemize}
	we have a partial evaluation from $\pi$ to the standard measure $\hat{f}_p$.
\end{theorem}

\begin{proof}
	Showing that there exists a partial evaluation from $\pi$ to
	$P(f^\dag_p)\circ  \delta_X\circ q$ is the same as showing that
	\begin{equation}\label{eq:unit}
		\pi \;\le\; P(f^\dag_p)\circ \delta_X\circ q
	\end{equation}
	in the order of stochastic dominance. Thus consider the experiment
	$\samp^\dag_\pi:\Theta\to P\Theta$, which, by	\Cref{stdmeas}, has $\pi$ as standard measure. Let us compose it with the
	map $Pf:P\Theta\to PX$ to obtain a new experiment $Pf\circ\samp^\dag_\pi:\Theta\to PX$, so that, taking $h:=Pf$ in the definition of the Blackwell order, $Pf\circ\samp^\dag_\pi\le
		\samp^\dag_\pi$. By Blackwell's theorem, to prove \eqref{eq:unit}
	it suffices to show that the standard measure of the composite experiment
	$Pf\circ\samp^\dag_\pi$ is the right-hand side of \eqref{eq:unit}.
	We thus want to show that
	\[
		\widehat{(\samp^\dag_{\delta_X \circ q} \circ f)}_p = P(f^\dag_p) \circ \delta_X\circ q
	\]
	which explicitly says
	\[
		\left( \left( \samp^\dag_{\delta_X \circ q} \circ f \right)^\dag_p \right)^\# \circ \samp^\dag_{\delta_X \circ q} \circ f \circ p = P(f^\dag_p) \circ \delta_X\circ q .
	\]
	We can now rewrite the left-hand side of the equation above as follows.
	\begin{align*}
		\left( \left( \samp^\dag_{\delta_X \circ q} \circ f \right)^\dag_p \right)^\# \circ \samp^\dag_{\delta_X \circ q} \circ f \circ p
		 & = \left( f^\dag_p\circ\samp_X \right)^\sharp \circ \samp^\dag_{\delta_X \circ q} \circ f \circ p \\
		 & = \left( f^\dag_p\circ\samp_X \right)^\sharp \circ \samp^\dag_{\delta_X \circ q} \circ q         \\
		 & = \left( f^\dag_p\circ\samp_X \right)^\sharp \circ \delta_X\circ q                               \\
		 & = P\left( f^\dag_p\circ\samp_X \right)\circ\delta_{PX} \circ \delta_X\circ q                     \\
		 & = P(f^\dag_p)\circ P(\samp_X) \circ\delta_{PX} \circ \delta_X\circ q                             \\
		 & = P(f^\dag_p)\circ \mu_X \circ\delta_{PX} \circ \delta_X\circ q                                  \\
		 & = P(f^\dag_p)\circ \delta_X\circ q ,
	\end{align*}
	which concludes the proof.
\end{proof}

Let $S_{p,f}$ denote the collection of states $\pi:I\to P\Theta$ such that
\[
	\samp_\Theta\circ\pi=p,
	\qquad
	Pf\circ\pi=\delta_X\circ q.
\]
By \Cref{thm:uni-std}, every $\pi\in S_{p,f}$ admits a partial evaluation to the standard measure $\hat{f}_p$, or equivalently
\[
	\pi\leq\hat{f}_p
	\qquad\text{for every }\pi\in S_{p,f}.
\]

Moreover, $\hat{f}_p$ itself belongs to $S_{p,f}$. Indeed, marginalising the Bayesian inverse condition gives
\[
	\samp_\Theta\circ\hat{f}_p=f^\dag_p\circ q=p.
\]
By \Cref{lem:f-as-det}, $f\circ f^\dag_p=_q\mathrm{id}_X$. A.s.-compatible representability and the identity $g^\#=Pg\circ\delta_X$ therefore give
\[
	P(f\circ f^\dag_p)\circ\delta_X
	=_q\delta_X.
\]
Composing with the state $q$ turns this into strict equality, as explained in \Cref{rem:stdmeas-welldef}. Using the equivalent form $\hat{f}_p=P(f^\dag_p)\circ\delta_X\circ q$ and $P$-functoriality, we obtain
\[
	Pf\circ\hat{f}_p
	=P(f\circ f^\dag_p)\circ\delta_X\circ q
	=\delta_X\circ q.
\]
Hence $\hat{f}_p$ belongs to $S_{p,f}$ and lies above every element of $S_{p,f}$. This is the maximum condition used below. In particular, \Cref{thm:uni-std} gives more than maximality. We can therefore give the following definition without assuming the existence of conditionals or Bayesian inverses.

\begin{definition}[alternative]\label{alt-def}
	Let $(\Theta,p)$ be a probability space in an a.s.-compatibly representable Markov category, let $f:\Theta\to X$ be $p$-a.s.\ deterministic, and write $q:=f\circ p$.
	The \textbf{hypernormalization} of $p$ with respect to $f$, if it exists, is a state $\pi:I\to P\Theta$ which
	\begin{itemize}
		\item satisfies $\samp_\Theta\circ\pi=p$ (i.e.~it is a decomposition of $p$);
		\item satisfies $Pf\circ\pi=\delta_X\circ q$ (i.e.~each sample of $\pi$ is concentrated on a single fibre of $f$);
		\item is a maximum with respect to stochastic dominance among the states satisfying the two conditions above, meaning that every such state $\rho$ satisfies $\rho\leq\pi$.
	\end{itemize}
\end{definition}

When a Bayesian inverse of $f$ exists, the preceding argument shows that the standard measure $\hat{f}_p$ satisfies \Cref{alt-def}. The definition itself does not require a Bayesian inverse. Without one, the existence of the maximum has to be established by other means.

The construction-based \Cref{def:hypernormalization} determines $\hat{f}_p$ canonically, since the standard measure is independent of the choice of Bayesian inverse by \Cref{rem:stdmeas-welldef}. If $\pi$ and $\pi'$ both satisfy \Cref{alt-def}, applying the defining property of each state to the other gives $\pi\leq\pi'$ and $\pi'\leq\pi$. Thus the alternative definition determines a unique state whenever stochastic dominance is antisymmetric on $S_{p,f}$, and otherwise determines the hypernormalization only up to mutual stochastic dominance. When stochastic dominance is transitive, as under the hypotheses of \Cref{cor:dominance}, mutual dominance is an equivalence relation and the alternative definition determines a unique equivalence class.

Even when the partial-evaluation relation is a preorder, it need not be antisymmetric. For monads on $\cat{Set}$ admitting indiscrete algebras --- for instance, the action monad of a group --- the partial-evaluation relation on those algebras is the kernel pair of the algebra map~\cite[Proposition~3.4.1]{Constantinetal-2020-Partialevaluationsan}. Although this is not an example of hypernormalization in a representable Markov category, it shows how distinct states can mutually dominate one another under partial evaluation.

When stochastic dominance on $S_{p,f}$ is a preorder, one could instead ask only for a maximal element, namely a state $\pi\in S_{p,f}$ such that
\[
	\pi\leq\rho\Longrightarrow\rho\leq\pi
	\qquad\text{for every }\rho\in S_{p,f}.
\]
Such a state need not lie above every compatible decomposition, and distinct maximal equivalence classes may be incomparable. We leave the existence and probabilistic meaning of maximal compatible decompositions to future work.

\section*{Acknowledgements}
Research for Mika Bohinen was part of his MSc thesis at the University of
Oxford, supported in part by Aker Scholarship, with additional support from the
Norwegian Government's L\aa{}nekasse (student loans/stipend program).

Research for Paolo Perrone was funded at the time of writing by Sam Staton's
European Research Council (ERC) Consolidator Grant ``BLaSt: Better Languages for
Statistics'' (Grant agreement ID: 864202) under the European Union's Horizon
2020 research and innovation programme.

\appendix

	\section{Some results about Markov categories}\label{appendix}
	
	\begin{lemma}\label{lem:f-as-det}
		Consider the following Bayesian inversion.
		\begin{equation}\label{eq:bayesian22}
			\begin{tikzpicture}
	\begin{pgfonlayer}{nodelayer}
		\node [style=bn] (0) at (3, 0.25) {};
		\node [style=none] (3) at (4.5, 3) {};
		\node [style=none] (4) at (4.5, 3.5) {$X$};
		\node [style=none] (5) at (2, 1.25) {};
		\node [style=none] (6) at (2, 3.5) {$Y$};
		\node [style=bn] (8) at (-3, -0.5) {};
		\node [style=morphism] (9) at (-4, 1.5) {$f$};
		\node [style=none] (10) at (-4, 3) {};
		\node [style=none] (11) at (-4, 3.5) {$Y$};
		\node [style=none] (12) at (-2, 0.5) {};
		\node [style=none] (13) at (-2, 3.5) {$X$};
		\node [style=none] (15) at (0, -0.25) {$=$};
		\node [style=morphism] (16) at (3, -2.25) {$p$};
		\node [style=none] (17) at (4, -4.5) {};
		\node [style=none] (18) at (4, -5) {$A$};
		\node [style=bn] (19) at (4, -3.75) {};
		\node [style=none] (20) at (5, -2.75) {};
		\node [style=none] (21) at (4, 1.25) {};
		\node [style=medium box] (22) at (4.5, 1.75) {$f^\dagger_p$};
		\node [style=none] (23) at (-3, -4.5) {};
		\node [style=none] (24) at (-3, -5) {$A$};
		\node [style=morphism] (25) at (-3, -2.5) {$p$};
		\node [style=none] (26) at (-2, 3) {};
		\node [style=none] (27) at (2, 3) {};
		\node [style=none] (28) at (5, 1.25) {};
		\node [style=none] (29) at (-4, 0.5) {};
		\node [style=none] (30) at (3, -2.75) {};
		\node [style=morphism] (32) at (3, -0.75) {$f$};
	\end{pgfonlayer}
	\begin{pgfonlayer}{edgelayer}
		\draw [in=-90, out=165] (0) to (5.center);
		\draw [in=-90, out=15] (8) to (12.center);
		\draw (9) to (10.center);
		\draw (16) to (0);
		\draw (19) to (17.center);
		\draw [in=15, out=-90] (20.center) to (19);
		\draw [in=15, out=-90] (21.center) to (0);
		\draw (3.center) to (22);
		\draw (25) to (23.center);
		\draw (8) to (25);
		\draw (12.center) to (26.center);
		\draw (27.center) to (5.center);
		\draw (28.center) to (20.center);
		\draw [in=165, out=-90] (29.center) to (8);
		\draw (29.center) to (9);
		\draw [in=165, out=-90] (30.center) to (19);
		\draw (16) to (30.center);
	\end{pgfonlayer}
\end{tikzpicture}
		\end{equation}
		If $f$ is almost surely deterministic, then $f\circ f^\dagger_p$ is almost surely equal to the marginalization $1_Y\otimes\texttt{del}_A$ for the measure $(f\circ p,1_A)$, i.e.\ the following equation holds.
		\begin{equation}\label{eq:f-as-det}
			\begin{tikzpicture}
	\begin{pgfonlayer}{nodelayer}
		\node [style=bn] (0) at (-5, 0) {};
		\node [style=none] (3) at (-3, 4.5) {};
		\node [style=none] (4) at (-3, 5) {$Y$};
		\node [style=none] (5) at (-6.25, 1.25) {};
		\node [style=none] (6) at (-6.25, 5) {$Y$};
		\node [style=none] (15) at (0, 0) {$=$};
		\node [style=morphism] (16) at (-5, -2.5) {$p$};
		\node [style=none] (17) at (-4.25, -4.75) {};
		\node [style=none] (18) at (-4.25, -5.25) {$A$};
		\node [style=bn] (19) at (-4.25, -3.75) {};
		\node [style=none] (20) at (-3.5, -3) {};
		\node [style=none] (21) at (-3.75, 1.25) {};
		\node [style=medium box] (22) at (-3, 1.75) {$\;f^\dagger_p\;$};
		\node [style=none] (27) at (-6.25, 4.5) {};
		\node [style=none] (30) at (-5, -3) {};
		\node [style=morphism] (32) at (-5, -1) {$f$};
		\node [style=morphism] (33) at (-3, 3.5) {$f$};
		\node [style=bn] (34) at (3, 0) {};
		\node [style=none] (35) at (4.25, 4.5) {};
		\node [style=none] (36) at (4.25, 5) {$Y$};
		\node [style=none] (37) at (1.75, 1.25) {};
		\node [style=none] (38) at (1.75, 5) {$Y$};
		\node [style=morphism] (39) at (3, -2.5) {$p$};
		\node [style=none] (40) at (3.75, -4.75) {};
		\node [style=none] (41) at (3.75, -5.25) {$A$};
		\node [style=bn] (42) at (3.75, -3.75) {};
		\node [style=none] (43) at (4.5, -3) {};
		\node [style=none] (44) at (4.25, 1.25) {};
		\node [style=none] (46) at (1.75, 4.5) {};
		\node [style=none] (48) at (3, -3) {};
		\node [style=morphism] (49) at (3, -1) {$f$};
		\node [style=bn] (50) at (4.5, 0) {};
		\node [style=none] (51) at (3.25, 1.25) {};
		\node [style=none] (52) at (5.75, 1.25) {};
		\node [style=none] (53) at (3.25, 4.5) {};
		\node [style=bn] (54) at (5.75, 2.75) {};
		\node [style=none] (55) at (3.25, 5) {$A$};
		\node [style=bn] (56) at (-3.5, 0) {};
		\node [style=none] (57) at (-4.75, 1.25) {};
		\node [style=none] (58) at (-4.75, 4.5) {};
		\node [style=none] (59) at (-4.75, 5) {$A$};
		\node [style=none] (60) at (-2.25, 1.25) {};
		\node [style=none] (61) at (-3.75, 1.75) {};
		\node [style=none] (62) at (-2.25, 1.75) {};
	\end{pgfonlayer}
	\begin{pgfonlayer}{edgelayer}
		\draw [in=-90, out=165] (0) to (5.center);
		\draw (16) to (0);
		\draw (19) to (17.center);
		\draw [in=15, out=-90] (20.center) to (19);
		\draw [in=15, out=-90] (21.center) to (0);
		\draw (3.center) to (22);
		\draw (27.center) to (5.center);
		\draw [in=165, out=-90] (30.center) to (19);
		\draw (16) to (30.center);
		\draw [in=-90, out=165] (34) to (37.center);
		\draw (39) to (34);
		\draw (42) to (40.center);
		\draw [in=15, out=-90] (43.center) to (42);
		\draw [in=15, out=-90] (44.center) to (34);
		\draw (46.center) to (37.center);
		\draw [in=165, out=-90] (48.center) to (42);
		\draw (39) to (48.center);
		\draw (35.center) to (44.center);
		\draw (50) to (43.center);
		\draw [in=165, out=-90] (51.center) to (50);
		\draw [in=270, out=15] (50) to (52.center);
		\draw (54) to (52.center);
		\draw (53.center) to (51.center);
		\draw (58.center) to (57.center);
		\draw [in=165, out=-90] (57.center) to (56);
		\draw [in=15, out=-90] (60.center) to (56);
		\draw (56) to (20.center);
		\draw (61.center) to (21.center);
		\draw (62.center) to (60.center);
	\end{pgfonlayer}
\end{tikzpicture}
		\end{equation}
		In particular, $f\circ f^\dag_p$ is almost surely deterministic.
		Note that the equation~\eqref{eq:f-as-det} itself, not just the
		determinism conclusion, is used in the proof of \Cref{thm:main}.
	\end{lemma}
	\begin{proof}
		We have
		\[
			\resizebox{\linewidth}{!}{\begin{tikzpicture}
	\begin{pgfonlayer}{nodelayer}
		\node [style=bn] (0) at (-12.75, 0) {};
		\node [style=none] (3) at (-10.75, 5) {};
		\node [style=none] (4) at (-10.75, 5.5) {$Y$};
		\node [style=none] (5) at (-14, 1.25) {};
		\node [style=none] (6) at (-14, 5.5) {$Y$};
		\node [style=none] (15) at (-9, 0) {$=$};
		\node [style=morphism] (16) at (-12.75, -2.75) {$p$};
		\node [style=none] (17) at (-12, -5.25) {};
		\node [style=none] (18) at (-12, -5.75) {$A$};
		\node [style=bn] (19) at (-12, -4.25) {};
		\node [style=none] (20) at (-11.25, -3.5) {};
		\node [style=none] (21) at (-11.5, 1.25) {};
		\node [style=medium box] (22) at (-10.75, 2.5) {$\;f^\dagger_p\;$};
		\node [style=none] (27) at (-14, 5) {};
		\node [style=none] (30) at (-12.75, -3.5) {};
		\node [style=morphism] (32) at (-12.75, -1.125) {$f$};
		\node [style=morphism] (33) at (-10.75, 4) {$f$};
		\node [style=bn] (56) at (-11.25, 0) {};
		\node [style=none] (57) at (-12.5, 1.25) {};
		\node [style=none] (58) at (-12.5, 5) {};
		\node [style=none] (59) at (-12.5, 5.5) {$A$};
		\node [style=none] (60) at (-10, 1.25) {};
		\node [style=none] (61) at (-11.5, 2.5) {};
		\node [style=none] (62) at (-10, 2.5) {};
		\node [style=bn] (63) at (-6.75, 0.5) {};
		\node [style=none] (64) at (-4.75, 5) {};
		\node [style=none] (65) at (-4.75, 5.5) {$Y$};
		\node [style=none] (66) at (-8, 1.75) {};
		\node [style=none] (67) at (-8, 5.5) {$Y$};
		\node [style=morphism] (68) at (-6.75, -2) {$p$};
		\node [style=none] (69) at (-5, -5.25) {};
		\node [style=none] (70) at (-5, -5.75) {$A$};
		\node [style=bn] (71) at (-6, -3.25) {};
		\node [style=none] (72) at (-5.25, -2.5) {};
		\node [style=none] (73) at (-5.5, 1.75) {};
		\node [style=medium box] (74) at (-4.75, 2.5) {$\;f^\dagger_p\;$};
		\node [style=none] (75) at (-8, 5) {};
		\node [style=none] (76) at (-6.75, -2.5) {};
		\node [style=morphism] (77) at (-6.75, -0.5) {$f$};
		\node [style=morphism] (78) at (-4.75, 4) {$f$};
		\node [style=none] (80) at (-6.5, 1.75) {};
		\node [style=none] (81) at (-6.5, 5) {};
		\node [style=none] (82) at (-6.5, 5.5) {$A$};
		\node [style=none] (83) at (-4, 1.75) {};
		\node [style=none] (84) at (-5.5, 2.5) {};
		\node [style=none] (85) at (-4, 2.5) {};
		\node [style=bn] (86) at (-5, -4.25) {};
		\node [style=none] (87) at (-4, -3.25) {};
		\node [style=none] (88) at (-4, -1.5) {};
		\node [style=none] (89) at (-5.25, -1.5) {};
		\node [style=none] (90) at (-3, 0) {$=$};
		\node [style=bn] (91) at (-0.75, -0.25) {};
		\node [style=none] (92) at (1.25, 5) {};
		\node [style=none] (93) at (1.25, 5.5) {$Y$};
		\node [style=none] (94) at (-2, 1) {};
		\node [style=none] (95) at (-2, 5.5) {$Y$};
		\node [style=morphism] (96) at (-0.75, -2) {$p$};
		\node [style=none] (97) at (0.25, -5.25) {};
		\node [style=none] (98) at (0.25, -5.75) {$A$};
		\node [style=none] (101) at (1.25, 1.75) {};
		\node [style=none] (103) at (-2, 5) {};
		\node [style=none] (104) at (-0.75, -2.75) {};
		\node [style=morphism] (106) at (1.25, 3) {$f$};
		\node [style=none] (107) at (-0.5, 1.75) {};
		\node [style=none] (108) at (-0.5, 5) {};
		\node [style=none] (109) at (-0.5, 5.5) {$A$};
		\node [style=bn] (113) at (0.25, -3.75) {};
		\node [style=none] (114) at (1.25, -2.75) {};
		\node [style=none] (115) at (1.25, -1.5) {};
		\node [style=bn] (118) at (11, 0.5) {};
		\node [style=none] (119) at (12.25, 5) {};
		\node [style=none] (120) at (12.25, 5.5) {$Y$};
		\node [style=none] (121) at (9.75, 1.75) {};
		\node [style=none] (122) at (9.75, 5.5) {$Y$};
		\node [style=morphism] (123) at (11, -2.25) {$p$};
		\node [style=none] (124) at (11.75, -5.25) {};
		\node [style=none] (125) at (11.75, -5.75) {$A$};
		\node [style=bn] (126) at (11.75, -3.75) {};
		\node [style=none] (127) at (12.5, -3) {};
		\node [style=none] (128) at (12.25, 1.75) {};
		\node [style=none] (129) at (9.75, 5) {};
		\node [style=none] (130) at (11, -3) {};
		\node [style=morphism] (131) at (11, -0.625) {$f$};
		\node [style=bn] (132) at (12.5, 0.5) {};
		\node [style=none] (133) at (11.25, 1.75) {};
		\node [style=none] (134) at (13.75, 1.75) {};
		\node [style=none] (135) at (11.25, 5) {};
		\node [style=bn] (136) at (13.75, 3.25) {};
		\node [style=none] (137) at (11.25, 5.5) {$A$};
		\node [style=none] (138) at (8.5, 0) {$=$};
		\node [style=none] (139) at (2.75, 0) {$=$};
		\node [style=morphism] (166) at (-2, 3) {$f$};
		\node [style=bn] (167) at (5.25, 0.5) {};
		\node [style=none] (168) at (7.25, 5) {};
		\node [style=none] (169) at (7.25, 5.5) {$Y$};
		\node [style=none] (170) at (4, 1.75) {};
		\node [style=none] (171) at (4, 5.5) {$Y$};
		\node [style=morphism] (172) at (5.25, -2.25) {$p$};
		\node [style=none] (173) at (6.25, -5.25) {};
		\node [style=none] (174) at (6.25, -5.75) {$A$};
		\node [style=none] (175) at (7.25, 2.5) {};
		\node [style=none] (176) at (4, 5) {};
		\node [style=none] (177) at (5.25, -2.75) {};
		\node [style=none] (179) at (5.5, 1.75) {};
		\node [style=none] (180) at (5.5, 5) {};
		\node [style=none] (181) at (5.5, 5.5) {$A$};
		\node [style=bn] (182) at (6.25, -3.75) {};
		\node [style=none] (183) at (7.25, -2.75) {};
		\node [style=none] (184) at (7.25, -1.5) {};
		\node [style=morphism] (185) at (5.25, -0.625) {$f$};
	\end{pgfonlayer}
	\begin{pgfonlayer}{edgelayer}
		\draw [in=-90, out=165] (0) to (5.center);
		\draw (16) to (0);
		\draw (19) to (17.center);
		\draw [in=15, out=-90] (20.center) to (19);
		\draw [in=15, out=-90] (21.center) to (0);
		\draw (3.center) to (22);
		\draw (27.center) to (5.center);
		\draw [in=165, out=-90] (30.center) to (19);
		\draw (16) to (30.center);
		\draw (58.center) to (57.center);
		\draw [in=165, out=-90] (57.center) to (56);
		\draw [in=15, out=-90] (60.center) to (56);
		\draw (56) to (20.center);
		\draw (61.center) to (21.center);
		\draw (62.center) to (60.center);
		\draw [in=-90, out=165] (63) to (66.center);
		\draw (68) to (63);
		\draw [in=15, out=-90] (72.center) to (71);
		\draw [in=15, out=-90] (73.center) to (63);
		\draw (64.center) to (74);
		\draw (75.center) to (66.center);
		\draw [in=165, out=-90] (76.center) to (71);
		\draw (68) to (76.center);
		\draw (81.center) to (80.center);
		\draw (84.center) to (73.center);
		\draw (85.center) to (83.center);
		\draw [in=165, out=-90] (71) to (86);
		\draw (86) to (69.center);
		\draw [in=15, out=-90] (87.center) to (86);
		\draw (88.center) to (87.center);
		\draw [in=90, out=-90] (80.center) to (88.center);
		\draw (89.center) to (72.center);
		\draw [in=90, out=-90] (83.center) to (89.center);
		\draw [in=-90, out=165] (91) to (94.center);
		\draw (96) to (91);
		\draw [in=15, out=-90] (101.center) to (91);
		\draw (103.center) to (94.center);
		\draw (96) to (104.center);
		\draw (108.center) to (107.center);
		\draw (113) to (97.center);
		\draw [in=15, out=-90] (114.center) to (113);
		\draw (115.center) to (114.center);
		\draw [in=90, out=-90] (107.center) to (115.center);
		\draw [in=-90, out=165] (118) to (121.center);
		\draw (123) to (118);
		\draw (126) to (124.center);
		\draw [in=15, out=-90] (127.center) to (126);
		\draw [in=15, out=-90] (128.center) to (118);
		\draw (129.center) to (121.center);
		\draw [in=165, out=-90] (130.center) to (126);
		\draw (123) to (130.center);
		\draw (119.center) to (128.center);
		\draw (132) to (127.center);
		\draw [in=165, out=-90] (133.center) to (132);
		\draw [in=270, out=15] (132) to (134.center);
		\draw (136) to (134.center);
		\draw (135.center) to (133.center);
		\draw [in=165, out=-90] (104.center) to (113);
		\draw (92.center) to (101.center);
		\draw [in=-90, out=165] (167) to (170.center);
		\draw (172) to (167);
		\draw [in=15, out=-90] (175.center) to (167);
		\draw (176.center) to (170.center);
		\draw (172) to (177.center);
		\draw (180.center) to (179.center);
		\draw (182) to (173.center);
		\draw [in=15, out=-90] (183.center) to (182);
		\draw (184.center) to (183.center);
		\draw [in=90, out=-90] (179.center) to (184.center);
		\draw [in=165, out=-90] (177.center) to (182);
		\draw (168.center) to (175.center);
	\end{pgfonlayer}
\end{tikzpicture}}
		\]
		using the comonoid axioms, \eqref{eq:bayesian22}, almost sure determinism of $f$, and explicitly copying and discarding.
	\end{proof}
	
	\begin{lemma}\label{lem:delta-q-as-det}
		Let $\cat{C}$ be an a.s.-compatibly representable Markov category with
		conditionals. If  $p \in \cat{C}(I, X)$,
		then $\texttt{samp}_X$ is $(\delta_X \circ p)$-a.s.
		deterministic.
	\end{lemma}
	\begin{proof}
		We want to show that the following equality holds
		\begin{equation}\label{eq:delta-q-as-det1}
			\begin{tikzpicture}
	\begin{pgfonlayer}{nodelayer}
		\node [style=state] (0) at (-3, -3.25) {$p$};
		\node [style=bn] (1) at (-3, 0) {};
		\node [style=bn] (2) at (-1.5, 2.75) {};
		\node [style=none] (3) at (-2.5, 3.75) {};
		\node [style=none] (4) at (-2.5, 4.25) {$X$};
		\node [style=none] (5) at (-0.5, 3.75) {};
		\node [style=none] (6) at (-0.5, 4.25) {$X$};
		\node [style=none] (7) at (-4.5, 3.75) {};
		\node [style=none] (8) at (-4.5, 4.25) {$PX$};
		\node [style=none] (9) at (-4.5, 1.5) {};
		\node [style=morphism] (10) at (-3, -1.5) {$\delta_X$};
		\node [style=morphism] (11) at (-1.5, 1.5) {$\texttt{samp}_X$};
		\node [style=none] (12) at (0, 0) {$=$};
		\node [style=state] (13) at (3, -3.25) {$p$};
		\node [style=bn] (14) at (3, 0) {};
		\node [style=bn] (15) at (5, 1.5) {};
		\node [style=none] (16) at (3.5, 3.75) {};
		\node [style=none] (17) at (3.5, 4.25) {$X$};
		\node [style=none] (18) at (6.5, 3.75) {};
		\node [style=none] (19) at (6.5, 4.25) {$X$};
		\node [style=none] (20) at (1.5, 3.75) {};
		\node [style=none] (21) at (1.5, 4.25) {$PX$};
		\node [style=none] (22) at (1.5, 1.5) {};
		\node [style=morphism] (23) at (3, -1.5) {$\delta_X$};
		\node [style=morphism] (24) at (3.5, 2.75) {$\texttt{samp}_X$};
		\node [style=morphism] (25) at (6.5, 2.75) {$\texttt{samp}_X$};
	\end{pgfonlayer}
	\begin{pgfonlayer}{edgelayer}
		\draw (9.center) to (7.center);
		\draw [in=270, out=180] (2) to (3.center);
		\draw [in=270, out=0] (2) to (5.center);
		\draw [in=180, out=-90, looseness=1.25] (9.center) to (1);
		\draw (1) to (10);
		\draw (10) to (0);
		\draw [in=0, out=-90, looseness=1.25] (11) to (1);
		\draw (2) to (11);
		\draw (22.center) to (20.center);
		\draw [in=180, out=-90, looseness=1.25] (22.center) to (14);
		\draw (14) to (23);
		\draw (23) to (13);
		\draw [in=180, out=-90] (24) to (15);
		\draw [in=270, out=0] (15) to (25);
		\draw [in=0, out=-90] (15) to (14);
		\draw (24) to (16.center);
		\draw (25) to (18.center);
	\end{pgfonlayer}
\end{tikzpicture}
		\end{equation}
		Using that $\delta_X$ is deterministic we get that the left hand side
		of \Cref{eq:delta-q-as-det1} can be written as
		\begin{equation}\label{eq:delta-q-as-det2}
			\begin{tikzpicture}
	\begin{pgfonlayer}{nodelayer}
		\node [style=state] (0) at (-6, -3.75) {$p$};
		\node [style=bn] (1) at (-6, 0) {};
		\node [style=bn] (2) at (-4.5, 2.75) {};
		\node [style=none] (3) at (-5.5, 3.75) {};
		\node [style=none] (4) at (-5.5, 4.25) {$X$};
		\node [style=none] (5) at (-3.5, 3.75) {};
		\node [style=none] (6) at (-3.5, 4.25) {$X$};
		\node [style=none] (7) at (-7.5, 3.75) {};
		\node [style=none] (8) at (-7.5, 4.25) {$PX$};
		\node [style=none] (9) at (-7.5, 1.5) {};
		\node [style=morphism] (10) at (-6, -1.75) {$\delta_X$};
		\node [style=morphism] (11) at (-4.5, 1.5) {$\texttt{samp}_X$};
		\node [style=none] (12) at (-3, 0) {$=$};
		\node [style=state] (13) at (0.25, -3.75) {$p$};
		\node [style=bn] (14) at (0.25, -2.25) {};
		\node [style=bn] (15) at (1.75, 2.75) {};
		\node [style=none] (16) at (0.75, 3.75) {};
		\node [style=none] (17) at (0.75, 4.25) {$X$};
		\node [style=none] (18) at (2.75, 3.75) {};
		\node [style=none] (19) at (2.75, 4.25) {$X$};
		\node [style=none] (20) at (-1.25, 3.75) {};
		\node [style=none] (21) at (-1.25, 4.25) {$PX$};
		\node [style=none] (22) at (-1.25, 1.5) {};
		\node [style=morphism] (23) at (-1.25, -0.5) {$\delta_X$};
		\node [style=morphism] (24) at (1.75, 1.5) {$\texttt{samp}_X$};
		\node [style=morphism] (25) at (1.75, -0.5) {$\delta_X$};
		\node [style=none] (26) at (4.25, 0) {$=$};
		\node [style=state] (27) at (7.5, -3.75) {$p$};
		\node [style=bn] (28) at (7.5, -2.25) {};
		\node [style=bn] (29) at (9, 2.75) {};
		\node [style=none] (30) at (8, 3.75) {};
		\node [style=none] (31) at (8, 4.25) {$X$};
		\node [style=none] (32) at (10, 3.75) {};
		\node [style=none] (33) at (10, 4.25) {$X$};
		\node [style=none] (34) at (6, 3.75) {};
		\node [style=none] (35) at (6, 4.25) {$PX$};
		\node [style=none] (36) at (6, 1.5) {};
		\node [style=morphism] (37) at (6, -0.5) {$\delta_X$};
		\node [style=none] (38) at (9, -0.5) {};
	\end{pgfonlayer}
	\begin{pgfonlayer}{edgelayer}
		\draw (9.center) to (7.center);
		\draw [in=270, out=180] (2) to (3.center);
		\draw [in=270, out=0] (2) to (5.center);
		\draw [in=180, out=-90, looseness=1.25] (9.center) to (1);
		\draw (1) to (10);
		\draw (10) to (0);
		\draw [in=0, out=-90, looseness=1.25] (11) to (1);
		\draw (2) to (11);
		\draw (22.center) to (20.center);
		\draw [in=270, out=180] (15) to (16.center);
		\draw [in=270, out=0] (15) to (18.center);
		\draw (15) to (24);
		\draw (23) to (22.center);
		\draw (24) to (25);
		\draw [in=270, out=180] (14) to (23);
		\draw [in=270, out=0] (14) to (25);
		\draw (14) to (13);
		\draw (36.center) to (34.center);
		\draw [in=270, out=180] (29) to (30.center);
		\draw [in=270, out=0] (29) to (32.center);
		\draw (37) to (36.center);
		\draw [in=270, out=180] (28) to (37);
		\draw (28) to (27);
		\draw (29) to (38.center);
		\draw [in=0, out=-90] (38.center) to (28);
	\end{pgfonlayer}
\end{tikzpicture}
		\end{equation}
		Similarly, the right hand side of \Cref{eq:delta-q-as-det1} can be
		written as
		\begin{equation}\label{eq:delta-q-as-det3}
			\begin{tikzpicture}
	\begin{pgfonlayer}{nodelayer}
		\node [style=state] (0) at (-4.5, -3.25) {$p$};
		\node [style=bn] (1) at (-4.5, 0) {};
		\node [style=bn] (2) at (-2.75, 1.5) {};
		\node [style=none] (3) at (-4.25, 3.75) {};
		\node [style=none] (4) at (-4.25, 4.25) {$X$};
		\node [style=none] (5) at (-1.25, 3.75) {};
		\node [style=none] (6) at (-1.25, 4.25) {$X$};
		\node [style=none] (7) at (-6.25, 3.75) {};
		\node [style=none] (8) at (-6.25, 4.25) {$PX$};
		\node [style=none] (9) at (-6.25, 1.5) {};
		\node [style=morphism] (10) at (-4.5, -1.5) {$\delta_X$};
		\node [style=morphism] (11) at (-4.25, 2.75) {$\texttt{samp}_X$};
		\node [style=morphism] (12) at (-1.25, 2.75) {$\texttt{samp}_X$};
		\node [style=none] (13) at (0, -0.5) {$=$};
		\node [style=state] (14) at (3.5, -3.25) {$p$};
		\node [style=bn] (15) at (3.5, -1.5) {};
		\node [style=bn] (16) at (5.25, 1.5) {};
		\node [style=none] (17) at (3.75, 3.75) {};
		\node [style=none] (18) at (3.75, 4.25) {$X$};
		\node [style=none] (19) at (6.75, 3.75) {};
		\node [style=none] (20) at (6.75, 4.25) {$X$};
		\node [style=none] (21) at (1.75, 3.75) {};
		\node [style=none] (22) at (1.75, 4.25) {$PX$};
		\node [style=none] (23) at (1.75, 1.5) {};
		\node [style=morphism] (24) at (1.75, 0) {$\delta_X$};
		\node [style=morphism] (25) at (3.75, 2.75) {$\texttt{samp}_X$};
		\node [style=morphism] (26) at (6.75, 2.75) {$\texttt{samp}_X$};
		\node [style=none] (27) at (8, -0.5) {$=$};
		\node [style=morphism] (28) at (5.25, 0) {$\delta_X$};
		\node [style=state] (29) at (11.75, -3.25) {$p$};
		\node [style=bn] (30) at (11.75, -1.5) {};
		\node [style=bn] (31) at (13.5, 1.5) {};
		\node [style=none] (32) at (12, 3.75) {};
		\node [style=none] (33) at (12, 4.25) {$X$};
		\node [style=none] (34) at (15, 3.75) {};
		\node [style=none] (35) at (15, 4.25) {$X$};
		\node [style=none] (36) at (10, 3.75) {};
		\node [style=none] (37) at (10, 4.25) {$PX$};
		\node [style=none] (38) at (10, 1.5) {};
		\node [style=morphism] (39) at (10, 0) {$\delta_X$};
		\node [style=none] (40) at (13.5, 0) {};
	\end{pgfonlayer}
	\begin{pgfonlayer}{edgelayer}
		\draw (9.center) to (7.center);
		\draw [in=180, out=-90, looseness=1.25] (9.center) to (1);
		\draw (1) to (10);
		\draw (10) to (0);
		\draw [in=180, out=-90] (11) to (2);
		\draw [in=270, out=0] (2) to (12);
		\draw [in=0, out=-90] (2) to (1);
		\draw (11) to (3.center);
		\draw (12) to (5.center);
		\draw (23.center) to (21.center);
		\draw [in=180, out=-90] (25) to (16);
		\draw [in=270, out=0] (16) to (26);
		\draw (25) to (17.center);
		\draw (26) to (19.center);
		\draw (23.center) to (24);
		\draw [in=180, out=-90] (24) to (15);
		\draw [in=270, out=0] (15) to (28);
		\draw (28) to (16);
		\draw (15) to (14);
		\draw (38.center) to (36.center);
		\draw (38.center) to (39);
		\draw [in=180, out=-90] (39) to (30);
		\draw (30) to (29);
		\draw [in=0, out=-90] (40.center) to (30);
		\draw (40.center) to (31);
		\draw [in=270, out=180] (31) to (32.center);
		\draw [in=270, out=0] (31) to (34.center);
	\end{pgfonlayer}
\end{tikzpicture}
		\end{equation}
		Hence \Cref{eq:delta-q-as-det1} holds and we are done.
	\end{proof}
	
	\begin{corollary}\label{cor:samp-det}
		Let $\cat{C}$ be an a.s.-compatibly representable Markov category with
		conditionals. Then, for all $p \in \cat{C}(I, X)$ we have
		\begin{equation}
			\texttt{samp} \circ \texttt{samp}_{\delta \circ p}^\dag =_p 1_X.
		\end{equation}
	\end{corollary}

\bibliographystyle{alpha}
\bibliography{./refs-paper}

\end{document}